\documentclass[11pt]{article}
\usepackage{amsfonts}
\usepackage{amssymb}
\def\c{{ c}}
\def\d{{ d}}
\def\k{{ k}}
\def\L{{ L}}
\def\RR{{\bf R}}
\def\BB{{\bf B}}
\def\CC{{\bf C}}
\def\DD{{\bf D}}

\def\HH{{\bf H}}
\def\PP{{\bf P}}
\def\Chat{\hat{\bf C}}
\def\ZZ{{\bf Z}}

\def\tr{\mathop{\rm Tr}}
\def\mmod{\mathop{\rm mod}}
\def\ax{\mathop{ Ax}}
\def\teich{\mathop{\rm  Teich}}
\def\torus{{\cal T}_1}

\def\QF{{\cal Q}{\cal F}}

\def\bch{{\partial{\cal C}}}

\newtheorem{thm}{Theorem}[section]
\newtheorem{thmI}{Theorem}
\newtheorem{lemma}[thm]{Lemma}

\newtheorem{cor}[thm]{Corollary}
\newtheorem{prop}[thm]{Proposition}
\newenvironment{genericem}[1]{\smallskip\noindent{\bf
#1}\em}{\rm
\smallskip}
\newtheorem{remark}[thm]{Remark}
\newenvironment{proof}{{\sc Proof.}}{$\;\square$ \vskip .2in}

\def\H#1{\mathop{{\bf H}^#1}}

\def\Pr#1{\mathop{{\cal P}_{#1}}}

\def\F{{\cal F}}
\def\E{{\cal E}}
\def\Q{{\cal Q}}
\def\P{{\cal P}}
\def\C{{\cal C}}

\def\S{{\cal S}}
\begin{document}

\title{ Pleating invariants for  punctured torus groups }
\author{Linda Keen
\thanks{Research
partially supported by NSF, PSC-CUNY, EPSRC}
\\Mathematics Department\\CUNY Lehman College\\Bronx,
NY 10468, U.S.A. \and Caroline
Series\\Mathematics Institute\\ Warwick University\\Coventry CV4
7AL,
U.K.}

\maketitle
\bibliographystyle{plain}

\begin{abstract}
In this paper we give a complete description of the space $\QF$
of
quasifuchsian punctured torus groups in terms of what we call
{\em pleating
invariants}.  These are natural invariants of the boundary
$\bch$ of the
convex core of the associated hyperbolic $3$-manifold $M$ and
give
coordinates for the non-Fuchsian groups $\QF - \F$.  The
pleating
invariants of a component of $\bch$ consist of the projective
class of its
bending measure, together with the lamination length of a fixed
choice of
transverse measure in this class.  Our description complements
that of
Minsky in \cite{MinskyPT}, in which he describes the space of
all punctured
torus groups in terms of {\em ending invariants} which
characterize the
asymptotic geometry of the ends of $M$.

Pleating invariants give a quasifuchsian analog of the
Kerckhoff-Thurston
description of Fuchsian space by critical lines and earthquake
horocycles.
The critical lines extend to {\em pleating planes} on which the
pleating
loci of $\bch$ are constant and the horocycles extend to {\em
BM-slices} on
which the pleating invariants of one component of $\bch$ are
fixed.

We prove that the pleating planes corresponding to rational
laminations are
dense and that their boundaries can be found {\em explicitly}.
This means,
answering questions posed by Bers in the late 1960's, that it is
possible
to compute an arbitrarily accurate picture of the shape of any
embedding of
$\QF$ into $\CC^2$.
\end{abstract}

\section{Introduction}
\label{sec:introduction}

In his recent paper \cite{MinskyPT}, Minsky gave a full
description of the
space of punctured torus groups in terms of their {\em ending
invariants}.
These invariants are the conformal structures of the quotient
surfaces of
the regular set of the group acting on the Riemann sphere, or,
if a
component is degenerate, the corresponding {\em ending
lamination} of
Thurston.

 In this paper we give an alternative description of
quasifuchsian space
$\QF$ in terms of what we call {\em pleating invariants}.  These
replace
conformal structures at infinity by natural invariants of the
geometry of
the boundary of the convex core of the associated three
manifold.  These
invariants again extend naturally to ending laminations for
groups on the
boundary of $\QF$.  Pleating invariants have considerable
computational
advantages: we show how they can be used to explicitly locate
the group
with given invariants, and to compute the shape and boundary of
$\QF$, for
any embedding into $\CC^2$.

\smallskip

A {\em punctured torus group} $<G; A,B> $ is a free marked two
generator
discrete subgroup of $PSL(2,\CC)$ such that the commutator of
the
generators is parabolic. Such a group is the image of a faithful
representation $\rho$ of the fundamental group of a punctured
torus
$\torus$ with presentation $ \pi_1(\torus)=<\alpha,\beta>$; the
commutator
of the generators represents a loop around the puncture and the
ordered
pair $(A ,B) = (\rho(\alpha),\rho(\beta))$ is the marking.  The
group $G$
acts as a discrete group of isometries of hyperbolic space $\H3$
and the
quotient hyperbolic manifold $M=\H3/G$ is a product $\torus
\times (-1,1)$.

A punctured torus group also acts as a group of conformal
automorphisms of
the Riemann sphere $\Chat$ and partitions it into two invariant
subsets,
the open (possibly empty) regular set $\Omega$ and the closed
limit set
$\Lambda$.  The group $G$ is {\em quasifuchsian} if $\Omega$
consists of
two non-empty simply connected invariant components denoted
$\Omega^{\pm}$.
The quotients $\Omega^{\pm}/G$ are punctured tori with conformal
structures
inherited from $\Chat$.

{\em Quasifuchsian space} $\QF$ is the space of quasifuchsian
marked
punctured torus groups modulo conjugation in $PSL(2,\CC)$; {\em
Fuchsian
space} $\F$ is the subset such that the components
$\Omega^{\pm}$ are round
disks.

The convex hull $\C$ of $\Lambda$ in $\H3$ is also invariant
under $G$.
The hyperbolic manifold $\C/G$ is called the {\em convex core}
of $G$.  If
$G$ is quasifuchsian, but not Fuchsian, $\partial{\C/G}$
consists of two
components, $\bch^{\pm}/G$.  Each component is homeomorphic to
$\torus$ and
admits an intrinsic hyperbolic structure making it a {\em
pleated surface}
in the sense of Thurston.  Such a surface is a hyperbolic
surface ``bent''
along a geodesic lamination called the {\em pleating locus} or
{\em bending
lamination}. The pleating locus carries a natural transverse
measure, the
bending measure $pl^{\pm}(G)$.

For any measured geodesic lamination $\mu$ on a hyperbolic
surface
$\Sigma$, we denote the projective class of $\mu$ by $[\mu]$ and
the
underlying lamination by $|\mu|$. Writing $l_{\mu}$ for the
lamination
length of $\mu$, we note that if $\mu, \mu'$ are in the same
projective
class, so that $\mu=c\mu', c>0$, then $l_{\mu}= cl_{\mu'}$.  We
define the
{\em pleating invariants} for $G \in \QF-\F$ to be the
projective class of
the pair $(\mu^{\pm}, l_{\mu^{\pm}})$ for any choice of measured
laminations $\mu^{\pm}$ in $[pl^{\pm}]$.

We prove
\begin{thmI}
\label{thm:plinvarsI} A non-Fuchsian quasifuchsian marked
punctured torus group is determined by its pleating invariants,
uniquely up
to conjugacy in $PSL(2,\CC)$.
\end{thmI}

The essential idea is to study the sets in $\QF$ on which some
or all of
the pleating invariants are constant; in particular, we study
the set
$\Pr{\mu,\nu} \subset \QF$ for which $[pl^+]=[\mu],
[pl^-]=[\nu]$.  Clearly
$\Pr{\mu,\nu}$ depends only on the projective classes $[\mu],
[\nu]$ of
$\mu,\nu$.  We prove that these sets are connected real two
dimensional
submanifolds of $\QF$ whose boundaries meet $\F$ and $\partial
\QF$ in
specific analytic curves; as the projective classes vary, the
sets
$\Pr{\mu,\nu}$, which for obvious reasons we call pleating
planes, foliate
$\QF-\F$.  We are also able to describe exactly how the closure
of
$\Pr{\mu,\nu}$ meets $\F$.

\medskip

The space $\QF$ has a natural $\CC^2$-holomorphic structure
induced from
$PSL(2,\CC)$.  Let $U \subset \QF$.  An {\em $\RR^2$-locus} in
$U$ is a set
$ f^{-1}(\RR^2) \cap U$ where $f:U \to \CC^2$ is a non-constant
holomorphic
function defined on $U$.  A {\em singularity} is a point where
$Det(Jac_f(z))=0$.  For example, Fuchsian space is an
$\RR^2$-locus in
$\QF$, (see section~\ref{sec:cfncoords}).

The starting point for our analysis of $\Pr{\mu,\nu}$ is to
prove that for
 $\mu \in ML$, the length function $l_{\mu}$ on $\F$ extends to
a
 holomorphic function $\lambda_{\mu}$, called the {\em complex
length} of
 $\mu$, on $\QF$, and that $\lambda_{\mu} $ is real valued at
points where
 the projective class of $pl^{\pm}$ is $[\mu]$.  Thus
$\Pr{\mu,\nu}$ is
 contained in the $\RR^2$-locus of the holomorphic function
$\L_{\mu,\nu} =
 \lambda_{\mu} \times \lambda_{\nu}$ from $\QF $ to $ \CC^2$.

To describe $\Pr{\mu,\nu}$ more precisely, we recall some facts
about
Fuchsian space $\F$.  Let $\mu$ be a measured geodesic
lamination on a
hyperbolic surface $\Sigma$.  The distance $t$ earthquake
$\E_{\mu}(t)$
along $\mu$ gives a one parameter family of deformations of $\F$
which
generalize Fenchel-Nielsen twists along simple closed geodesics.
 For a
point $p \in \F$, we denote the earthquake path $\{ \E_{\mu}(t)
(p): t \in
\RR \}$ through $p$ by $\E_{\mu}^p$.  The earthquake path is
contained in
$\F$ and meets $\partial\F$, the Thurston boundary of $\F$, in
the point
$[\mu]$.  Kerckhoff proved that for each measured lamination
$\nu$ whose
intersection $i(\mu,\nu)$ with $\mu$ is non-zero, the length
function
$l_{\nu}$ has a unique minimum along $\E_{\mu}^p$.

In the special case of the punctured torus, it is an easy
consequence of
Kerckhoff's results that for each $\c>0$, there is a unique
earthquake path
$\E_{\mu,\c}$ on which $l_{\mu} \equiv \c$.  We denote the point
at which
$l_{\nu}$ is minimal on this path by $p_{\mu,\nu,\c}$, and set
$f_{\mu,\nu}(\c)=l_{\nu}(p_{\mu,\nu,\c})$.  For fixed $\mu,\nu$
and
variable $\c$, the points $p_{\mu,\nu,\c}$ define an analytic
path
$\F_{\mu,\nu}$, which we call a {\em critical line}; it meets
$\partial\F$
in the points $[\mu],[\nu]$.  The length functions
$l_{\mu},l_{\nu}$ are
monotonic on $\F_{\mu,\nu}$ and $f_{\mu,\nu}(\c)$ is continuous,
decreasing
monotonically from $\infty$ to $0$ on its domain $(0,\infty)$.

The following result completely describes the pleating planes
$\Pr{\mu,\nu}$; in particular it shows that $\Pr{\mu,\nu}$ can
be viewed as
an extension into $\QF$ of the critical line $\F_{\mu,\nu}$.

\begin{thmI}
\label{thm:plplanesI}
Let $(\mu,\nu)$ be measured laminations on $\torus$ with
$i(\mu,\nu)>0$.
Then $\Pr{\mu,\nu} $ is a non-empty connected non-singular
component of the
$\RR^2$-locus in $\QF-\F$ of the function $\L_{\mu,\nu}$.  The
restriction
of $\L_{\mu,\nu}$ to $\Pr{\mu,\nu}$ is a diffeomorphism to the
open region
under the graph of the function $f_{\mu,\nu}$ in $\RR^{+} \times
\RR^{+}$.
\end{thmI}

 The closure of $\Pr{\mu,\nu}$ in $\QF$ is the critical line
$\F_{\mu,\nu}
\subset \F$; it is mapped homeomorphically by $L_{\mu,\nu}$ to
the graph of
$f_{\mu,\nu}$.  The planes $\Pr{\nu,\mu}$ and $\Pr{\mu,\nu}$ are
disjoint
with common boundary $\F_{\mu,\nu}$ in $\QF$.  The set
$\Pr{\mu,\nu} \cup
\Pr{\nu,\mu} \cup \F$ is an $\RR^2$-locus in $\QF$ and the union
$\Pr{\mu,\nu} \cup \Pr{\nu,\mu} \cup\F_{\mu,\nu}$ may be
regarded as the
extension of the $\mu,\nu$ critical line to $\QF$.

The three components of the boundary of the image of
$\Pr{\mu,\nu}$ in $\RR^{+} \times \RR^{+}$ correspond to three
distinct parts of its closure in the set of algebraic limits of
groups in $\QF$.  As above, the component corresponding to the
graph of $f_{\mu,\nu}$ represents groups on the critical line
$\F_{\mu,\nu} \subset \F$.    For limit groups corresponding to
the axis $\lambda_{\mu}=0$ the component $\Omega^+$ has
degenerated and the support $|\mu|$ of $\mu$ is an ending
lamination; the bending measure of $\bch^-$, however, is still in
the projective class of $\nu$. Likewise, for limit groups
corresponding to the axis $\lambda_{\nu}=0$, the component
$\Omega^-$ has degenerated and the ending lamination is $|\nu|$.
The boundary point $(0,0)$ represents a doubly degenerate group,
unique by the results of~\cite{MinskyPT} (or~\cite{KMS} in the
rational case), with the two ending laminations $|\mu|$ and
$|\nu|$.

\smallskip
Theorems~\ref{thm:plinvarsI} and~\ref{thm:plplanesI} together
show that
we have a nice coordinate system on $\QF - \F$:
theorem~\ref{thm:plinvarsI}
 shows that the map to pleating invariants
is injective and theorem~\ref{thm:plplanesI} describes the
image.

\smallskip
The measured lamination $\mu$ is called {\em rational} if its
support is a
simple closed geodesic.  Such a geodesic can only belong to the
pleating
locus $|pl^{\pm}|$ if its representatives $V \in G$ are purely
hyperbolic
and hence have real trace.  Given any embedding $\QF$ into $
\CC^2$, the
generators of $G$ are holomorphic functions of the embedding
parameters and
$\tr V$ is a polynomial in the entries of the generators. In
particular,
given any elements $V,W \in G$ representing distinct simple
closed curves
$|\mu|,|\nu|$ on $\torus$, one can compute the position of the
critical
line $\F_{\mu,\nu}$.  If both laminations $\mu,\nu$ are
rational, we call
$\Pr{\mu,\nu}$ a {\em rational pleating plane}.
Theorem~\ref{thm:plplanesI} implies

\begin{thmI}
\label{thm:bdryI} Let $\mu,\nu $ be
rational laminations represented by non-conjugate elements $V,W
\in G$.
Then $\Pr{\mu,\nu}$ and $\Pr{\nu,\mu}$ are the unique components
of the
$\RR^2$-locus of the function $\tr V \times \tr W$ in $\QF-\F$
whose
closures meet $\F$ in $\F_{\mu,\nu}$.  On $\Pr{\mu,\nu} \cup
\Pr{\nu,\mu}$
the function $\tr V \times \tr W$ is non-singular and the
boundary of
$\Pr{\mu,\nu} \cup \Pr{\nu,\mu}$ can be computed by solving $\tr
V=\pm 2$
and $\tr W= \pm 2$ on this component.
\end{thmI}

We also prove
\begin{thmI}
\label{thm:ratlplanesdenseI}
The rational pleating planes are dense in $\QF$.
\end{thmI}

\medskip

In the late 1960's, Bers asked whether it was possible to find
the shape of
quasifuchsian space by explicit computation; one would expect
the punctured
torus to be the easiest case.  Partial results were obtained by
a number of
people, some using computational methods, among them
\cite{Kn,Porter,Wright}, others developing new tools and
techniques
\cite{Jorg,McM}.  For the punctured torus, the above results
give an
effective means of finding the boundary of the image of any
chosen
embedding of $\QF$ into $\CC^2$, answering Bers' question in
full.

\medskip

We also study the way in which the pleating planes fit together
 transversally to the real locus of $\L_{\mu,\nu}$.  This is
done by fixing
 the pleating invariants of one side of $\bch$; one can regard
this as
 analogous to fixing the ending invariant on one side in $\QF$,
to obtain
 the classical {\em Bers slice}~\cite{Berssim}.  Thus for a
fixed measured
 lamination $\mu$ and $\c > 0$, we define the {\em BM-slice}
 $BM^{+}_{\mu,\c}$ as the subset of $\QF$ on which $[pl^+]
=[\mu]$ and
 $\lambda_{\mu}=\c$.  The BM-slices are subsets of the {\em
quakebend
 planes} $\Q_{\mu,\c}$ obtained by Thurston's quakebend
construction along
 the measured lamination $\mu$ (see~\cite{EpM} and
 section~\ref{sec:twistsandqbs} below).  These are extensions of
the
 earthquake path $\E_{\mu,\c}$ into $\QF$.  Unlike the path
$\E_{\mu,\c}$
 which is completely contained in $\F$, the quakebend plane
$\Q_{\mu,\c}$
 is not totally contained in $\QF$.  We prove

\begin{thmI}
\label{thm:qbplanesI}
Let $\mu$ be a measured lamination on $\torus$ and let $\c > 0$.
 Then the
closures in $\QF$ of exactly two of the connected components of
$\Q_{\mu,\c} \cap (\QF-\F)$ meet $\F$.  These components are the
slices
$BM^{\pm}_{\mu,\c}$ and the closure of each slice meets $\F$
precisely in
the earthquake path $\E_{\mu,\c}$.  Furthermore, each slice is
simply
connected and retracts onto $\E_{\mu,\c}$.
\end{thmI}

Thus, just like the Bers slices, the $BM$-slices are complex
planes in
$\QF$ and like them, they foliate $\QF -\F$.  We note that while
the
boundary of the pleating planes consists of smooth curves, the
boundary of
a $BM$-slice is typically a fractal-like curve.  Pictures of
such curves
may be found in~\cite{KStop,PP,Wright}.

\medskip
The basis of the proofs of the above results are two important
theorems which control the local behavior of pleating invariants.
We call these the {\em limit pleating theorem } and {\em local
pleating theorem } respectively.  Roughly, the limit pleating
theorem states that if the pleating invariants of a sequence of
groups in $\QF$ converge, then the groups converge to an algebraic
limit;  furthermore the limit group is in $\QF$ provided the limit
pleating lengths are non-zero. It is closely related to Thurston's
double limit theorem~\cite{ThuH}, and also to the `Lemme de
fermeture' in \cite{BonO}.

The local pleating theorem makes essential use of the complex
length
function $\lambda_{\mu}$. As mentioned above, if $q \in \QF$,
then
$\lambda_{pl^+}(q) \in \RR$.  In general, the converse of this
result is
false; however the local pleating theorem gives a partial
result: if $q \in
\Pr{\mu}$ so that $\lambda_{\mu}(q) \in \RR$, then for $q' $
near $q$, the
condition $\lambda_{\mu}(q') \in \RR$ implies that $q' \in
\Pr{\mu}$. (As
discussed in the introduction of \cite{KSbend} this result does
not hold
for higher genus.)

\smallskip

The theory of quakebends as developed in~\cite{EpM} allows us to
extend the
earthquake paths $\E_{\mu,\c}$ into a family of holomorphic
planes
$\Q_{\mu,\c}$ in $\QF$.  We reduce the problem of studying the
sets
$\Pr{\mu,\nu}$ by restricting to the subset $\Pr{\mu,\nu,\c}$ of
$\Pr{\mu,\nu}$ on which the value of $\lambda_{\mu}$ is fixed at
$\c \in
\RR^+$.  For reasons that will be clear below, we call such a
set a {\em
pleating ray}.  In $\Q_{\mu,\c}$, the complex length
$\lambda_{\nu}$
restricts to a holomorphic function of one variable and it
follows from the
limit and local pleating theorems that $\Pr{\mu,\nu,\c}$ is both
open and
closed in the $\RR$-locus of $\lambda_{\nu}$ in $\Q_{\mu,\c}
\cap \QF$.

The fact that the pleating rays are non-empty and the discussion
of how
they meet Fuchsian space $\F$ results from the detailed study of
the
situation near $\F$ which was carried out in~\cite{KSbend}. We
also have
detailed information from~\cite{PS} about rays for which the
laminations
$\mu,\nu$ are rational and correspond to a pair of generators of
$\torus$.
Combining this information allows us to prove

\begin{thmI}
\label{thm:plraysI}
Let $\mu,\nu$ be measured laminations on $\torus$ with
$i(\mu,\nu)>0$ and
let $\c>0$.  Then the set $\Pr{\mu,\nu,\c} \subset \QF$ on which
$[pl^+]=[\mu], [pl^-]=[\nu]$ and $l_{\mu} = \c$, is a non-empty
connected
non-singular component of the $\RR$-locus of the restriction of
$\lambda_{\nu}$ to $\Q_{\mu,\c}$.  This restriction is a
diffeomorphism
onto its image $(0,f_{\mu,\nu}(\c)) \subset \RR^{+}$.
\end{thmI}

Theorem~\ref{thm:plplanesI}, and hence also
theorem~\ref{thm:plinvarsI},
are immediate consequences of this result.  We also easily
deduce
theorem~\ref{thm:qbplanesI}.

\smallskip

For groups on the boundary of $\QF$, at least one of the
components
$\Omega^{\pm}$ degenerates and it is clear that our pleating
invariants
extend naturally to the corresponding ending laminations for
which the
length (and also the complex length) is always $0$.  It is also
clear that
these invariants should also characterize boundary groups;
careful analysis
requires the study of generalized Maskit slices in which the
fixed ending
lamination is irrational, see~\cite{McMeq}.

\smallskip
 The reader is referred to \cite{MinskyPT} for a good outline of
the
history relating to the study of punctured torus groups.

\smallskip
Some of the ideas of this paper, in particular the relation of
pleating planes to the Kerckhoff picture of $\F$ and the idea of
looking at the $BM$-slices, grew out of discussions with John
Parker, and we should like to thank him for his input into this
work.  We should also like to thank our referees for their
detailed reading  of earlier versions of this paper, in
particular, for having signalled, in view of the examples in
\cite{KerckTh}, a gap in our proof of theorem~\ref{thm:thmq}, as
well as having suggested a more direct proof of
lemma~\ref{lemma:retoinfty} and a simplification of the proof of
theorem~\ref{thm:realI}. We would also like to thank Yair Minsky
for conversations which helped us precisely locate the above
mentioned gap, and Francis Bonahon and Cyril Lecuire for very
useful discussions about how to rectify it.

\smallskip
The paper is organized as follows. Section~\ref{sec:background}
contains
background on the punctured torus, geodesic laminations and
surfaces.
Section~\ref{sec:fuchsian} explains the picture of earthquake
paths and
critical lines in $\F$ and in section~\ref{sec:hypmflds} we
review results
on pleated surfaces and the convex hull boundary. We prove the
limit
pleating theorem in section~\ref{sec:closure}.  In
section~\ref{sec:cxlenlam} we show how to complexify the length
functions
and show that the complex length of the pleating locus is real.
In
section~\ref{sec:twistsandqbs} we review results about
quakebends and the
convex hull boundary and then in section~\ref{sec:openness}
prove the local
pleating theorem. We also derive various important consequences
of this
result, including the proof of
theorem~\ref{thm:ratlplanesdenseI}. In
section~\ref{sec:Pleating Rays} we prove our main results,
theorem~\ref{thm:plraysI} on pleating rays and
theorem~\ref{thm:plplanesI}
on pleating planes.  In section~\ref{sec:BMslices} we study
BM-slices,
proving theorem~\ref{thm:qbplanesI}, and we conclude in
section~\ref{sec:ratlplplanes} with a discussion of rational
pleating
planes, computation, and some explicit examples.  For
readability, the
proofs of three technical results are deferred to the appendix.

\section{Background}
\label{sec:background}

\subsection{Punctured torus groups and markings}
 \label{sec:normalization}

Let $\torus$ be a torus with one puncture and a fixed
orientation.  Any
pair of simple closed loops on $\torus$ that intersect exactly
once are
free generators of $\pi_1(\torus)$.  Let $(\alpha,\beta)$ be
such an
ordered pair of free generators, chosen so that their commutator
$\alpha
\beta \alpha ^{-1} \beta ^{-1}$ represents a loop around the
puncture that
is positively oriented around the component of $\torus$ not
containing the
puncture.  The ordered pair $(\alpha,\beta)$ is called a {\em
marking}.

A {\em punctured torus group} is a discrete subgroup
$G \subset  PSL(2,\CC)$ that  is
the image of a faithful
representation $\rho$ of  $\pi_1(\torus)$
such that the image of the loop around the puncture is
parabolic.
If  $(\alpha,\beta)$ is a   marking of $\torus$,
and if  $
A = \rho(\alpha),B =\rho(\beta)$, then
the commutator  $K = ABA^{-1}B^{-1}$
 is parabolic and  the ordered pair
$(A ,B) = (\rho(\alpha),\rho(\beta))$ is called a
{\em marking} of $G$.
If  $\gamma$ is any simple closed curve on $\torus$, then we can
always
choose a curve $\delta$ such that $(\gamma, \delta)$ is a
marking
of $\torus$.
Setting
$\rho(\gamma)=V, \rho(\delta)=W$, then
  all possible markings $(V,W'), W' \in G$ of $G$ are of
the form $(V,V^mW), m \in \ZZ$.

The group $G$ is {\em quasifuchsian} if the regular set $\Omega$
consists
of two non-empty simply connected invariant components
$\Omega^{\pm}$.
The limit set $\Lambda(G)$ is topologically a circle.
{\em Quasifuchsian space} $\QF$ is the space of marked
quasifuchsian
punctured torus groups modulo conjugation in $PSL(2,\CC)$; it
has
 a holomorphic structure induced from the
natural holomorphic structure of $SL(2,\CC)$.
{\em Fuchsian
space} $\F$ is the subset such that the components
$\Omega^{\pm}$
are round
disks.  It is canonically isomorphic to the Teichm\"uller space
of
marked
conformal structures on
$\torus$.

The quotients
$\Omega^{\pm}/G$ are punctured tori with conformal structures,
and
hence
also orientations, inherited from $\Chat$; the orientations of
$\Omega^+/G$
and $\torus$ agree whereas those of $\Omega^-/G$ and $\torus$
are
opposite.
 This means  $\Omega^+(G)$ is the
component such that  $A^-,B^+, A^+,B^-$ occur in
counterclockwise order around its boundary $\Lambda(G)$, where,
for
a loxodromic
 $g \in SL(2,\CC)$, $g^+$ and $g^-$ denote its attracting
and repelling fixed
points  respectively.  Thus an alternative
way to choose a marking of $G$ is to choose {\em any} pair of
generators
$X,Y$ of $G$, and to specify the choice of $\Omega^+$ by
choosing
it to be the
component such that the fixed points $X^-,Y^+, X^+,Y^-$ run
counterclockwise around its boundary.

A point $q \in \QF$ represents an equivalence class of marked
groups
in $ PSL(2,\CC)$.
 We choose once and for all a  triple of distinct
 points in $\Chat$ and let $G= G(q)$ denote the
representative normalized by choosing
$A^-,A^+,K^{\infty}$  to be this this fixed triple, where
$K^{\infty}$ is the fixed point of the parabolic $K$. We will
refer
to this
as the {\em standard normalization}.
 If it is clear from the
context, for readability, we suppress the dependence on $q$.

 Note that throughout this paper, $\QF$ and $\F$ refer to the
special case of the once punctured
torus $\torus$ only.

 \subsection{ Laminations}
\label{sec:laminations}

Let $\Sigma$ be a hyperbolic surface.  We denote by ${\cal S}$
the set of all simple closed geodesics on $\Sigma$.
There is one such geodesic
in each free homotopy class of simple closed non-boundary
parallel
loops,
and the  set $\S$  is
independent of the hyperbolic structure on $\Sigma$.

{\em Geodesic laminations} were introduced by
Thurston~\cite{ThuN}
as a
generalization of simple closed geodesics.  A geodesic
lamination
on
$\Sigma$ is a closed
set that is a union of pairwise disjoint simple geodesics called
its {\em
leaves}.  We denote by $GL= GL(\Sigma)$
 the set of of
all geodesic laminations on $\Sigma$; $GL(\Sigma)$ is also
independent of the hyperbolic structure, see
e.g.  \cite{CEpG} section 4.1.4 and \cite{KSconvex} section 3.7.

The
Hausdorff
topology on the set of closed subsets of $\Sigma$ induces a
topology on
$GL$.   Two laminations are
close in this topology if any long segment of a leaf of either
one
is
closely approximated by a long segment of a leaf of the other.
See
\cite{CEpG,EpM,Penner} for a complete discussion.

A {\em measured lamination} $\mu$ on $\Sigma$ is a geodesic
lamination,
called the {\em support} of $\mu$ and denoted $|\mu|$, together
with a
transverse measure, also denoted $\mu$.  We denote the set of
all
measured
laminations on $\Sigma$ by $ML(\Sigma)$.  The space $ML$ is
topologized by
defining laminations to be close in $ML$ if the measures they
assign to
any finite set of transversals are close, for details see
\cite{CEpG} or
\cite{KSconvex}. Notice that the support of any measured
lamination
always avoids a definite neighborhood of each cusp.
The relationship between the   topologies on $ML$ and $GL$
 is discussed in
section~\ref{sec:convergencelemma}  below.

Any element $\gamma \in \S$ carries a natural transverse measure
$\delta_{\gamma}$  which assigns unit mass to each intersection
with
$\gamma$.
We call a measured geodesic lamination on $\Sigma$ {\em
rational} if its support is a union of curves in $\S$.
The maximum number of disjoint loops in $\S$
on the punctured torus $\torus$ is  one, so that
 rational
measured laminations are  of the form $\mu=
\k \delta_{\gamma}, \k>0$.  We   denote the
set of all rational measured laminations on $\Sigma$ by
$ML_Q(\Sigma)$;
the set $ML_Q$ is dense in $ML$.

Two measured laminations $\mu, \mu' \in ML$ are {\em
projectively
equivalent} if $|\mu|=|\mu'|$ and if there exists $\k >0$ such
that
for any
arc $\sigma$ transverse to the leaves of $|\mu|$, $\mu'(\sigma)
=
\k
\mu(\sigma)$.  We write $[\mu]$ for the projective class of $\mu
\in
ML(\Sigma)$.
We denote the set of projective equivalence classes on $\Sigma$
by
 $PML(\Sigma)$.
  It is well known that
 $PML(\torus)$ is homeomorphic to $S^1 \simeq \RR
 \cup \{\infty\}$ (see for example \cite{ThuN}).

   The length $l_{\gamma}$ of a geodesic $\gamma \in {\cal S}$
generalizes
  to arbitrary laminations. Let $\phi$ represent a hyperbolic
structure on
  $\Sigma$.  For $\mu \in ML$, the length $l_{\mu}(\phi)$ is the
total
  mass, on the surface with structure $\phi$, of the measure
that
is the
  product of hyperbolic distance along the leaves of $\mu$ with
the
  transverse measure $\mu$.  In particular, if $\mu \in
ML_Q(\Sigma)$ with
  $\mu= \delta_{\gamma}$, then $l_{\gamma}=\int_{\Sigma}
d\delta_{\gamma}
  ds$ is just the hyperbolic length of $\gamma$.

Clearly, if $\mu' = \k \mu$ then $l_{\k \mu} = \k l_{\mu}$. We
define $$[\mu,l_{\mu}]\stackrel {\rm def}{=} \{\k \mu,\k l_{\mu}
\in ML \times
\RR^+: \k >0 \} $$ and call it the {\em projective class} of the
pair
$(\mu,l_{\mu})$.

The geometric intersection number $i(\gamma,\gamma')$ of two
geodesics
$\gamma,\gamma' \in {\cal S}$ extends to a continuous function
$i(\mu,\nu)$
on $ML(\Sigma)$ (see for example~\cite{KerckN}).  For
$\Sigma=\torus$,
  $i(\mu,\nu) > 0 $ is equivalent to $[\mu] \neq [\nu]$.  We
also recall the well known fact that on $\torus$, measured
laminations are
uniquely ergodic; that is, if $\mu,\mu' \in ML(\torus)$ with
$|\mu|=|\mu'|$, then $[\mu]=[\mu']$.

\subsection{The convergence lemma}
\label{sec:convergencelemma}

In general, laminations which are close in $ML$ may not be close
in
the
Hausdorff topology on $GL$.  For example,  one can put a
transverse
measure
$\nu'$ on a  long closed geodesic $\gamma'$
  spiralling in to a   closed geodesic $\gamma$ with transverse
measure
  $\nu$, such that
  $\nu,\nu'$ are close in $ML$
  but $\gamma'$ has arcs far from  $\gamma$.
 A sequence of laminations  may converge in $ML$ to  a
  measured lamination $\nu_0$ with support
  in one part of $\Sigma$,
  while simultaneously limiting on a closed curve with support
disjoint
   from $|\nu_0|$.

The following lemma  gives conditions  under which
 Hausdorff convergence is a consequence of  convergence in $ML$.
We note that the lemma depends  crucially on the fact that
on $\torus$, any irrational
measured lamination is maximal.
As stated, it is false for more general surfaces, and it is
false
if  $\nu_0 \in ML_Q(\torus)$.

\begin{lemma}
\label{lemma:convergence}  Suppose that $\nu_0 \in ML(\torus)
 -  ML_Q(\torus)$, and that
$\nu $ and $\nu_0$ are close in $ML(\torus)$. Then
$|\nu|$ and $|\nu_0|$ are close in the Hausdorff topology on
$GL(\torus)$.
\end{lemma}

This lemma is proved in appendix~\ref{app:convergence}.

From now on, unless specifically stated, $GL, ML, PML$ will always
refer to $\torus$.

\section{Fuchsian space}
\label{sec:fuchsian}

Kerckhoff and Thurston used earthquake deformations to study the
set of
hyperbolic structures on a surface $\Sigma$.  For $\torus$ the
description
is especially simple.  For an unpunctured torus, the
Teichm\"uller space is
a disk. Thinking of this disk as the hyperbolic plane $\DD$ with
boundary
circle $S^1$, for each boundary point $\xi$ there is a foliation
of $\DD$
by horocycles tangent to ${\partial \DD}$ at $\xi$.  Joining
each pair of
distinct boundary points $\xi,\eta$ is a unique geodesic
$\gamma_{\xi,\eta}$ which, for fixed $\xi$ and varying $\eta$,
give another
foliation of $\DD$.  It follows from Kerckhoff's results
\cite{KerckN,KerckLM} and Thurston's compactification of
Teichm\"uller
space \cite{FLP} (see also \cite{GardMas}), that there is an
analogous
picture for $\F$, the Teichm\"uller space of $\torus$.  This
picture is
certainly well known and described for Teichm\"uller spaces of
compact
surfaces in \cite{KerckLM}.  As it is of central importance for
us we
explain it in detail here.

Since the torus is homogeneous, $\F$ is holomorphically the same
as the
Teichm\"uller space of the unpunctured torus, namely $\DD$.  The
Thurston
boundary of $\F$ is naturally identified with the circle $S^1$.
The
classical Fenchel-Nielsen coordinates for $\F$ are the length
$l_{\alpha}$
of a generating curve $\alpha$ and a corresponding twist
parameter
$t_{\alpha}$.  In \cite{KerckN, KerckEA}, the Fenchel-Nielsen
deformation
defined by varying the twist parameter $t_{\alpha}$ is
generalized to a map
$\E_{\mu}(t):\F \rightarrow \F$ defined relative to a measured
lamination
$\mu \in ML$.  The map $\E_{\mu}(t)$ is called the time $t$ {\em
earthquake} along $\mu$; when needed for clarity we write the
parameter $t$
as $ t_{\mu}$.  The family $\E_{\mu}(t), t \in \RR$ is a one
parameter
family of deformations of $\F$; in particular $\E_{\mu}(0) =
{\rm id}$.

For $p \in \F$, we define the {\em earthquake path} along $\mu$
through $p$
by $$\E^p_{\mu} = \{ \E_{\mu}(t)(p) \in \F : t \in \RR \}.$$
Clearly,
$\E^p_{\mu}$ is invariant under the earthquakes $\E_{\mu}(t)$.
In
\cite{KerckEA}, Kerckhoff showed that $\E^p_{\mu}$ is a real
analytic path
in $\F$.  Along $\E^p_{\mu}$, the length $l_{\mu}$ is constant.
Thus for
every $p \in \F$, $\E_{\mu}(t)(p)$ tends to the same point
$[\mu] \in
\partial\F$ as $t \to \pm\infty$.

In \cite{KerckN}, Kerckhoff showed that if $\nu \in ML$ with
$i(\mu,\nu) >
0$, then along an earthquake path $\E^p_{\mu}$, the length
$l_{\nu}$ is a
strictly convex real analytic function of $t$ and $l_{\nu}(t)
\to \infty$
as $t \to \pm\infty$. Thus $l_{\nu}$ has a unique minimum on
$\E^p_{\mu}$;
at this point we say that $l_{\nu}$ is minimal with respect to
$\E^p_{\mu}$.  Wolpert~\cite{Wolp} showed in addition, that at
the minimum,
$d^2l_{\nu}/dt_{\mu}^2 > 0$.

It follows from the anti-symmetry of the derivative formula
$$dl_{\nu}/dt_{\mu} = \int_{\torus} \cos \theta d\mu d\nu,$$
(where
$\theta$ is the angle, measured counterclockwise, from a leaf of
$|\mu|$ to
a leaf of $|\nu|$ at each intersection point of the laminations
$|\mu|,|\nu|$), that the minimum points for $l_{\nu}$ along
$\E_{\mu}$ and
$l_{\mu}$ along $\E_{\nu}$ coincide, and that at this minimum
point $p$ we
have $D\E_{\nu}(t_{\nu})(p) = -D\E_{\mu}(t_{\mu})(p)$.

The results which follow are simple consequences of Kerckhoff's
results
applied to $\torus$.

\begin{prop}
\label{prop:earthqkpath} For any $\c \in \RR^+$ and $\mu \in
ML$, there is at most one earthquake path $\E^p_{\mu}$ along
which
$l_{\mu}=\c$.
\end{prop}

\begin{proof}
Suppose that there are two such paths, $\E_1,\E_2$.  They are
clearly
disjoint, moreover since $\F \cup \partial\F$ is a closed disk
and both
$\E_1$ and $\E_2$ meet $\partial\F$ at the same point $[\mu]$,
one path,
$\E_1$ say, separates $\F \cup \partial\F$ so that one component
of the
complement contains both $\E_2$ and $\partial\F - \{[\mu] \}$.
Choose $\nu
\in ML$ with $i(\mu,\nu)>0$ and let $p$ be the minimum point for
$l_{\nu}$
on $\E_1$.  Then $\E_{\nu}^p$ must also cut $\E_2$ at a point
$p'$.  Since
$p$ is the unique minimum point for $l_{\mu}$ on $\E_{\nu}^p$,
and since
$l_{\mu}(p)=l_{\mu}(p')$ we have a contradiction.
\end{proof}

\medskip
We denote the unique earthquake path on which $l_{\mu}=\c$ by
 $\E_{\mu,\c}$.  It follows easily from
proposition~\ref{prop:monotone}
 below that $\E_{\mu,\c} \neq \emptyset$.  Since for $s > 0$,
$\E_{s\mu}(t)
 = \E_{\mu}(st)$ and $l_{s\mu}=sl_{\mu}$, we have
 $\E_{s\mu,s\c}=\E_{\mu,\c}$.  For $\nu \notin [\mu]$, we denote
the
 minimum point for $l_{\nu}$ on $\E_{\mu,\c}$ by
$p(\mu,\nu,\c)$.  We
 define a function $f_{\mu,\nu}: \RR^+ \rightarrow \RR^+$ by
 $f_{\mu,\nu}(\c) = l_{\nu}(p(\mu,\nu,\c))$.

\smallskip
  For each pair $\mu,\nu \in ML \times ML $, $\mu \not\in
[\nu]$, set
  $$\F_{\mu,\nu} = \{ p \in \F | dl_{\nu}/dt_{\mu}(p)=0 \}$$
Note that
  $\F_{\mu,\nu}$ depends only on $[\mu],[\nu]$.  We call
$\F_{\mu,\nu}$ the
  $\mu,\nu$-{\em critical line}.  This is justified by the
following
  proposition.

 \begin{prop}
\label{prop:monotone} For each pair $\mu,\nu$, $i(\mu,\nu) >0$,
 the locus $\F_{\mu,\nu}$ is a real analytic path in ${\F}$ with
endpoints
  at $[\mu]$ and $[\nu]$ in $\partial\F$.  Both $l_{\mu}$ and
$l_{\nu}$ are
  strictly monotonic on $\F_{\mu,\nu}$ and vary from $0$ to
$\infty$ in
  opposite directions.  \end{prop}

 \begin{proof} By Wolpert's result, $d^2l_{\nu}/dt_{\mu}^2 > 0$
at every
point of $\F_{\mu,\nu}$.  Therefore $\F_{\mu,\nu}$ is a union of
real
analytic arcs.

We claim the function $l_{\mu}$ is strictly monotonic on each
component of
$\F_{\mu,\nu}$.  If not, there is an earthquake path
$\E_{\mu,\c}$ that
meets $\F_{\mu,\nu}$ in two distinct points.  Both these points
are
critical for $l_{\nu}$ on $\E_{\mu,\c}$ which is impossible.

Since $l_{\mu}$ is real analytic, its restriction to
$\F_{\mu,\nu}$ is open
and proper and hence its range must be $(0,\infty)$.  Clearly,
as
$l_{\mu}(p) \to 0$ along $\F_{\mu,\nu}$, we have $p \to [\mu]
\in
\partial\F$.  Thus each component of $\F_{\mu,\nu}$ is an
embedded arc with
endpoints $[\mu]$ and $[\nu]$ in $\partial\F$.

If $\F_{\mu,\nu}$ had two components, then, for some $\c>0$, we
could find
a path $\E_{\mu,\c}$ intersecting both components of
$\F_{\mu,\nu}$.  Thus
$l_{\nu}$ would be minimal at two points on $\E_{\mu,\c}$ which
is
impossible.

 By the anti-symmetry in the formulas, we see that $l_{\nu}$
also varies
monotonically from $0$ to $\infty$ along $\F_{\mu,\nu}$ but in
the opposite
direction.
\end{proof}

\begin{cor}
\label{cor:earthquakepaths}
For any $\c \in \RR^+$ and $\mu \in ML$ there is a unique
earthquake path
$\E^p_{\mu}$ along which $l_{\mu}=\c$.
\end{cor}

\begin{remark} {\rm In~\cite{KerckLM}, Kerckhoff proves that
given
$(\mu,\nu) \in ML$ with $i(\mu,\nu)>0$ and such that $\mu,\nu$
fill up the
surface (that is, the complement of their union consists of
pieces which
are either simply connected or a neighborhood of the puncture),
then for
each $t \in (0,1)$ there is a unique $p \in \F$ at which the
function $t
l_{\mu}(p) + (1-t) l_{\nu}(p)$ attains minimum.  As $t$ varies
keeping
$\mu, \nu$ fixed, the set of these minima is a line.  For the
punctured
torus, any pair $(\mu,\nu) \in ML$ with $i(\mu,\nu)>0$ fills up
the
surface.  While not strictly needed for our development, the
following
lemma confirms that for the punctured torus, Kerckhoff's {\em
line of
minima} is identical with our critical line, see
also~\cite{KerckLM}
theorem 3.4.}
\end{remark}

\begin{lemma}
\label{lem:lineofminima}
 Suppose that $i(\mu,\nu)>0$. Then
$ p \in \F_{\mu,\nu}$ if and only if $p$ is the
global minimum for some function
$tl_{\mu}(p) + (1-t) l_{\nu}(p)$ for some $t \in (0,1)$.
\end{lemma}

\begin{proof}
At a minimum of $tl_{\mu}(p) + (1-t) l_{\nu}(p)$, since
$l_{\mu}$ is
constant along the earthquake path $\E_{\mu}(p)$, we find
$dl_{\nu}/dt_{\mu}(p) = dl_{\mu}/dt_{\nu}(p)= 0$ so that $ p \in
\F_{\mu,\nu}$.  Conversely, if $dl_{\nu}/dt_{\mu}(p) = 0$, the
earthquake
paths $\E_{\mu}(p)$ and $\E_{\nu}(p)$ must be tangent at $p$
because $p$ is
the unique minimum of $l_{\nu}$ on $\E_{\mu}(p)$.  Thus
$\E_{\mu}'(p) = - k
\E_{\nu}'(p)$ for some $k \neq 0$, where $'$ denotes the tangent
vector to
the corresponding earthquake path. From the derivative formula
$dl_{\nu}/dt_{\mu} = -dl_{\mu}/dt_{\nu}$ it follows that $k>0$.
We get $ d
l_{\eta} /d\tau_{\mu}(p) = - kd l_{\eta} /d\tau_{\mu}(p)$ for
any $\eta \in
ML$, which, using the derivative formula again, gives $ d
l_{\mu}
/d\tau_{\eta}(p) = - kd l_{\nu} /d\tau_{\eta}(p)$.  Since the
tangent
vectors $\E_{\eta}'(p), \eta \in ML$ certainly span the tangent
space to
$\F$ at $p$, we must be at a critical point of $l_{\mu} + k
l_{\nu} $.
\end{proof}

Using the identification of the critical line $\F_{\mu,\nu}$
with the Kerckhoff  line of minima,
the following proposition follows immediately
from~\cite{KerckLM} theorem~2.1.
Here is another proof.

\begin{prop}
\label{prop:fgammafol} Fix $[\mu] \in PML$.  Then
the arcs $\F_{\mu,\nu}, [\nu] \in PML-\{ [\mu]\}$
 are pairwise disjoint and foliate $\F$.
\end{prop}

\begin{proof}
 Given  $p \in \F$, following Kerckhoff we
define  $\beta = \beta_p:ML  \rightarrow T_p\F $ to be the map
 which  takes $\mu \in ML$ to
$D\E_{\mu}(t_{\mu})(p)|_{t_{\mu}=0}$,
the derivative  with respect to $t_{\mu}$ of the earthquake path
$\E_{\mu}(t_{\mu})(p)$ through $p$ evaluated at $p$.
By~\cite{KerckLM}, Theorem 3.5 the
map $\beta$ is a homeomorphism.
Clearly, $\beta$ induces a homeomorphism between $PML$ and the
set
of rays through the origin in $T_p\F $.

 Suppose $[\mu], [\nu],[\nu'] \in PML$ are distinct, and suppose
that $p
\in \F_{\mu,\nu} \cap \F_{\mu,\nu'}$.  Pick representatives
$\mu,\nu,\nu'$
of $[\mu],[\nu],[\nu']$ and let $c= l_{\mu}(p),
d=f_{\mu,\nu}(\c),
d'=f_{\mu,\nu'}(\c)$.  The earthquake paths $ \E_{\nu,\d}$ and
$\E_{\nu',\d'}$ both go through $p$ and, because $l_{\mu}$ is
minimal at
$p$ with respect to both $\E_{\nu}$ and $\E_{\nu'}$, from the
derivative
formula we see that $D\E_{\nu}(t_{\nu})(p)|_{t_{\nu}=0} =
D\E_{\nu'}(t_{\nu'})(p)|_{t_{\nu'}=0}$.  By the injectivity of
$\beta_p$ on
$PML$, $[\nu]=[\nu']$.

Now let $p \in \F$.  By the surjectivity of $\beta$, there is
some $\nu \in
 ML$ such that $D\E_{\nu}(t_{\nu})(p)|_{t_{\nu}=0} =
 -D\E_{\mu}(t_{\mu})(p)|_{t_{\mu}=0}$.  Therefore the earthquake
paths
 $\E_{\nu,l_{\nu}}(p)$ and $\E_{\mu,l_{\mu}}(p)$ are tangent at
$p$.  Since
 earthquake paths can intersect in at most two points it follows
that
 $l_{\nu}$ is minimal at $p$ with respect to $\E_{\mu}$, so that
$p \in
 \F_{\mu,\nu}$.

These two facts show that the sets $\F_{\mu,\nu}$ foliate $\F$.
\end{proof}

We shall also need
\begin{cor}
\label{cor:X}
For fixed $\mu \in ML$, $\c \in \RR^+$, the map
$\psi:PML - \{[\mu]\} \rightarrow  \E_{\mu,\c}$, $\psi([\nu]) =
p(\mu,\nu,\c)$, is a homeomorphism.
\end{cor}

\begin{proof}
Proposition~\ref{prop:fgammafol} shows that $\psi$ is
well defined and a bijection.  It is also clear, thinking of
$PML - \{[\mu]\}$ and $ \E_{\mu,\c}$ as intervals,
that $\psi$ is monotonic.  The result follows.
\end{proof}

Corollary~\ref{cor:earthquakepaths} implies that for $\mu \in
ML$, the
paths $\E_{\mu,\c}$, $\c \in \RR^+$ are pairwise disjoint and
foliate $\F$.
This is the analogue of the foliation of the hyperbolic disk
$\DD$ by
horocycles tangent at to a point on the boundary.  Likewise, the
critical
lines   $\F_{\mu,\nu}$ are the analogue of the geodesics in
$\DD$ joining a
pair of distinct points in $S^1$.  For fixed $[\mu]$ the
foliation by
leaves $\F_{\mu,\nu}$, $[\nu] \neq [\mu]$ is clearly transverse
to that by
the earthquake paths $\E_{\mu,\c}$.

This is the picture that we shall extend to $\QF$ below.

\section{Hyperbolic 3-manifolds}
\label{sec:hypmflds}

 \subsection{The pleating locus}

\label{sec:pllocus}

Let $q \in \QF$ and let $G=G(q)$ be a group representing $q$
with
the
standard normalization of
 section~\ref{sec:normalization}.  The group $G$ acts as a
discrete
group
of isometries of hyperbolic space $\H3$ and the quotient
hyperbolic
manifold $M=\H3/G$ is a product $\torus \times (-1,1)$.  If $G$
is
quasifuchsian, but not Fuchsian, the boundary $\bch$ of the
hyperbolic
convex hull $\C$ of $\Lambda$ in $\H3$ has two components
$\bch^{\pm}$ each
of which is also $G$-invariant. Each quotient $\bch^{\pm}/G$ is
homeomorphic to $\torus$, see for example~\cite{KSconvex}
proposition 3.1.
The metric induced on the components $\bch^{\pm}$ from $\H3$
makes
them
{\em pleated surfaces}.  This means, see for example~\cite{EpM},
that there
are surjective isometric maps $\psi^{\pm} : \DD \rightarrow
\bch^{\pm}$
such that for each point $z$ in $\DD$ there is at least one
geodesic
segment through $z$ that is mapped to a geodesic segment in
$\bch^{\pm}$.
The group G acts as a discrete group of isometries on each
component
$\bch^{\pm}$.  Since $\bch^{\pm}/G$ are both homeomorphic to
$\torus$,
these two groups of isometries are both isomorphic to
$\pi_1(\torus)$ and
inherit a marking in the obvious way. (The marking on
$\bch^{-}/G$
has its orientation reversed.)
 The isometries $\psi^{\pm}$ induce
isomorphisms to marked Fuchsian punctured torus groups
$F^{\pm}=F^{\pm}(q)$
acting on $\DD$, which we may again take to have the standard
normalization. We refer to both the marked groups
$F^{\pm}(q)$ and the quotients $\DD/F^{\pm}(q)$ as the {\em flat
structures} of either the surfaces $\bch^{\pm}/G(q)$ or of their
universal covers $\bch^{\pm}(q)$.

The bending laminations of $\bch^{\pm}/G$
 carry  natural transverse measures,
 the {\em bending measures} $pl(q)^{\pm}$, see
\cite{EpM,KSconvex}.
 The underlying laminations $|pl(q)^{\pm}|$
 are the {\em pleating loci} of $G$.
If $G \in \QF$ is a Fuchsian group acting
 on
the hyperbolic disk  $\DD \subset \H3$, then $\C=\DD$ is
degenerate
and
 we regard  $\bch$
and $\bch/G$ as
2-sided surfaces,
each side of which is a pleated surface with empty pleating
locus
(and zero measure).

The following proposition  follows immediately
from ~\cite{KSbend}  proposition~3.3 and corollary~3.4.

\begin{prop}
\label{prop:neq} Suppose that $q \in \QF-\F$.
Then the projective class of the bending measure
cannot be the same on both sides
of the convex core; that is, $[pl^+(q)] \neq [pl^-(q)]$.
\end{prop}

 \begin{remark} {\rm The work in~\cite{KSbend} depends heavily on the $\lambda$-lemma
 and the theory of holomorphic motions
 which is usually stated in the context of one complex variable. In the
 present case we shall be studying families of groups parameterized by a
 two dimensional complex manifold; in fact the theory of holomorphic
 motions extends to motions over any complex manifold.  see~\cite{McS}.}
 \end{remark}

In  \cite{KSconvex} we prove:
\begin{thm}
\label{thm:plcont} The maps $\QF \rightarrow \F$,
$q \mapsto F^{\pm}(q)$, and $\QF-\F \rightarrow ML$,
$q \mapsto pl^{\pm}(q)$  are
continuous.
\end{thm}

\subsection{Pleating Varieties}
\label{sec:plvars}

Given $\mu \in ML$ we set
$$\textstyle \Pr{\mu}^{\pm}= \{ q \in \QF : [pl^{\pm}(q)]=[\mu]
\}
\mbox{ and }
 \Pr{\mu}= \Pr{\mu}^+ \cup \Pr{\mu}^-.$$
We call these sets the $\mu$-{\em pleating varieties}.

\noindent Given the ordered pair $(\mu,\nu) \in ML \times ML$,
we set $$
\Pr{\mu,\nu} = \{ q \in \QF : [pl^+(q)]=[\mu], [pl^-(q)]=[\nu]
\}.$$ We
call this set the $\mu,\nu$-{\em pleating plane}.  Note that two
these
definitions depend only on the projective classes $[\mu],[\nu]$.

Finally, given the ordered pair $\mu,\nu \in ML \times ML$, and
$\c>0$ we
set $$ \Pr{\mu,\nu,\c} = \{ q \in \Pr{\mu,\nu} :l_{\mu}(q) =\c
\}.$$ We
call this set a {\em pleating ray}.  Note that for $s \in
\RR^+$,
$\Pr{\mu,\nu,\c} = \Pr{s\mu,\nu,s\c}$.  Thus $\Pr{\mu,\nu,\c}$
depends on
the projective class of the {\em pair} $(\mu, \c)$, (recall
section~\ref{sec:laminations}), and on the projective class
$[\nu]$.

Theorems sections~\ref{sec:plplanes} and~\ref{sec:Pleating Rays}
below will
justify the terminology {\em planes} and {\em rays}.

Proposition~\ref{prop:neq}
 implies $\Pr{\mu,\mu} = \emptyset$.
 It is also clear that $\Pr{\mu,\nu} \cap \Pr{\mu',\nu'} =
\emptyset$ unless $[\mu]=[\mu'],[\nu]=[\nu']$.  In particular
 $\Pr{\mu,\nu} \neq \Pr{\nu,\mu}$ whenever $i(\mu,\nu)>0$.

\begin{remark} {\rm Whether a  group is in $\Pr{\mu,\nu}$ or in
$\Pr{\nu,\mu}$
 depends on our conventions in labelling the  sides $\bch^{\pm}$
of
 $\bch$. This is based on the labelling of
  the  components of the regular set $\Omega^{\pm}$.
 The point  here is that two groups
 which differ only in the labelling of their  $+$ side and their
$-$ side
 are {\em not the same} as marked groups in $\QF$.}
\end{remark}

The main result of  \cite{KSbend} is that the
pleating varieties are non-empty.  Precisely, we prove

\begin{thm}
\label{thm:bend} Let $\mu,\nu  \in ML$,
$[\mu] \neq [\nu]$.  Then $\Pr{\mu,\nu} \neq \emptyset$.
\end{thm}

We shall need to study the ideas in the proof of
this result in some detail; see~\ref{sec:genqbs} below.

\subsection{Lamination length in  $M= \H3/G$.} \label{sec:length}

For the proof of theorem~\ref{thm:thmq} below, we need also to
discuss briefly the length $l_{\mu}(M )$ of a measured lamination
$\mu \in ML$ in the hyperbolic $3$-manifold $M= \H3/G$.  First,
suppose that $\mu = \delta_{\gamma}$ where $\gamma \in {\cal S}$
is represented by an element $V \in G $.  The multiplier
$\lambda_V$ is related to its trace by the formula $\tr V = 2
\cosh \lambda_{V }/2$.  The translation length of $V$, $ \Re
{\lambda_{V }}$, is the minimum distance that $V$ moves a point in
$\H3$.  Equivalently it is the length of the geodesic
representative of $\gamma$ in $M$, so that $l_{\delta_{\gamma}}(M
)= \Re {\lambda_{V }}$.

In \cite{ThuN}, p.9.21 and \cite{BonV3}, p.117, it is shown that
this definition can be extended by linearity and continuity to
define the lamination length $l_{\mu}(M )$ for an arbitrary $\mu
\in ML$. In the proof of theorem~\ref{thm:thmq} below, we shall
need to make crucial use of the fact that one can extend this
definition continuously to the algebraic closure of $\QF$.

Suppose $G$ is a (discrete) punctured torus group associated to
the faithful representation $\rho: \pi_1(\torus) \to G \subset
PSL(2,\CC)$. This representation {\it marks} the associated
hyperbolic 3-manifold $M=\HH^3/G$. One says  that a lamination
$|\mu|$ on $\torus$ is {\it realized} in $M$ relative to the
marking $\rho$, if there is a Fuchsian group $\Gamma$, a
homeomorphism $h: \torus \to S= \HH^2/\Gamma$, and a pleated
surface  $f: S  \to M$ with pleating locus containing $|\mu|$,
such that $fh$ induces $\rho$.

Let $AH(\torus)$ denote the set of Kleinian once punctured torus
groups as defined in section~\ref{sec:normalization},  modulo
conjugation in $PSL(2,\CC)$. By abuse of notation, we also denote
by $AH(\torus)$ the set of hyperbolic $3$-manifolds $\{M =
\HH^3/H: [H] \in AH(\torus)\}$, where $[H]$ is the conjugacy class
of $H$ in $PSL(2,\CC)$.

Clearly, whether or not a lamination is realized is a conjugacy
invariant.  Simple closed geodesics are always realized in any
hyperbolic $3$-manifold $M  \in AH(\torus)$, and are dense in the
set of realizable laminations,~\cite{CEpG} theorem 5.3.11. Since
length is a conjugacy invariant,  the above definition of
lamination length $l_{\mu}(M )$ extends by continuity to any $M
\in AH(\torus)$ containing a realization of $|\mu|$. If $|\mu|$ is
connected and not realized in $M$, set $l_{\mu}(M ) = 0$.  (If
$|\mu|$ is not connected one has to be more careful with this
definition since some components of $\mu$ may be realized and
others not; for example on a general surface $|\mu|$ might consist
of disjoint loops some but not all of whose components are
accidentally parabolic. In the case of a punctured torus all
laminations are connected and this difficulty does not occur.) In
the next section, we shall make important use of the following
result.

\begin{prop} The function $L:AH(\torus) \times ML \to \RR$, $L(H,\mu) =
l_{\mu}(\HH^3/H )$ is continuous. \end{prop} \begin{proof} This
result was  asserted by Thurston in \cite{ThuH}; detailed proofs
appear in~\cite{Ohshika} Lemma 4.2 and \cite{Brock} Theorem~5.1.
We remark that the proof in~\cite{Ohshika} seems to have
overlooked the above mentioned difficulties when $|\mu|$ is not
connected. See \cite{Brock} section 7 for a discussion of the
general case. \end{proof}

Note that if a lamination $\mu \in ML$ is realized in $M \in
AH(\torus)$, then the length of $\mu$ in $M$ is equal to the
hyperbolic length of $\mu$ on the surface $\Sigma$, where $\psi:
\Sigma \to M$ is the pleated surface map realizing $|\mu|$, and so
is strictly positive.

In general, the lamination lengths $l_{\mu}(\bch )$ on $\bch$ and
$l_{\mu}(M )$ in $M$ are not the same, and we shall take care to
indicate which length we mean.  In the special case in which $q
\in \Pr{\mu}^+$, however, the lengths $l_{\mu}(\bch^+)$ and
$l_{\mu}(M)$ coincide, and may be safely denoted by $l_{\mu}=
l_{\mu}(q)$.  This is the situation we are discussing in
theorem~\ref{thm:thmq} below.

In section~\ref{sec:cxlenlam}, we shall show how to extend the
holomorphic multiplier $\lambda_{V}$ to a holomorphic function
called the {\em complex length} $\lambda_{\mu}$ of $\mu$ on $\QF$.
Again by linearity and continuity, we have $l_{\mu}(M )= \Re
\lambda_{\mu}$. We also prove in section~\ref{sec:cxlenlam} that
$q \in \Pr{\mu}^+$ implies $\lambda_{\mu} \in \RR$.  Combining
these observations gives that $q \in \Pr{\mu}^+$ implies
$\lambda_{\mu} = l_{\mu}(\bch^+)=l_{\mu}(M)$.

\section{The limit pleating  theorem} \label{sec:closure}

Classically, the ending invariants of a quasifuchsian group are
the marked conformal structures $\omega^{\pm}(q)$ of the tori
$\Omega^{\pm}(q)/G(q)$ and so are points in the Teichm\"uller
space $\teich$.  Suppose we have a sequence $q_n \in \QF$ with
$\omega^{\pm}(q_n) \to \omega^{\pm} \in \teich$.  It then follows
from Bers' simultaneous uniformization theorem that the groups
$G(q_n)$ have an algebraic limit in $ \QF$.  If both of the
sequences $\omega^{\pm}(q_n)$ converge to distinct points in the
Thurston boundary of $\teich$, then Thurston's double limit
theorem~\cite{ThuH} again asserts the existence of an algebraic
limit $G_{\infty}$; the intermediate situation works in a similar
way and is discussed in~\cite{MinskyPT}.

We need an analogous result which asserts the existence of a limit
group when our pleating invariants converge. We also need to
understand the behavior of the pleating invariants when an
algebraic limit exists.  The results we need are collected in the
following {\em limit pleating theorem}, which will be a key factor
in the proof of our main results in section~\ref{sec:Pleating
Rays}.

\begin{thm}{\bf Limit Pleating Theorem.} \label{thm:thmq} Let $\mu,\nu
\in ML$, $[\mu] \neq [\nu]$ and suppose that  $\{q_n\} \in
\Pr{\mu,\nu}$. Then \begin{enumerate} \item  if $l_{\mu}(q_n) \to
\c \ge 0$ and $ l_{\nu}(q_n) \to \d \ge 0$, then there is a
subsequence of the groups $\{ G(q_n)\}$ with an  algebraic limit
$G_{\infty};$ \item  if the  sequence $\{G(q_{n})\}$  has
algebraic limit $G_{\infty}$, then the sequences $\{l_{\mu}(q_{n})
\}$ and $\{l_{\nu}(q_{n})\}$ have finite limits $c \ge 0, d \ge 0$
respectively.  The group $G_{\infty}$ represents a point in $\QF$
if and only if $c >0$ and  $d >0$.
\end{enumerate} \end{thm}

We remark that in the case of a more general surface, the second
statement as it stands is false, as is seen by taking $|\mu|$ to
be a multiple loop such that one, but not all of its components,
becomes accidentally parabolic. It works in our case because any
measured lamination  on $\torus$ is automatically connected. The
result is closely related to, but not the same as, the `Lemme de
fermeture' in~\cite{BonO}, which concerns the existence of the
limit  groups under hypotheses on the limits of bending measures
as opposed to lengths.

The first statement, the existence of the algebraic limit, follows
from a deep estimate of Thurston's about lengths of geodesics in
hyperbolic $3$-manifolds, \cite{ThuH}, theorem~3.3 (Efficiency of
pleated surfaces). The same estimate is fundamental in Thurston's
proof of the double limit theorem in~\cite{ThuH}.  A detailed
discussion and proof of Thurston's estimate is to be found
in~\cite{Canary}, where a limit theorem similar to our first
statement in the context of Schottky groups is proved.

To prove  the second statement we use continuity of lamination
length described in section~\ref{sec:length} above. This allows us
to deduce that the laminations $\mu,\nu$ must be realized in the
{\em algebraic} limit. We complete the proof by showing that the
pleated surfaces which realize $\mu$ and $\nu$ are in fact
components of the convex hull boundary of the algebraic limit.
This idea is in essence the same as that used in~\cite{BonO}, and
we would like to thank F. Bonahon for suggesting this approach.

The statement, and the theorem on continuity of lamination length,
conceals much subtlety. The hypothesis that $G_n \in \Pr{\mu,\nu}$
is crucial;  examples like the one described in~\cite{KerckTh}
show that it is not enough just to  require that some fixed curve
on $\bch^+$ have bounded length. Again, if one takes a varying
sequence $\mu_n \to \mu$ as in~\cite{BonO}, then it is essential
to add the hypothesis that the laminations converge in the
Hausdorff topology as well as in measure, otherwise examples
similar to the one in~\cite{KerckTh} again show that the
convergence may not be strong.

\begin{proof} First we suppose that $l_{\mu}(q_n) \to \c \ge 0$ and $
l_{\nu}(q_n) \to \d \ge 0$, and show that there is some
subsequence of $\{q_n\}$, along which an algebraic limit exists.
Choose and fix an ideal triangulation $\lambda$ on $\torus$;
specifically, take $\lambda$ as the lines from the cusp to itself
in the homotopy classes of the curves $\alpha,\beta$ and
$\alpha\beta$, where $<\pi_1(\torus);\alpha,\beta>$ corresponds to
$<G; A,B>$.

Let $M_n = \HH^3/G_n$ and realize $\lambda$ as the pleating locus
of a pleated surface $S_{n}$ in $M_n$.  The lamination $\lambda$
has no closed leaves and its complement is a pair of ideal
triangles.  Pick $\xi \in ML$. When an oriented arc on a leaf
$|\xi|$ cuts two consecutive sides of one of these complementary
triangles $T$, the two sides meet in an ideal vertex which is
either to its left or its right. The arc of leaf containing an
intersection point $P$ of $|\xi|$ and $\lambda$ goes from one
triangle $T_1$ to another $T_2$. Following Thurston, \cite{ThuH},
we call $P$ a {\em boundary intersection} if the right-left
location of the ideal vertex switches as we cross from $T_1$ to
$T_2$, and we define the {\em alternation number} $a(\xi,\lambda)$
as the total $\xi$-measure of the set of boundary intersection
points. Recall from section~\ref{sec:length} that $l_{\xi}(S_n)$
denotes the length of the lamination $\xi$ measured in the flat
structure of $S_n$ and $l_{\xi} (M_n)$ denotes the length of the
lamination $\xi$ in $M_n$. Then by \cite{ThuH}~theorem~3.3, there
exists a constant $C > 0$, depending only on a fixed choice of
structure for $\torus$, such that
$$l_{\xi} (S_n) \leq l_{\xi} (M_n) + Ca(\xi,\lambda). $$ (We remark that
since $a(\xi,\lambda) \le i(\xi,\lambda)$ the usual intersection
number would be just as good a bound in the present case.)
Applying this inequality in our case to the pleating laminations
$|\mu|$ and $|\nu|$ we find, $$l_{\mu} (S_n) \leq l_{\mu}(q_n) +
Ca(\mu,\lambda), \ \ l_{\nu} (S_n) \leq l_{\nu}(q_n) +
Ca(\nu,\lambda). $$ It follows that the sequences $\{l_{\mu}(S_n)
\}$ and $\{l_{\nu}(S_n) \}$ are bounded.

Since $[\mu] \neq [\nu]$, the laminations $|\mu|, |\nu|$ fill up
$\torus$ and we conclude from \cite{ThuH} proposition~2.4 that the
hyperbolic structures of the surfaces ${S_n}$ lie in a bounded
subset of $\F$ and thus that the lengths $l_{\alpha}(S_n)$ and
$l_{\beta}(S_n)$ of the geodesic representatives of the marking
curves $\alpha$ and $\beta$ on $S_n$ are bounded.  From the
discussion in section~\ref{sec:length}, we conclude that, since
$l_{\alpha} (M_n) \leq l_{\alpha} (S_n)$ and $l_{\beta} (M_n) \leq
l_{\beta} (S_n)$, the sequences $\{|\tr A_n|\},\{|\tr B_n|\} $ are
also bounded.  Therefore we can find a convergent subsequence
along which $\tr A_n$ and $\tr B_n $ converge and thus, (because
from the Markov identity $\tr A$ and $\tr B $ determine at most
two normalized punctured torus groups up to conjugation) we
conclude that a subsequence of $\{G_n\}$ has an algebraic limit
$G_{\infty}$.  This proves statement~1.

\medskip

Now suppose that $G_{\infty}$ is the algebraic limit of a sequence
$G_n = G(q_n) \in \Pr{\mu,\nu}$.  By the continuity of lamination
length on $AH(\torus)$, the sequences $\{l_{\mu}(q_n)\
\{l_{\nu}(q_n)\}$ converge to $\{l_{\mu}(G_{\infty})\},
\{l_{\nu}(G_{\infty})\}$, in particular the limits exist.  We have
to prove that that $G_{\infty} \in \QF$ if and only if both limits
are non-zero.  We note immediately that if $G_{\infty} \in
\QF$,
then, using our assumption that $q_n \in \Pr{\mu,\nu}$, we have
$\{l_{\mu}(q_{n}) \} \to c>0$ and $\{l_{\nu}(q_{n})\} \to d>0$ by
the continuity theorem~\ref{thm:plcont}. This can also be seen
from the fact that all laminations, in particular $\mu$ and $\nu$,
are realized in $G_{\infty}$, see~\cite{ThuN}, \cite{CEpG} theorem
5.3.11.

Suppose that one of the laminations $\mu$ or $\nu$, for
definiteness say $\mu$,  is not realized in $G_{\infty}$.  Since
$|\mu|$ is connected, $l_{\mu}(G_{\infty}) = 0$  and by the
continuity of lamination length on $AH(\torus)$ we deduce that
$c=0$.  Thus we need only prove that if $\mu,\nu$ are both
realized in $G_{\infty}$, and if $c>0,d>0$, then $G_{\infty} \in
\QF$.

\smallskip

Our strategy is to show that the pleated surfaces  which realize
$|\mu|$ and $|\nu|$ are in fact invariant components of
$\bch(G_{\infty})$ which face simply connected invariant
components of the regular set $\Omega(G_{\infty})$. The key point
is to show that if $|\mu|$ is realized in the algebraic limit
$M_{\infty}= \HH^3/G_{\infty}$, then the lift of any leaf of
$|\mu|$ to $\HH^3$ is the limit of corresponding lifts of leaves
of $|\mu|$ in their realizations in $M_n$. To see this we use a
criterion for algebraic convergence to be found in
\cite{McMre} Sec 3.1, see also~\cite{BP} Theorem E.1.13 and~\cite{JorM}
Prop.3.8:

{\em A sequence of groups  $G_n \to G_{\infty}$ algebraically if
and only if there are smooth marking preserving homotopy
equivalences $q_n: M_{\infty} \to M_n$ such that on any compact
subset of $M_{\infty}$, $q_n$ is $C^{\infty}$ close to a local
isometry for all large enough $n$.}

\smallskip

We also have to be careful about markings. Our normalizations are
fixed in such a way that $\rho_n (g)  \to \rho_{\infty}(g)$ for
each $g \in \pi_1(\torus)$.  Let  $\Gamma_0$  be a fixed Fuchsian
group acting on $\DD$, and choose a fixed normalized
representation $\rho_0 : \pi_1(\torus) \to \Gamma_0$ of the {\em
marked} torus  $\torus$.  The action of $G_n=
\rho_n(\pi_1(\torus))$ on $\bch^+_n$ pulls back to the action of a
correspondingly normalized Fuchsian group    $\Gamma_n$  on $\DD$.
This induces a pleated surface map $f_n: \DD \to \HH^3$ with image
$\bch^+_n $, intertwining the action of $\Gamma_n$ on $\DD$ and
$G_n$ on  $\bch^+_n$. Let $h_n: \DD \to \DD$ denote the
homeomorphism which intertwines the actions of  $\Gamma_0$ and
$\Gamma_n$, so that $ f_n   h_n: \DD \to \HH^3$ induces the
representation $\rho_n$. Since $|\mu|$ is realized in
$M_{\infty}$,  there is also a marked Fuchsian group
$\Gamma_{\infty}$,  and a pleated surface  $f: \DD \to \HH^3$
intertwining the actions of $\Gamma_{\infty}$ and $G_{\infty}$
with pleating locus containing $|\mu|$, together with a
homeomorphism $h:\DD \to \DD$ intertwining the actions of $G_0$
and $\Gamma_{\infty}$ such that
 $ fh$ induces $\rho_{\infty}$.
 In this setup,  McMullen's marking preserving homotopy equivalence
$q_n$ lifts to a map  $\tilde q_n: \HH^3 \to \HH^3$ such that
$f_nh_n = \tilde q_n fh$.

Let $l$ be the lift to
$\DD$ of some leaf  of $\mu$ in $\torus$, and suppose that, in the
structure induced  by $\Gamma_0$, it has endpoints  $l^{\pm}$ on
$\partial \DD$.
The corresponding leaves for the structures induced  by
$\Gamma_n,\Gamma_{\infty}$ are
 the geodesics $l_n$, $l_{\infty}$ whose  endpoints on $\partial \DD$
are $h_n(l^{\pm}), h(l^{\pm})$ respectively. From the definition
of pleated surfaces, under $f_n$ and $f$  these leaves are mapped
to geodesics in $\HH^3$.  To make precise the statement  that
leaves of $|\mu|$ in $M_{\infty} = \HH^3/G_{\infty}$ are close to
leaves of the corresponding realizations in $M_n$, we shall prove
that $f_n (l_n) \to   f(l_{\infty})$.

Since  the projection of $ f(l_{\infty})$ to $M_{\infty}$ is
carried on a train track  (see \cite{BonV3} or \cite{Brock} lemma
5.2), it follows that any small neighborhood of the projection is
contained in some compact subset of $M_{\infty}$. Fix an origin $O
\in \HH^3$ and let $x \in  f(l_{\infty})$ be the point nearest
$O$. For any $\epsilon>0, L>0$ we can find $g \in \pi_1(\torus)$
such that the axis of $\rho_{\infty}(g)$ is within $\epsilon$ of
$l$ for a distance $L$ on either side of $x$. The projection of
this long segment $\sigma$ of $\ax \rho_{\infty}(g)$ is contained
in some compact set $V$ in $M_{\infty}$. The restriction to $V$ of
the map $\tilde q_n$ is close to a local isometry for large $n$.
Now the image of a geodesic arc under a map which is $C^{\infty}$
near a local isometry is clearly a quasi-geodesic with small
constants. Thus the images $\tilde q_n(f(l_{\infty}))$ and $\tilde
q_n(\sigma)$ are close to each other and to their corresponding
geodesic representatives. Since $f_nh_n = \tilde q_n fh$, the
geodesic representative of $\tilde q_n(\ax \rho_{\infty}(g))$ is $
\ax \rho_{n}(g)$ and in addition,  since the endpoints of
$l_{\infty}$  are  $h(l^{\pm})$ and 
the endpoints of $l_{n}$ are $h_n(l^{\pm})$,   the  geodesic representative of $\tilde
q_n(f(l_{\infty}))$ is $f_n(l_{n}) $. From the algebraic
convergence, $ \ax \rho_{n}(g)$ is close to $ \ax
\rho_{\infty}(g)$ for all  sufficiently large $n$. Putting this together we see that 
 $f_n(l_n)$ is close to $ f(l_{\infty})$ as required. 

\smallskip
  We now use this fact to prove that the image of the pleated
surface $ f:\DD  \to \HH^3$  is a component of the convex hull boundary of $G_{\infty}$. The projection of the pleating locus of $f$ to $\DD/\Gamma_{\infty}$ is a geodesic lamination
$\hat \mu$ which contains $|\mu|$; we also use $\mu$  to denote the lift  to $\DD$. If the pleating locus actually
equals $|\mu|$,  add an extra leaf to make a maximal lamination
$\hat \mu$. Otherwise let $\hat \mu$ be the pleating locus of $f$.
In either case, by area considerations, $\hat \mu$ contains
exactly  one extra leaf, one end of which goes out to the cusp and
the other end of which spirals onto boundary leaves of $|\mu|$.
Notice that since the pleating locus of $f_n$ actually equals
$|\mu|$ (since the pleating locus of  the convex hull boundary
cannot contain any leaf going out to the cusp) the additional leaf
of $\hat \mu - |\mu|$ is necessarily mapped to a geodesic by
$f_n$. Moreover the endpoint of the additional leaf is a cusp and
hence any lift moves continuously as $n \to \infty$.

\smallskip

We call any ideal triangle in $\HH^3$ formed by the lifts of the
images of the boundary leaves of a complementary region of $\hat
\mu$ under a pleated map  a {\it plaque}. The vertices of such a
triangle are either the endpoints of leaves of the lamination or
parabolic fixed points. For clarity,  denote the images in $\HH^3$  of $\hat \mu$ under the pleated surface maps $f_n, f$ by
$\hat \mu_n$, $\hat \mu_{\infty}$ respectively.   We have just
shown  that any plaque of $\hat \mu_{\infty}$  is arbitrarily
closely approximated
in $\HH^3$ by a plaque of $\hat \mu_n$ for
all sufficiently large $n$. Notice also that any plaque of $\hat
\mu_n$ is contained in a support plane for $\bch^+_n$.

\smallskip Denote the image of $f$ by $\Pi^{+}$. We want to show that
$\Pi^{+} $ is a component of $\bch(G_{\infty})$. Let $X$ be a
plane containing a plaque of $|\mu_{\infty}|$ and let $X_n$ be a
sequence of planes containing approximating plaques for
$\bch^+(G_n)$.  We claim that all of $\Pi^+$ lies on the same side
of $X$ so that $X$ is a support plane for $\Pi^+$. If not, we can
find points  $y,y' \in \Pi^+$ on opposite sides of $X$ so that the
geodesic joining $y$ to $y'$ crosses $X$ transversally. By
choosing $n$ sufficiently large, we can find $y_n,y'_n$ near to
$y,y'$ in $\bch^+_n$, and a support plane $X_n$ to $\bch^+_n$
close to $X$, such that the geodesic from $y_n$ to $y'_n$ crosses
$X_n$, which is impossible.

Denote by $H_X$ the closed half space bounded by $X$ containing
$\Pi^+$ and set $K=\cap_X H_X$ where $X$ runs through all planes
containing plaques of $\Pi^+$. By the above, $\Pi^+ \subset K$ so
$K \neq \emptyset$. By its construction, $K$ is convex and closed.
Moreover $K$ is  $G_{\infty}$ invariant since the same is true of
$\Pi^+$. Let $g \in \Gamma_{0}$.  Pick $y \in K$; then
$\rho_{\infty}(g)^{\pm m}(y) \in K$ for  $m=1,2,\ldots$.  By
convexity $K$  contains the geodesic joining
$\rho_{\infty}(g^{-m})(y)$ and $\rho_{\infty}(g^{m})(y)$ for all
$m$;  by closure, it contains the axis of $\rho_{\infty}(g)$. The
axes of elements of $G_{\infty}$ are dense in the geodesics
joining all its limit points and we conclude $\C(G_{\infty})
\subset K$.

We claim $\Pi^+ \subset \bch(G_{\infty})$. Let $P$ be a plaque of
$\Pi^+$. Clearly $P \subset {\cal C}(G_{\infty})$ and so  $P
\subset K$. Since $P$ is by definition contained in a support
plane for $K$, we conclude  $P \subset \bch(G_{\infty})$. Since
$\Pi^+$ is the closure in $\HH^3$ of the union of its plaques, the
claim follows.

We prove in lemma~\ref{lem:embedded} below that $\Pi^+$ is
embedded in $\HH^3$. (This rules out the possibility that, for
example, $|\mu|$ is rational and the bending angle is $\pi$.)
Thus $\Pi^+$ is isometric to a complete hyperbolic surface and
hence is both open and closed in $\bch(G_{\infty})$. Since $\Pi^+$
is connected, it must be  a component of $\bch(G_{\infty})$. As
such, it faces a component $\Omega^+$ of $\Omega(G_{\infty})$.
Moreover since $\Pi^+$ is simply connected, so is $\Omega^+$.
(This also follows from the fact that the limit representation $
\rho: \pi(\torus) \to G_{\infty}$ is faithful.) Also $ G_{\infty}$
invariance of $\Omega^+$ follows from that of $\Pi^+$.

Now there is a similar image $\Pi^{-}$ for the pleated surface map
which realizes $|\nu|$, from which we deduce the existence of
another simply connected invariant component $\Omega^-$ of
$\Omega(G_{\infty})$. We conclude, see for example~\cite{Marden}
lemma 3.2, that $G_{\infty} \in \QF$. \end{proof}

\begin{lemma} \label{lem:embedded} With the notation and conditions
above, the image $\Pi^{+}$ of the pleated surface map $f$ is
embedded in $\HH^3$. \end{lemma}
 \begin{proof} If $\Pi^{+}$ is not
embedded then $  f(x) =  f(y)$ for some distinct points $x,y \in
\DD$; these cannot be in the same plaque since $f$ is an isometry
on plaques. We begin by reducing to the case in which $x$ and $y$
are both contained in leaves of $\hat \mu$. If not, suppose that
$x$ is in a complementary region of $\hat \mu$, and let $P_x$ be
the image plaque containing $ f(x)$. Now $y$ is either in a
distinct complementary region with image plaque $P_y$, or on a
leaf with image  a geodesic $L$. If  $P_y$ or $L$ cuts $P_x$
transversally, then the same is true for all nearby pleated
surfaces $f_n$, since the endpoints which determine plaques and
leaves move continuously. This is impossible since $f_n(\DD) =
\bch^+_n$ is embedded. Thus $P_y$ (or $L$) and $P_x$ are in a
common plane. In the first case there is some point on boundary
leaves of both $P_y$ and $P_x$, and in the second $L$ meets some
boundary point of $P_x$.

Now we use the uniform injectivity theorem~\cite{ThuH1, MinskyH}:
for any $\epsilon>0$, and for any $M \in AH(\torus)$, there exists
$\delta>0$ such that for any pleated surface  $F:S \to M$, if the
distance in $S$ between points $u$ and $v$ in leaves $l_u$ and $
l_v$ of the pleating lamination is at least $\epsilon>0$, and if
$u$ and $v$ are in the thick part of $S$, then the distance in the
projective unit tangent bundle to $M$ between the unit tangent
vectors to  $F(  l_u)$ and $F(  l_v)$ at $F(u)$ and $F(v)$ is at
least $\delta$.

Now it is a standard argument that because the non-cuspidal part
${ \DD / \Gamma_{\infty}}^{{\rm nc}}$ of $ \DD / \Gamma_{\infty}$
has bounded diameter,  the  injectivity radius of $M_{\infty}$  is
bounded below in a neighborhood of $\bar f(\DD /
\Gamma_{\infty})$, where  $\bar f $  denotes the induced map on
quotients.) For otherwise $\bar f({ \DD / \Gamma_{\infty}}^{{\rm
nc}})$ would contain loops corresponding to pairs of non-commuting
loxodromics contained in a Margulis tube of $M_{\infty}$ which is
impossible, see~\cite{MinskyH} section 2.1. Choose $\epsilon$ less
than this injectivity radius.  Then if  $ f(x)= f(y)$,  the
distance between  $x$ and $y$ (in $ \DD$) must be at least
$\epsilon$. Let $l_x,l_y$ be leaves of the lift of $\hat \mu$ to
$\DD$ through $x,y$ respectively. It follows that the image leaves
$ f (l_x)$ and $  f (l_y)$ meet at a definite angle in $\HH^3$.

Consider the plane $P$ containing these two leaves. It meets
$\Chat$ in a circle $C$. Notice that any circle through the
endpoints of $f(l_x)$ other than $C$ separates the endpoints of
$f(l_y)$. Now for any nearby group $G_n$, there are leaves
$f_n(l_x),f_n(l_y)$ near $f(l_x), f(l_y)$. Any support plane to
$\bch^+_n$ through either of these leaves meets $\Chat$ in a
circle which cannot  separate the other pair of endpoints. One
deduces easily that any pair of support planes for $\bch^+_n$ must
meet $\Chat$ in circles both of which are close to $C$, and that
$\Lambda(G_n)$ is contained in the thin ring or crescent between
them. It follows that every support plane of $\bch_n^-$ has very
small diameter, and hence that the distance of any such support
plane to $  f(x)$ tends to  $\infty$ with $n$. On the other hand,
any support plane for $\bch_n^-$ contains points  close to some
plaque of the pleated surface which realizes $|\nu|$ in
$M_{\infty}$. Pick a point $z \in \HH^3$ on a lift of a leaf of
$|\nu|$, at distance $D$ say from $ f(x)$. Since $z$ is on a
plaque of $|\nu|$ it is close to a support plane of $\bch_n^-$.
This shows there are points in $\bch_n^-$ which stay at bounded
distance, with bound close to $D$, from $ f(x)$. This
contradiction completes the proof.
\end{proof}

\section{Complex length}
\label{sec:cxlenlam}

In this section we introduce the {\em complex length} of a
measured
lamination.  Just as lamination length as defined in
section~\ref{sec:laminations} is a real analytic function on
$\F$, the
complex lamination length is a holomorphic function on $\QF$.
The
relationship of this holomorphic function to pleating varieties,
in
particular theorem~\ref{thm:realI}, is a central tool in
everything which
follows.  Complex lamination length has also been introduced
using somewhat
different techniques by Bonahon~\cite{Bon}.

\subsection{Complex length of a loxodromic }
\label{sec:cxlenlox}

Let $M \in PSL(2,\CC)$.  Its {\em complex translation length}
$\lambda_M
 \in \CC/2 \pi i \ZZ$ is given by the equation
\begin{equation}
\label{eqn:trmult} {\pm} \tr M = 2 \cosh {\lambda_M/2}
\end{equation}
where $\tr M$ is the trace of $M$ and we
choose the sign so that $\Re \lambda_M \ge 0$.

Complex length is invariant under conjugation by M\"obius
transformations
and has the following geometric interpretation, provided $M$ is
not
parabolic.  Let $x \in \ax M$ and let $\bar v$ be a vector
normal to $\ax
M$ at $x$.  Then $\Re \lambda_M$ is the hyperbolic distance
between $x$ and
$M(x)$ and $ \Im \lambda_M$ is the angle $\mmod 2 \pi$ between
$M(\bar v)$
and the parallel transport of $\bar v$ to $M(x)$, measured
facing the
attracting fixed point $M^+$ of $M$.  In particular, if $M$ is
loxodromic
then $ \Re \lambda_M > 0$ and if $M$ is purely hyperbolic then
in addition
$ \Im \lambda_M \in 2\pi\ZZ$; equivalently $\tr M \in \RR, |\tr
M| > 2$.
(We refer to \cite{PS} for a detailed discussion of the sign
ambiguity in
equation~\ref{eqn:trmult}; note that in our notation here
$\lambda_M$ is
{\em twice} the multiplier denoted by $\lambda_M$ in \cite{PS}.)

Let $q \in \QF$, let $\gamma \in \S$ and denote the element
representing
$\gamma$ in the group $G(q)$ by $W(q)$.  Because the trace is a
conjugation
invariant, the complex translation length $\lambda_W(q)$ depends
only on
$q$ and is independent of the normalization of $G(q)$.  We want
to define
the complex length $\lambda_{\gamma}(q)=\lambda_W(q)$ as a
holomorphic
function on $\QF$ with values in $\CC$, not $\CC/2 \pi i \ZZ$.
To do this,
we choose the branch that is real valued on $\F$.  Since
$\lambda_{\gamma}
\neq 0$ on $\QF$ this choice uniquely determines a holomorphic
function
$\lambda_{\gamma}:\QF \rightarrow \CC$.  From now on, the term
``complex
length'' will always refer to this branch.

 We define the {\em complex length} of
 the rational lamination $\mu = c \delta_{\gamma} \in ML_Q$, $
c>0$, as
$\lambda_{\mu}(q) = c \lambda_{\gamma}(q)$.

To define the complex length $\lambda_{\mu}(q)$ for arbitrary
$\mu \in ML$
and $q \in \QF$, we would like to choose $\mu_n \in ML_Q$,
$\mu_n \to \mu$
and set  $$\lambda_{\mu}(q) = \lim_{n \to \infty}
\lambda_{\mu_n}(q).$$
To justify this, we need to show these limits exist
 and are independent of the sequence
$\{\mu_n\}$.

We do this using the following theorem which summarizes the
results of
\cite{KerckEA}, lemma~2.4 and \cite{KerckN}, theorem~1.  In the
statement,
$l_{\mu}$ denotes lamination length defined in
section~\ref{sec:laminations}.

\begin{thm}
\label{thm:Kerck2}
The function $( c\delta_{\gamma},p) \mapsto c
l_{\delta_{\gamma}}(p)$ from
$ML_Q \times \F$ to $\RR^+$ extends to a continuous function
$(\mu,p)
\mapsto l_{\mu}(p)$ from $ML \times \F $ to $ \RR^+$.  If $\mu_n
\in ML_Q$,
$\mu_n \to \mu$, and $p \in \F$, then $ l_{\mu_n}(p) \to
l_{\mu}(p)$.  The
limit functions $ p \mapsto l_{\mu}(p)$ are non-constant, and
the limit is
uniform on compact subsets of $\F$.
\end{thm}

We also need an elementary lemma about holomorphic functions.
\begin{lemma}
\label{lemma:elem} If $f:\QF \rightarrow \CC$ is holomorphic and
if $f \equiv c$ on $\F$ for some constant $c$, then $f \equiv c$
on $\QF$.
\end{lemma}

\begin{proof}
Because $\F$ is the $\RR^2$-locus of the complex Fenchel-Nielsen
coordinates $(\lambda,\tau)$ in $\QF$, see \cite{Christos}, and
section~\ref{sec:twistsandqbs} below, the conclusion follows
directly from
the Cauchy-Riemann equations applied to each variable
separately.
\end{proof}

\begin{thm}
\label{thm:welldefined}
The function $(\mu,q) \mapsto
\lambda_{\mu}(q)$ from $ ML_Q\times  \QF$ to $ \CC$   extends to
a  continuous function from $ML \times \QF$ to
$\CC$, also denoted $\lambda_{\mu}(q)$.  The function
$q \mapsto \lambda_{\mu}(q)$ is holomorphic and
non-constant for all $\mu$ and the family $\{\lambda_{\mu}\}$ is
bounded and equicontinuous on compact subsets of $\QF$.
\end{thm}

\begin{proof}
By construction, the functions $\{\lambda_{\mu}\}$, $\mu \in
ML_Q$, omit
the half plane $\Re z < 0$ and thus form a normal family on
compact subsets
of $\QF$.  It follows that if $\mu_n \to \mu$, $\mu_n \in ML_Q$,
then
suitable subsequences of $\{\lambda_{\mu_n}\}$ converge to limit
functions
that are holomorphic.

We note that on $\F$, if $\mu \in ML_Q$, then $\lambda_{\mu}$ is
real and
coincides with $l_{\mu}$.  By theorem~\ref{thm:Kerck2}, if
$\mu_n \to \mu$,
$\mu_n \in ML_Q$, then $\{l_{\mu_n}\}$ is uniformly convergent
on compact
subsets of $\F$; further, the limit function $l_{\mu}$ is
finite,
non-constant and independent of the choice of the sequence
$\{\mu_n\}$.
The result now follows from lemma~\ref{lemma:elem}.  \end{proof}

 For $\mu \in ML$, we call $\lambda_{\mu}$ the {\em complex
length} of
$\mu$.  Throughout this paper, the complex length functions are
a
fundamental tool.

We remark that
\begin{enumerate}
\item Suppose $q \in \QF$ and let $F^{\pm}(q) \in \F$ denote the
flat
structures (see section~\ref{sec:pllocus}) on the convex core
boundary
$\bch^{\pm}(q)/G(q)$.  If $\mu \in [pl^{\pm}(q)]$, then
$l_{\mu}(F^{\pm}(q)) = \Re \lambda_{\mu}(q)$.
\item
For $q \in \QF$, $\mu \in ML$, $\Re \lambda_{\mu}$ coincides
with
the
lamination length $l_{\mu} (M(q))$ in the 3-manifold $M(q)
=\H3/G(q)$ as
discussed in section~\ref{sec:length} above.
\end{enumerate}

For $\mu \in ML_Q$, $\Re \lambda_{\mu}(q) = l_{\mu}(M(q))$, so
by
continuity, both statements hold for all $\mu \in ML$.

\subsection{Complex length and Pleating varieties}
\label{sec:clandxrs}

The first step in proving our main theorems
 is to show that for any $\mu \in ML$,
the complex length $\lambda_{\mu}$ is real valued on $\Pr{\mu}$.

First
consider the case $\mu \in ML_Q$.  We have

\begin{prop}
\label{thm:realtrace} Suppose $\mu \in ML_Q$.  Then $\Pr{\mu}
\subset
\lambda_{\mu}^{-1}(\RR^+)$.
\end{prop}

\begin{proof}  This is just a reformulation of the easy
observation, proved in \cite{KStop}, lemma~4.6, that if a
geodesic $\gamma$
 is contained in $|pl^{\pm}(G)|$, then any representative in $G$
is purely
 hyperbolic.
\end{proof}

We now extend
proposition~\ref{thm:realtrace} to arbitrary laminations.

\begin{thm}
\label{thm:realI} Let $q \in \QF$ and suppose $ pl^+(q)=\mu $.
Then $\lambda_{\mu}(q) \in \RR$.
\end{thm}

\begin{proof}
For $[\mu] \in ML_Q$ this is proposition~\ref{thm:realtrace}, so
suppose
$[\mu] \notin ML_Q$.

The map $\QF - \F \rightarrow ML \times \F$ that takes $q \in
\QF - \F$ to
$( pl^+(q), F^+(q))$ where $F^+(q)$ is the flat structure of
$\bch^+/G(q)$
is continuous by theorem~\ref{thm:plcont}.  The map is also
injective
because the hyperbolic structure $F^+$ together with the bending
data
$pl^+$ determine the group $G = G(q)$.  Let $U \subset \QF -\F$
be an open
ball containing $q$; if $[pl^+(q')]$ were constant on $ U$, a
four
dimensional neighborhood would have a three dimensional image,
violating
the invariance of domain for a continuous injective map.

 By the continuity theorem~\ref{thm:plcont}, since $PML$ is one
dimensional, we may find a sequence $q_n \to q$ in $U$ such that
$pl^+(q_n)=\mu_n$ with $\mu_n \in ML_Q$.  By
proposition~\ref{thm:realtrace}, $\lambda_{\mu_n}(q_n) \in \RR$.
 By the
continuity theorem again, $\mu_n \to \mu $ and hence
$\lambda_{\mu_n} \to
\lambda_{\mu}$ uniformly on compact subsets of $\QF$.  Thus
taking a
diagonal limit we have $\lambda_{\mu_n}(q_n) \to
\lambda_{\mu}(q)$ and
$\lambda_{\mu}(q) \in \RR$.
\end{proof}

\section{Twists and Quakebends}
\label{sec:twistsandqbs}
In this section we briefly discuss complex Fenchel-Nielsen
coordinates and
quakebends, and the connection with the convex hull boundary
$\bch$.  This
circle of ideas is at the heart of the proof of the local
pleating
theorem~\ref{thm:open} in section~\ref{sec:openness}; some of
the ideas are
also needed in section~\ref{sec:Pleating Rays}, where we work in
quakebend
planes as defined in section~\ref{sec:qbplanes} below.

\subsection{Complex Fenchel Nielsen coordinates}
\label{sec:cfncoords}

 Complex Fenchel Nielsen parameters were introduced in
\cite{Christos,SerTan} (see also \cite{KSbend}) as a
generalization to
$\QF$ of the classical Fenchel Nielsen coordinates for Fuchsian
groups.
Here we briefly summarize the main points as applied to
$\torus$.

Let $\langle G; A,B \rangle$ be a marked quasifuchsian punctured
torus
group constructed from a pair of marked generators $\alpha,
\beta$ of
$\pi_1(\torus)$ as described in~\ref{sec:normalization}.
Complex Fenchel
Nielsen coordinates $(\lambda_A, \tau_{A,B})$ for $\langle G;
A,B \rangle$
are obtained as follows.  The parameter $\lambda_A \in \CC /2\pi
i \ZZ$ is
the complex translation length of the generator
$A=\rho(\alpha)$, or
equivalently the complex length $\lambda_{\alpha}$. The twist
parameter
$\tau_{A,B} \in \CC /2\pi i \ZZ$ measures the complex shear when
the axis
$\ax B^{-1}AB$ is identified with the axis $\ax A$ by $B$.  More
precisely,
if the common perpendicular $\delta$ to $\ax B^{-1}AB$ and $\ax
A$ meets
these axes in points $Y,X$ respectively, then $\Re \tau_{A,B}$
is the {\em
signed} distance from $X$ to $B(Y)$ and $\Im \tau_{A,B}$ is the
angle
between $\delta$ and the parallel translate of $B(\delta)$ along
$\ax A$ to
$X$, measured facing towards the attracting fixed point of $A$.
On the
critical line $\F_{\alpha,\beta}$, $\tau_{A,B} \equiv 0 \mmod
2\pi i $ and
$\ax A, \ax B$ intersect orthogonally.  Thus a point on this
line
corresponds to a rectangular torus with generators $(A,B)$.  The
conventions for measuring the signed distance and the angle are
explained
in more detail in \cite{KSbend} but are not important here.

As shown in~\cite{Fenchel,KSbend,Christos}, given the parameters
$
\lambda_A,\tau_{A,B}$, and a fixed a normalization, one can
explicitly
write down the matrix generators for a marked two generator
group $
G(\lambda_A,\tau_{A,B}) \subset PSL(2,\CC)$ in which the
commutator $[A,B]$
is parabolic.  This group may or may not be discrete.  The
matrix
coefficients of $G$ depend holomorphically on the parameters.
The
construction thus defines a holomorphic embedding of $ \QF$ into
a subset
of $\CC /2\pi i \ZZ \times \CC /2\pi i \ZZ$, in which Fuchsian
space $\F$
is identified with the image of $\RR^2$.

We want to lift this to an embedding into $\CC^2$.  In
section~\ref{sec:cxlenlam} we discussed how to lift the length
function
$\lambda_A$ on $\QF$ to a holomorphic function on $\CC$.  We can
similarly
lift the twist parameter $\tau_{A,B}$ by specifying that it be
real valued
on $\F$.

On $\F$, the real valued parameters $ \lambda_A,\tau_{A,B}$
reduce to the
 classical Fenchel Nielsen parameters $ l_A,t_{A,B}$ defined by
the above
 construction with $\lambda_A$ the hyperbolic translation length
$ l_A$ of
 $A$ and $\tau_{A,B}$ the twist parameter $t_{A,B}$.

Clearly, the complex Fenchel Nielsen construction can be made
relative to
any marking $V,W$ of $G$.  As described in detail in section~5
of
\cite{KSbend}, for fixed $\lambda \in \RR^+$ and $\tau \in \CC$,
the
complex Fenchel Nielsen construction relative to $V,W$
determines a map $
\DD \to \HH^3$.  This map is the composition of the earthquake
$\E_{\gamma}(\Re \tau)$ along the geodesic $\gamma$ represented
by $V$ with
an isometry $\psi:\DD \to \H3$.  The earthquake $\E_{\gamma}(\Re
\tau): \DD
\to \DD$ intertwines the action of the rectangular torus group
$G(\lambda,0)$ with the group $ G(\lambda,\Re \tau)$.  The
isometry $\psi$
is a pleated surface map with with pleating locus $\gamma$ and
angle $\Im
\tau$ between the outward normals to adjacent flat planes.  It
conjugates
the actions of $ G(\lambda,\Re \tau)$ on $\DD$ and $ G(\lambda,
\tau)$ on
its image in $\HH^3$.  We set $ \DD_{\gamma}(\lambda, \tau)=
\psi(\DD) $.
We note for future use that the bending measure of a transversal
$\sigma$
is $i(\gamma, \sigma) \Im \tau$.

 \subsection{Quakebends}
\label{sec:genqbs}

 Quakebends are a complex version of earthquakes.  The
construction was
  introduced by Thurston and is explained in detail in
\cite{EpM} and also
  summarized in \cite{KSbend}.  An alternative discussion can be
found in
  \cite{McMeq}.

Let $p \in \F$ and let $G_0 = G(p)$ act on the disk $\DD \subset
\H3$.  For
$\mu \in ML$ and $\tau \in \CC$, the quakebend construction
defines an
isomorphism $\Q_{\mu}(\tau)$ from $G_0$ to its image
$\Q_{\mu}(\tau)(G_0)
=G^p_{\mu}(\tau)$, together with a pleated surface
$\psi^p_{\mu}(\tau):\DD
\to \H3$ conjugating the actions of $\E_{\mu}( \Re \tau) (G_0) =
G^p_{\mu}(
\Re \tau)$ on $\DD$ and $G^p_{\mu}(\tau)$ on the image
$\DD^p_{\mu}(\tau)=\psi^p_{\mu}(\tau)(\DD)$.  If $\Im \tau \neq
0$, then
$\DD^p_{\mu}(\tau)$ has pleating locus $|\mu|$.  When $\tau=0$,
$\psi^p_{\mu}(\tau)= id$ and $G^p_{\mu}(\tau)=G_0$.  When $\Im
\tau = 0$
and $\Re\tau = t$, $\Q_{\mu}(\tau)$ coincides with the
earthquake
$\E_{\mu}(t)$, $\DD_{\mu}^p(\tau)=\DD$ and $G^p_{\mu}(t)$ is
discrete and
Fuchsian for all $t \in \RR$.  If $\Re \tau =0$, we call the
quakebend a
{\em pure bend}.

If the lamination $\mu$ is rational, $\mu = \k \delta_{\gamma}$,
an
earthquake along $\mu$ reduces to a Fenchel Nielsen twist.  In
terms of
Fenchel Nielsen coordinates $ (l_V,t_{V,W})$ relative to a
marking $(V,W)$,
where $V \in G$ represents the geodesic $\gamma$, this is given
by the
formula $(\lambda_V,\tau_{V,W}) \mapsto G(\lambda_V,\tau_{V,W} +
\k t)$.
Likewise a quakebend along $\k \delta_{\gamma}$ is the complex
Fenchel
Nielsen twist given by the formula
$\Q_{\mu}(\tau):G(\lambda_V,\tau_{V,W})
\mapsto G(\lambda_V,\tau_{V,W} + \k\tau)$.  In particular, if
the base
point $p \in \F$ is the rectangular group $G(\lambda, 0)$
relative to its
marking $(V,W)$, the image pleated surface $\DD^p_{\mu}(\tau)$
is exactly
$\DD_{\gamma}(\lambda, \tau)$ as described
in~\ref{sec:cfncoords} above.
We shall make frequent use of this observation below.  Note that
the
bending measure of a transversal $\sigma$ to $\DD^p_{\mu}(\tau)$
is always
$i(\sigma,\mu) \Im \tau$.

So far, we have only discussed quakebends when the basepoint $p$
is in
$\F$.  Examining \cite{EpM}, however, it is clear that one can
make the
same construction starting from a basepoint $q \in \Pr{\mu}^+$.
More precisely, let $pl^+(q)$ be the bending measure on
$\bch^+(q)$, so
that (by the unique ergodicity of measured laminations on a
punctured
torus) $pl^+(q) = \k \mu$ for some $ \k >0$.  Let the
flat structure of $\bch^+(q)$ be represented by the Fuchsian
group $F^+(q)$
acting in $\DD$. One can define the quakebend $\Q^q_{\mu}(\tau)$
as the
group obtained by the quakebend $\Q^q_{\mu}(\tau + ik)$ acting
on $F^+(q)$;
in other words compose an earthquake along $\mu$ by $\Re \tau $
with a pure
bend by $\Im \tau + \k$.  In this case, we should consider the
time zero
pleated surface $\DD^q_{\mu}(0)$ to be the surface $\bch^+$.
(See also
\cite{Christosbend, McMeq} for other versions of this
construction.)

We shall not need to discuss here the problems associated with
defining a
quakebend from an arbitrary basepoint in $\QF$.

 \subsection{Quakebend planes}
\label{sec:qbplanes}

In what follows, we shall often want to regard the quakebend
parameter
$\tau$ as a holomophic function on the space of representations
$ \rho:
\pi_1(\torus) \rightarrow PSL(2,\CC)$, modulo conjugation in
$PSL(2,\CC)$.
When the basepoint is Fuchsian, this is justified by the
following
proposition, which is~\cite{EpM}, Lemma~3.8.1.

\begin{prop}
\label{prop:Q1}  Let $p \in \F$,  $\tau \in \CC$, $\mu \in ML$,
and let $G^p_{\mu}(\tau) = \Q_{\mu}(\tau)(G(p))$.  Then the
matrix
coefficients of the elements of $G^p_{\mu}(\tau)$ are
holomorphic functions
of $\tau$.
\end{prop}

 It is clear that the Epstein-Marden proof still works when the
basepoint
$q$ is in $\Pr{\mu}^+$.

This result enables us to
 introduce {\em quakebend planes}, which are the  device
used in section~\ref{sec:Pleating Rays}
to reduce the investigation of  pleating varieties to
a tractable problem  in one complex dimension.

For $q \in \Pr{\mu}^+ \cup \F$, we set $\Q_{\mu}^{q} = \{
G^q_{\mu} (\tau)
: \tau \in \CC \}$; we call $\Q_{\mu}^{q}$ the {\em
$\mu$-quakebend plane
based at $q$} and sometimes write $\Q_{\mu}^{q}(\tau) $ for
$G^q_{\mu}
(\tau)$.  By proposition~\ref{prop:Q2} below, a neighborhood of
$q$ in
$\Q_{\mu}^{q}$ is contained in $\Pr{\mu}$ --- but we emphasize
once again
that that in general the whole of $\Q_{\mu}^{q}$ is not
contained in $\QF$
(see proposition~\ref{prop:genbendinglimit} below
and~\cite{McMeq}).

 In the rational case $\mu \in ML_Q$, $\Q_{\mu}^{q}$ has a very
easy
description in terms of complex Fenchel Nielsen coordinates.
Suppose that
$\mu=\delta_{\gamma},\ \gamma \in \S$ and that
$(\gamma,\gamma')$ are a
pair of marked generators for $\pi_1(\torus)$.  Let
$(\lambda_V,\tau_{V,W})
\subset \CC^2$ be complex Fenchel Nielsen coordinates relative
to
corresponding marked generators $(V,W)$ of $G$. Let $\c =
\lambda_{\mu}(q)$.  Then it is clear from the discussion above
that
$\Q_{\mu}^{q}$ is just the slice $ \{(\c,\tau)\} \subset \CC^2$.
 We denote
this slice by $\Q_{\gamma, \c}$. Clearly, $\Q_{\gamma,\c}$ meets
$\F$ along
the earthquake path $\E_{\gamma,\c}$.

More generally, if $\mu \in ML$ and $p,p' \in \E_{\mu,\c}$, it
is clear
that $\Q_{\mu}^{p}= \Q_{\mu}^{p'}$; we denote this plane by
$\Q_{\mu,\c}$.
Clearly, $\Q_{\mu,\c}$ meets $\F$ along the earthquake path
$\E_{\mu,\c}$.
In general, however, if $ q,q' \in \Pr{\mu}$ and $
\lambda_{\mu}(q)=
\lambda_{\mu}(q')$, then it is not immediately clear whether or
not
$\Q_{\mu}^{q}= \Q_{\mu}^{q'}$.  It is a consequence of our main
results
that $ \lambda_{\mu}(q)= \lambda_{\mu}(q')$ always implies
$\Q_{\mu}^{q}=
\Q_{\mu}^{q'}$; this is proved in corollary~\ref{cor:uniqplane}
below.

As explained above, for a basepoint $q \in \Pr{\mu} \cup \F$,
the quakebend
plane $\Q_{\mu}^{q}$ is not, in general, contained in $\QF$.  We
note that
in the special case $ p \in \F$, since $\QF$ is an open
neighborhood of
$\F$ (in the space of representations into $PSL(2,\CC)$ modulo
conjugation), it follows that for small $\tau$,
$G^p_{\mu}(\tau)$ is
quasifuchsian.  The following stronger result shows that, as one
would
naively expect, as one quakebends along $\mu$ away from a
basepoint $q \in
\Pr{\mu}^+ \cup \F$ (for which $\bch^+= \DD^q_{\mu}(0)$), the
pleated
surface $\DD^q_{\mu}(\tau)$ remains equal to $\bch^+$ for all
small $\tau$.

\begin{prop}
\label{prop:Q2}  Given $q \in \Pr{\mu}^+ \cup \F$  and $\mu \in
ML$, there
exists
$\epsilon > 0$, depending on $\mu$ and $q$, such that if
$|\tau|<\epsilon$,
then $G^{q}_{\mu}(\tau) \in \QF$ and
  $\DD^q_{\mu}(\tau)$ is a component of
$\bch(G^q_{\mu}(\tau))$.
\end{prop}

\begin{proof} This is proved in~\cite{KSbend}, prop~8.10 for the
case
in which the basepoint
$q$ is in $\F$.  It is clear that the same proof works
in our more general case.
\end{proof}

We note that if $G^q_{\mu}(\tau) \in \QF$ and
 $\DD^{q}_{\mu}(\tau)=\bch^+(G^q_{\mu}(\tau))$, then the flat
structure of
 $\bch^+(G^q_{\mu}(\tau))$ is represented by the Fuchsian group
 $\E_{(\Re\tau) \mu}(F^+(q))$ obtained by earthquaking a
distance $\Re
 \tau$ along the pull-back of $\mu$ to $\DD$. This observation
will be
 important in section~\ref{sec:openness} below.

\section{The local pleating theorem}
\label{sec:openness}

In this section we prove the local pleating
theorem~\ref{thm:open}.  We
derive various consequences including the density theorem 5 of
the
introduction and a detailed description of how pleating
varieties meet
$\F$.  The statement of the theorem is as follows.

\begin{thm}{\bf Local Pleating Theorem}
\label{thm:open} Suppose that $\nu \in ML$ and $q_0 \in
\Pr{\nu} \cup \F$.  Then there exists a  neighborhood $U$ of $
q_0$
in $\QF$ such that if $q \in U$
and $\lambda_{\nu}(q) \in \RR^+$,  then $q \in \Pr{\nu} \cup
\F$.
\end{thm}

Our starting point for proving this theorem is proposition~7.6
 of~\cite{KSbend}, part of whose content can be stated in the
following
 way.  We write $G^{q_0}_{\gamma}(\tau) $ for
 $G^{q_0}_{\delta_{\gamma}}(\tau)
=\Q_{\delta_{\gamma}}(\tau)(G(q_0))$.

\begin{prop}
\label{prop:7.6}
 Suppose that $\gamma \in \S$ and $q_0 \in
 \F$.  Then there exists $ \eta > 0$ such that  if
$|\tau | < \eta$
and $\lambda_{\gamma}(G^{q_0}_{\gamma}(\tau)) \in \RR^+$,  then
$G^{q_0}_{\gamma}(\tau)  \in \Pr{\gamma} \cup
\F$.
\end{prop}

This proposition can be regarded as the special case
of~\ref{thm:open} in
 which $\nu= \k \delta_{\gamma} \in ML_Q$, the basepoint $q_0$
is Fuchsian
 and we restrict the discussion to the quakebend plane
$\Q^{q_0}_{\gamma}$
 through ${q_0}$.

We begin by reviewing the argument in~\cite{KSbend}.  Suppose
$\gamma \in
\S$, let $V \in G $ represent $\gamma$ and choose $W \in G $
such that
$(V,W)$ is a marking. Let $(\lambda_V,\tau_{V,W})$ be complex
Fenchel
Nielsen coordinates for $ \QF$ relative to $(V,W)$; thus we
regard
$(\lambda_V,\tau_{V,W})$ as holomorphic functions on $\QF$.  As
described
in section~\ref{sec:cfncoords}, whenever $\lambda_V =
\lambda_V(q) \in
\RR^+$, the complex Fenchel Nielsen construction determines a
pleated
surface map $ \DD \to \H3$ with pleating locus $\gamma$.  To
indicate more
clearly the relevant variables, we shall write $\PP_{\gamma}(q)$
for the
image $\DD_{\gamma}(\lambda_V,\tau_{V,W}) \subset \H3$.

If $q_0 \in \F$, then $\Im \tau_{V,W}(q_0) = 0$, hence for $q$
near $q_0$,
$\Im \tau_{V,W}(q)$ is small.  In~\cite{KSbend}, we argued that
for $\Im
\tau$ sufficiently small,
$\PP_{\gamma}(q)=\DD_{\gamma}(\lambda_V,\tau_{V,W})$ is embedded
and bounds
a convex half space in $\H3$.  It follows by proposition~7.2
of~\cite{KSbend}, that $\PP_{\gamma}(q)$ is a component of
$\bch(q)$.

There are two problems in applying this argument in the present
circumstances. First, we wish to include the case $q_0 \notin
\F$, and thus
can no longer assume that $\Im \tau_{V,W}$ is small.  Second, we
want to
prove theorem~\ref{thm:open} for an irrational lamination $\nu $
by taking
a limit of rational laminations.  Since the constant $\epsilon$
of
proposition~\ref{prop:Q2} depends on $\gamma$ and is not
uniform, (in fact
$\epsilon \sim 2\exp{(-l_{\gamma}/2)}$), the limiting process
fails,
indicating that we need to scale the approximating laminations
properly. To
resolve these problems, we digress to study the geometry of the
pleated
surfaces $\PP_{\gamma}(q)$ more carefully.

\smallskip

Fix $q_0 \in \QF, \gamma \in \S $ and a marking $(V,W)$ as
above.  Suppose
that $q \in \QF$ and that $\lambda_V(q) \in \RR^+$.  Let $
\phi_{\gamma}(q)
$ be the normalized Fuchsian group with (real) Fenchel Nielsen
coordinates
$(\lambda_V(q), \Re\tau_{V,W}(q))$.  The surface
$\PP_{\gamma}(q)$ is the
image of the pleated surface map $\DD \to \HH^3$ defined by a
pure bend
along $\delta_{\gamma}$ by $i\Im \tau_{V,W}(q)$.  We refer to
$\phi_{\gamma}(q)$ as the {\em flat structure} of
$\PP_{\gamma}(q)$.

We can associate a transverse measure $b_{\nu}(q)$ to
$\PP_{\gamma}(q)$ in
an obvious way: for any arc $\sigma$ on $\PP_{\gamma}(q)$
transverse to its
pleating locus $\gamma$, set $b_{\gamma}(q)(\sigma) =
i(\sigma,\gamma) \Im
\tau_{V,W}(q)$.  Thus we can also write $\PP_{\gamma}(q)=
\Q^p_{b_{\gamma}(q)}(i)$, where $p $ is the image of $
\phi_{\gamma}(q) $
in $\F$.

We remark that we are {\em not} making the assumptions that
$\PP_{\gamma}(q)$ is a component of $\bch(q)$, or that
$\phi_{\gamma}(q)$
is one of the flat structures $F^{\pm}(q)$ of $\bch(q)$ (see
section~\ref{sec:pllocus}); in fact, this is exactly what we
must prove.
In particular, we cannot assume that $b_{\gamma}(q)$ is the
bending measure
$pl^{\pm}(q)$.  The following result, however, gives information
about
$\phi_{\gamma}(q)$ and $b_{\gamma}(q)$ for $ q$ near $ q_0 \in
\Pr{\nu_0}$
for irrational $\nu_0$.

\begin{prop}
 \label{prop:P2} Given $\nu_0 \in ML - ML_Q$, and $q_0 \in
\Pr{\nu_0}^+\cup\F$, let $F^+(q_0)\in \F$ and $pl^+(q_0)$ be the
flat
structure and bending measure of $\bch^+(q_0)$ respectively. (If
$q_0 \in
\F$, then $pl^+(q_0)=0$ and $F^+(q_0)$ is the Fuchsian group
representing
$q_0$.)  Then, given neighborhoods $V$ of $F^+(q_0)$ in $\F$ and
$W$ of
$pl^+(q_0)$ in $ML$, there exist neighborhoods $U$ of $q_0$ in
$\QF$ and
$X$ of $[\nu_0]$ in $PML$ such that if $q \in U$,
$[\delta_{\gamma}] \in X
\cap PML_Q$ and $\lambda_{\gamma}(q) \in \RR^+$, then the flat
structure
$\phi_{\gamma}(q)$ of $\PP_{\gamma}(q)$ is in $V$ and the
transverse
measure $b_{\gamma}(q)$ is in $W$.
\end{prop}

 The idea of the proof of this proposition is that by the
convergence
lemma~\ref{lemma:convergence}, for $\nu_0 \in ML - ML_Q$, nearby
rational
laminations are close in the Hausdorff topology, so that the
bending loci
and hence the structures of the associated pleated surfaces are
also close.
The details are a technical modification of the arguments in
\cite{KSconvex} and are given in appendix~\ref{app:deferred2}.
(We remark
that the result is still true for $\nu_0 \in ML_Q$, however the
details of
the proof differ since the convergence lemma does not apply.  We
omit this
case since it is not needed here.)

The plan of the proof of theorem~\ref{thm:open} is the
following.  The hard
case to handle is $\nu \notin ML_Q$.  We shall show in
theorem~\ref{thm:P4}
below, that if $q_0 \in \Pr{\nu}$, then for $q$ in a
neighborhood of $q_0$,
if $[\delta_{\gamma}]$ is sufficiently close to $[\nu]$ in
$PML$, the
condition $\lambda_{\gamma}(q) \in \RR^+$ implies that
$\PP_{\gamma}(q)$ is
a also a component of $\bch$.  Theorem~\ref{thm:open} then
follows by an
easy limiting argument using the continuity
theorem~\ref{thm:plcont}.

We prove theorem~\ref{thm:P4} using an extension of
proposition~\ref{prop:Q2}, which we state as
proposition~\ref{prop:Q2'}.
Stated roughly it says that if $p \in \F$ and the pleated
surface
$\DD^p_{\mu}(\tau)$ associated to the quakebend
$\Q^p_{\mu}(\tau)$ is a
component of $\bch$, then the same is true of any surface
$\DD^{p'}_{\mu'}(\tau')$ obtained by quakebending a nearby
amount $\tau'$
from a nearby point $p' \in \F$ along a nearby lamination
$\mu'$.  Now, a
component of $\bch$ can be obtained from the Fuchsian group
representing
its flat structure by a pure bend along the pleating lamination
$|\mu|$.
Proposition~\ref{prop:P2} allows us to apply
proposition~\ref{prop:Q2'} to
$\PP_{\gamma}(q)$ for $[\delta_{\gamma}]$ close to $[\nu]$ and
$q$ close to
$q_0$, thus proving theorem~\ref{thm:P4}.

\begin{prop}
\label{prop:Q2'}
Let $p_0 \in \F$ be represented by $G_0 = G(p_0)$ and suppose
that $ \tau_0
\in \CC$ is such that $q_0=\Q^{p_0}_{\mu_0}(\tau_0) \in
\Pr{\mu_0}^+$.
Then there exist neighborhoods $ X,Y$ and $Z$ of $\mu_0$, $p_0 $
and
$\tau_0$ in $ ML$, $\F$ and $\CC$ respectively, such that if
$\mu \in X$,
$p \in Y$ and $\tau \in Z$, then $q=\Q^p_{\mu}(\tau) \in \QF$
and
$\DD^p_{\mu}(\tau)$ is a component of $\bch(q)$.
\end{prop}

The proof of this result is identical with the version in
\cite{KSbend}
once we note that the constants involved depend continuously on
$G$ and
$\mu$.  This follows from the following variant of lemma~8.2 of
\cite{KSbend}.

\begin{lemma}
\label{lemma:P3}
Let $X$ and $Y$ be compact sets in $ML$ and $\F$ respectively.
Then there
exist constants $d>0$ and $ K>0$ such that if $\mu \in X$ and $G
\in Y$,
and if $\sigma$ is any geodesic segment on $\DD/G$ of length
less than $d$,
then $\mu(\sigma) < K$.
\end{lemma}

We can now prove theorem~\ref{thm:P4}, which is important in its
own right.

\begin{thm}
\label{thm:P4} Suppose $\nu_0 \in ML - ML_Q$ and $q_0
\in \Pr{\nu_0}^+ \cup \F$.  Then there are neighborhoods $U$ of
$q_0$ in
$\QF$ and $X$ of $[\nu_0]$ in $ PML$ such that if
$[\delta_{\gamma}] \in
PML_Q \cap X$, $q \in U$ and if $\lambda_{\gamma}(q) \in \RR^+$,
then
$\PP_{\gamma}(q)$ is a component of $\bch(q)$.
\end{thm}

\begin{proof} By proposition~\ref{prop:P2}, there are
neighborhoods $X$ of $[\nu_0]$ in $ PML$ and $U$ of $q_0$ in
$\QF$ such
that for $ q \in U$ and $[\delta_{\gamma} ] \in X$, the flat
structures
$F^+(q_0)$ of $\bch^+(q_0)$ and $\phi(q )$ of $\PP_{\gamma}(q)$
are close
in $\F$, and the transverse measures $pl^{+}(q_0)$ and
$b_{\gamma}(q)$ are
close in $ML$.

As remarked earlier, $\bch^+(q_0)$ is just the pleated surface
obtained
from $F^+(q_0)$ under a pure bend by $i$ along the measured
lamination
$pl^+(q_0)$ while $\PP_{\gamma}(q)$ is obtained from
$\phi_{\gamma}(q)$ by
a pure bend by $i$ along $b_{\gamma}(q)$.  The result now
follows from
proposition~\ref{prop:Q2'}.
\end{proof}

We now prove theorem~\ref{thm:open}.

\begin{proof}
Suppose first that $\nu \in ML_Q$.  In this case the result is
just
 proposition~\ref{prop:Q2'}, using proposition~\ref{prop:7.6} as
a
 substitute for the condition $\Im \tau$ near $0$ when the base
point $q_0$
 is not Fuchsian.

  Suppose therefore that $\nu \notin ML_Q$, and pick $\nu_n \in
ML_Q$, $\nu_n
\to \nu$.  Find neighborhoods $U$ of $q_0$ in $\QF$ and $X$ of
$[\nu]$ in $
PML$ satisfying the conclusion of theorem~\ref{thm:P4}.

Assume $q \in U$ and $\lambda_{\nu}(q) \in \RR^+$. Since
$\lambda_{\nu_n}
\to \lambda_{\nu}$ uniformly on $U$, and since $\lambda_{\nu}$
is
non-constant on $U$, by Hurwitz's theorem we can find $q_n \in
U$, $ q_n
\to q$, such that $\lambda_{\nu_n}(q_n) =\lambda_{\nu}(q)$, and
in
particular such that $\lambda_{\nu_n}(q_n) \in \RR^+$.  Applying
theorem~\ref{thm:P4}, we see that for sufficiently large $n$,
$\PP_{|\nu_n|}(q_n)$ is one of the components of $\bch^+(q_n)$
so that $q_n
\in \Pr{\nu_n}^+$.  Hence, by the continuity
theorem~\ref{thm:plcont}, we
get $q \in \Pr{\nu}^+$.  This completes the proof of
theorem~\ref{thm:open}.
\end{proof}

\begin{cor}
\label{cor:openinqkbendplane}
 Suppose $\mu,\nu \in ML$, $[\mu] \neq [\nu]$.  Let $q_0 \in
\Pr{\mu,\nu}^+
\cup \F$ and let $\Q_{\mu}^{q_0}$ be the $\mu$-quakebend plane
based at
$q_0$. There exists a neighborhood $U$ of $ q_0$ in
$\Q_{\mu}^{q_0}$ such
that if $q \in U$ and $\lambda_{\nu}(q) \in \RR^+$, then $q \in
\Pr{\mu,\nu} \cup \Pr{\nu,\mu} \cup \F$.
\end{cor}

\begin{proof} This is just theorem~\ref{thm:open} applied in the
quakebend plane $\Q_{\mu}^{q_0}$.  We can prove it either by
applying
proposition~\ref{prop:Q2} to see that for $q \in \Q_{\mu}^{q_0}$
near
$q_0$, we have $q \in \Pr{\mu} \cup \F$, and then applying
theorem~\ref{thm:open} to $\nu$; or by noting that since
$\lambda_{\mu}$ is
constant on $\Q_{\mu}^{q_0}$ and real valued at $q_0$, we can
apply
theorem~\ref{thm:open} first to $\mu$ and then to $\nu$.
\end{proof}

\begin{remark} {\rm  The condition $\lambda_{\gamma}(q) \in
\RR^+$ is key in proposition~\ref{prop:P2} and in
theorem~\ref{thm:P4}.  We
can always find a pleated surface $\Pi$ whose pleating locus
$\sigma$
contains the geodesic $\gamma$.  In general, however, $\sigma$
properly
contains $\gamma$ and has leaves spiralling into $\gamma$, and
thus carries
no transverse measure.  Then, even though $[\delta_{\gamma}]$ is
near
$[pl^{\pm}]$ in $PML$, the pleated surface $\PP_{\gamma}$
realizing
$\gamma$ (see \cite{ CEpG, ThuN}) is not necessarily embedded;
moreover,
even if it is, neither of the half spaces it bounds in $\H3$
will be
convex.  The point is that the condition $\sigma = \gamma$ is
equivalent to
$\lambda_{\gamma}(q) \in \RR^+$.  }
\end{remark}

\subsection{Consequences of theorem~\ref{thm:open}}

\label{sec:consequences}

From theorem~\ref{thm:open} we obtain the following local
extension of the  picture of Fuchsian space described in
section~\ref{sec:fuchsian}.

\begin{thm}
\label{thm:bendingaway}
  Let $\mu,\nu \in ML$, $i(\mu,\nu)>0$, $p \in \F$.  Then there
is a
   neighborhood $U$ of $ p$ in $\QF$ such that \begin{enumerate}
\item if
   $p \notin \F_{\mu,\nu}$ then $\Pr{\mu,\nu} \cap U =
\emptyset$, while
   \item if $p \in \F_{\mu,\nu}$ then the $\RR$-locus of
$\lambda_{\nu}$ in
   $U$ is exactly
$$(\Pr{\mu,\nu} \cup \Pr{\nu,\mu} \cup \F) \cap U.$$
\end{enumerate}

In the second case, let  $p=p(\mu,\nu,\c) \in \F_{\mu,\nu}$,
let $\Q^{p}_{\mu}$ be the quakebend plane
along $\mu$ based at  $p$ and let $V =U \cap \Q^{p}_{\mu}$ .
 Then  ${\lambda_{\nu}}|_{V}$
 has a simple critical point at $p$
and  $\lambda_{\nu}^{-1}(\RR^+) \cap (V - \F)$ has exactly two
components, one lying in $\Pr{\mu,\nu}$ and the other in
$\Pr{\nu,\mu}$.
\end{thm}

\begin{proof}   Part 1 follows since
for $p \notin \F_{\mu,\nu}$,
there exists a neighborhood $U$ of $ p$ in the quakebend plane
$\Q_{\mu}^{p}$ based at $p$ such that
$\lambda_{\nu}^{-1}(\RR^+) \cap U \subset \F$.

By \cite{KerckN, Wolp}, $\lambda_{\nu}|_{\E_{\mu,\c}}$ has
exactly
one critical point  at $p$ and it is simple.  Thus part 2 is a
restatement of corollary~\ref{cor:openinqkbendplane} with
$q_0=p$.
\end{proof}

We note that this theorem provides an alternative proof of
theorem~\ref{thm:bend}.

\smallskip
We can also now prove the density
theorem~\ref{thm:ratlplanesdenseI} of the introduction.
First, we  need a bound on the bending angle in a quakebend
plane.

\begin{prop}
\label{prop:genbendinglimit} Suppose $\mu \in ML$,  $q \in
\Pr{\mu} \cup \F$ and let $\Q^q_{\mu}$ be the quakebend plane
along $\mu$
based at $q$ with parameter $\tau =\tau_{\mu}$.  Given $K > 0$,
there
exists $B > 0$ such that if $|\Re \tau| <K$ and $|\Im \tau| >
B$, then
$\Q^q_{\mu}(\tau) \notin \Pr{\mu}$.
\end{prop}

The statement $\Q^q_{\mu}(\tau) \notin \Pr{\mu}$ means that
either
$\Q^q_{\mu}(\tau) \notin \QF$ or that
$\Q^q_{\mu}(\tau) \in \QF$ but $ [pl^+(\Q^q_{\mu}(\tau) )] \neq
[\mu]$.
  We show that, under the hypotheses of
the proposition, $\Q^q_{\mu}(\tau)$ fails to be in $\Pr{\mu}$
because the
surface obtained by bending along $\mu$ is not embedded.  This
may or may
not imply that $\Q^q_{\mu}(\tau) \in \QF$.  The proof is given
in
appendix~\ref{app:deferred3}, see also~\cite{McMeq} theorem 6.2.

  As an immediate corollary we have
\begin{prop}
\label{prop:genlnonconst}  Suppose $q \in \QF$,  $q \in
\Pr{\mu,\nu} \cup
\F$.  Then the holomorphic function $\lambda_{\nu}(q)$ is
non-constant on
 $\Q^q_{\mu} \cap \QF$.
\end{prop}

\begin{proof}
Since $q \in \Pr{\nu}$ we know $\lambda_{\nu}(q)\in \RR^+$.  By
 construction $\lambda_{\mu}(q)=c>0$ for all $q \in \Q^q_{\mu}$.
 Suppose
 that $\lambda_{\nu}(q)=d>0$ for all $q \in \Q^q_{\mu} \cap
\QF$.  By
 theorem~\ref{thm:open}, $\Pr{\mu,\nu}$ is open in $\Q^q_{\mu}$.

Now suppose that $q_n= \Q^q_{\mu}(\tau_n) \in \Pr{\mu,\nu}$ and
that $
\tau_n \to \tau_{\infty}$.  Since $l_{\mu}(q_n)=c$ and $
l_{\nu}(q_n)=d$
for all $n$, it follows from theorem~\ref{thm:thmq} that $q_n
\to
q_{\infty} \in \QF$.  By theorem~\ref{thm:plcont}, $q_{\infty}
\in
\Pr{\mu,\nu} \cup \F$.  Clearly, $q_{\infty} =
\Q^q_{\mu}(\tau_{\infty}) $
and so $\Pr{\mu,\nu}$ is closed in $\Q^q_{\mu} - \F$.  Therefore
$\Pr{\mu,\nu}$ is a connected component of $\Q^q_{\mu} - \F$ and
must be
one of the half planes $\Re \tau_{\mu} >0$ or $\Re \tau_{\mu}
<0$,
contradicting proposition~\ref{prop:genbendinglimit}.
\end{proof}

Finally we can prove theorem~\ref{thm:ratlplanesdenseI}.

\begin{genericem}{Theorem~\ref{thm:ratlplanesdenseI}}
The rational pleating varieties $\Pr{\mu,\nu}, \mu,\nu \in ML_Q$
are dense in $\QF$.
\end{genericem}

\begin{proof}  Let $q \in \QF$ and let $\mu \in [pl^+(q)], \nu
\in [pl^-(q)] $.  By theorem~\ref{thm:realI},
$\lambda_{\mu}(q),\lambda_{\nu}(q) \in \RR^+$.  Clearly, we may
as well
assume $ \mu \notin ML_Q$.  Find a sequence $\{\mu_n\} \in ML_Q,
\mu_n \to
\mu$.  By Hurwitz's theorem in $\QF$, we can find points $q_n
\to q$ with
$\lambda_{\mu_n}(q_n) \in \RR^+$ and so by theorem~\ref{thm:P4},
$q_n \in
\Pr{\mu_n}$ for large enough $n$.  If $\nu \in ML_Q$ we are
done, otherwise
find $\{\nu_n\} \in ML_Q, \nu_n \to \nu$.  By
proposition~\ref{prop:genlnonconst}, $\lambda_{\nu_n}$ is
non-constant on
$\Q^{q_n}_{\mu_n} \cap \QF$ and we can apply Hurwitz's theorem
again in
$\Q^{q_n}_{\mu_n} \cap \QF$ to find $q_n'$ near $q_n$, such that
$q_n' \to
q$ and such that $\lambda_{\nu_n}(q_n') \in \RR^+$.  By
theorem~\ref{thm:P4} again, $q_n' \in \Pr{\mu_n,\nu_n}$ for
large enough
$n$.
\end{proof}

\section{Pleating rays and planes}
\label{sec:Pleating Rays}

In this section, we apply the local and limit pleating theorems
to prove our main results
theorems~\ref{thm:plraysI} and~\ref{thm:plplanesI}  of the
introduction.

Recall from section~\ref{sec:plvars} the definition of the
pleating ray $$
\Pr{\mu,\nu,\c} = \{ q \in \Pr{\mu,\nu} :l_{\mu}(q) =\c \},$$
where
$(\mu,\nu) \in ML \times ML$, and $\c>0$.  Pleating rays are the
basic
building blocks out of which we construct pleating planes and
the
$BM$-slices mentioned in the introduction.  Notice that, because
of
theorem~\ref{thm:realI}, we can equally well define
$$
\Pr{\mu,\nu,\c} = \{ q \in \Pr{\mu,\nu} :\lambda_{\mu}(q) =\c
\}.$$

Our results will justify the names ``rays'' and ``planes''.

The main work is in the study of the pleating rays.  Our
strategy is as
follows.  We begin by applying the limit pleating theorem and
the local
pleating theorem to obtain some general results about
$\Pr{\mu,\nu}$ for
arbitrary $\mu, \nu \in ML$.  We then prove
theorem~\ref{thm:plraysI} in
the case where $[\mu] = [\delta_{\gamma}], [\nu] = [
\delta_{\gamma'}]$ and
$(\gamma,\gamma')$ is a marking for $\torus$.  We show that in
this case
$\Pr{\delta_{\gamma},\delta_{\gamma'}, \c}$, which we call an
{\em integral
pleating ray}, is a straight line segment in the quakebend plane
$\Q_{\gamma,\c}$.  Using the integral rays we derive constraints
on the
rays $\Pr{\delta_{\gamma},\nu, \c} \subset \Q_{|\mu|, \c}$ for
arbitrary
$\nu$; using our general results we are then able to deduce
theorem~\ref{thm:plraysI} in the general case.  Finally, we
apply
theorem~\ref{thm:plraysI} to deduce theorem~\ref{thm:plplanesI}.

\subsection{Pleating rays}

\smallskip
In the four lemmas which follow, $\mu, \nu$ are arbitrary
laminations in
$ML$ and, as usual, $\Q^q_{\mu}$ denotes the $\mu$-quakebend
plane through
$q \in \Pr{\mu} \cup \F$.

\begin{lemma}
\label{lemma:conncpts} Let $q \in \Pr{\mu,\nu}$.
The set $\Pr{\mu,\nu} \cap \Q^q_{\mu}$ is a union of connected
components
of the $\RR$-locus of $\lambda_{\nu}$ in $(\QF-\F) \cap
\Q^q_{\mu}$.
\end{lemma}

\begin{proof}
We have to show that $\Pr{\mu,\nu} \cap \Q^q_{\mu}$ is open and
closed in
the $\RR$-locus of $\lambda_{\nu}$ in $(\QF - \F) \cap
\Q^q_{\mu}$.  The
openness is the local pleating theorem~\ref{thm:open} and
closure follows
by the continuity theorem~\ref{thm:plcont}.  \end{proof}

If $ \nu = \k \delta_{\gamma}, \gamma \in {\cal S}$, we obtain a
stronger
result.  Let $V \in G$ represent $\gamma$.  In this case, by
proposition~\ref{prop:Q1}, trace $\tr V$ is defined and
holomorphic on all
of $\Q^q_{\mu}$ (including the part outside $\QF$), and we
obtain a version
of lemma~\ref{lemma:conncpts} for the $\RR$-locus of
$\lambda_{\gamma}$ in
$\Q^q_{\mu}$.  Define the hyperbolic locus of $\gamma$ in
$\Q^q_{\mu}$ as
$\{ q \in \Q^q_{\mu}: \tr V \in \RR, |\tr V| >2 \}$.

\begin{lemma}
\label{lemma:ratconncpts} Let $ \nu=\k \delta_{\gamma} \in ML_Q$
and let $q \in \Pr{\mu,\nu}$. Let $V \in G$ represent $\gamma$.
Then the set $\Pr{\mu,\nu} \cap \Q^q_{\mu}$ is a
union of connected components of the hyperbolic locus of $\tr V$
in
$ \Q^q_{\mu} - \F$.
\end{lemma}

\begin{proof}
 The   openness
 follows as above, using
the local
pleating theorem~\ref{thm:open}.  The closure follows from
theorem~\ref{thm:thmq}.  The point is first, that  length and
trace
are related by the trace formula
  $\tr V = 2 \cosh( l_{\gamma} /2) $, and second, that if we
reach a limit
  point at which $|\tr V| >2$, then $l_{\gamma}  >0$ so that
  by the second part of theorem~\ref{thm:thmq} we must still be
in $\QF$.
  (See~\cite{KStop} proposition 5.4
 for a more elementary proof without using
theorem~\ref{thm:thmq}.)
\end{proof}

 This is a strong result.  The point is, that starting from a
point we know
is in $\QF$, the lemma asserts that if we move along branches of
the
hyperbolic locus, then we stay in $\QF$ until we reach a
boundary point of
$\partial \QF$ at which $|\tr V|=2$.  This observation is what
makes it
possible to use the pleating invariants for computations of
$\partial \QF$,
see theorem~\ref{thm:bdryI} of the introduction.

With the notation of lemma~\ref{lemma:conncpts}, set
$c=\lambda_{\mu}(q)$.
Clearly, $\Pr{\mu,\nu} \cap \Q^q_{\mu} = \Pr{\mu,\nu,c}$.  As
usual, we let
$p_{\mu,\nu,\c} \in \F_{\mu,\nu}$ be the minimal point for the
length
function $l_{\nu}$ on the earthquake path $\E_{\mu,\c}$.  The
following two
lemmas make essential use of theorem~\ref{thm:thmq}.

\begin{lemma}
\label{lemma:image}
 Let $ q \in \Pr{\mu, \nu}$ and let $c=\lambda_{\mu}(q)$.  The
image of
$\Pr{\mu,\nu} \cap \Q^q_{\mu}$ under the map $\lambda_{\nu}$ is
a union of
intervals of the form $(0, \infty)$, $(0,\d)$ and $(\d, \infty)$
where $\d
= f_{\mu,\nu}(\c) =l_{\nu}(p_{\mu,\nu,\c})$.  Moreover, there is
at most
one component of $\Pr{\mu,\nu} \cap \Q^q_{\mu}$ whose image is
$(0,\d)$;
the closure of such a component meets $\F$ exactly in $
p(\mu,\nu,\c)$.
\end{lemma}

\begin{proof}
Let $K$ be a connected component of $\Pr{\mu,\nu,\c}$. By
theorem~\ref{thm:realI}, $\lambda_{\nu}|_K$ is real valued and,
by
proposition~\ref{prop:genlnonconst}, it is non-constant on
$\Q^q_{\mu}$.
Since it is holomorphic, it is not locally constant and thus not
constant
on $K$.  Therefore by lemma~\ref{lemma:conncpts} the image $I_K$
of
$\lambda_{\nu}|_K$ is an open interval in $\RR^+$.

Suppose that $r \in \RR^+$ and that there is a sequence $\{q_n\}
\in K$
such that $\lambda_{\nu}(q_n) \to r$.  Since $\lambda_{\mu}(q_n)
= \c$, by
theorem~\ref{thm:thmq} a subsequence of $\{G(q_n)\}$ has an
algebraic limit
$G_{\infty}$.  Furthermore, since $\lambda_{\nu}(q_n) \to r>0$,
the group
$G_{\infty}$ is represented by a point $q \in \QF$ such that
$\lambda_{\nu}(q)=r$.  If $q \in \QF -\F$ then by
theorem~\ref{thm:plcont},
$q \in K$ so that $r \in I_K$.  On the other hand, if $q \in \F$
then by
theorem~\ref{thm:bendingaway}, $q=p(\mu,\nu,\c)$ and
$r=\lambda_{\nu}(q)=f_{\mu,\nu}(\c)=\d$.  Thus
$\lambda_{\nu}(K)$ is open
and closed in $(0,\d) \cup (\d,\infty)$. The result follows from
theorem~\ref{thm:bendingaway}.
\end{proof}

\begin{lemma}
\label{lemma:retoinfty}
Let $ q \in \Pr{\mu, \nu}$ and let $c=\lambda_{\mu}(q)$.  Let
$\tau$ denote
the quakebend parameter in the quakebend plane $\Q_{\mu}^q$.
Suppose that
the points $q_n \in \Pr{\mu,\nu,\c}$ are represented by the
quakebend
parameter $\tau_n$ and that $\lambda_\nu(q_n) \to \infty$.  Then
$|\Re(\tau_n)| \to \infty$.
\end{lemma}

\begin{proof}
Since $q_n \in \Pr{\mu,\nu,\c}$ we know $\lambda_{\nu}(q)$ is
real. Moreover, $\lambda_{\nu}(q) \le l_{\nu}(F^+(q_n))$; that
is,
$\lambda_{\nu}(q)$ is bounded above by the length of $\nu$ on
the flat
structure of $\bch^+/G(q_n)$.  This flat structure is determined
by the
length of $\mu$, which is fixed, and the earthquake parameter
$\Re(\tau_n)$.  Thus if $|\Re(\tau_n)|$ is bounded, so is
$\lambda_{\nu}(q_n)$.
\end{proof}

We can now start investigating the integral pleating rays.
Suppose that
$[\mu] = [\delta_{\gamma}], [\nu] = [ \delta_{\gamma'}]$ and
$(\gamma,\gamma')$ is a marking for $\torus$.  For simplicity,
we write
$\Pr{\gamma}$ for $\Pr{\delta_{\gamma}}$ and so on.  Let
$(\lambda_V,\tau_{V,W}) \in \CC^2$ be complex Fenchel Nielsen
coordinates
relative to a marked pair of generators $(V,W)$ corresponding to
$(\gamma,\gamma')$.  As in section~\ref{sec:qbplanes}, we denote
by
$\Q_{\gamma,\c}$ the slice $ \{(\c,\tau)\} \subset \CC^2$;
$\Q_{\gamma,\c}$
is the quakebend plane along $\gamma$ that meets $\F$ along the
earthquake
path $\E_{\gamma,\c}$.  We denote points in this slice simply by
the
parameter $\tau =\tau_{\delta_{\gamma}}$.  As usual, $ \tau = 0$
corresponds to the point $p(\gamma, \gamma',\c) \in \F$, while
$\Im \tau =
0$ is the earthquake path $\E_{\gamma,\c}$.

For $m \in \ZZ$, the pair $(\gamma,\gamma^m\gamma')$ is a pair
of marked
generators for $\pi_1(\torus)$ corresponding to the pair of
generators
$V,V^mW$ for $G$.  Clearly $\Pr{\gamma,\gamma^m\gamma',\c}
\subset
\Q_{\gamma, \c}$.  The generators $V,VW$ are obtained from the
pair $V,W$
by the map induced by a Dehn twist about $\gamma$.  The
basepoint relative
to which we measure the twist parameter changes and we find
$\tau_{V,VW} =
\tau_{V,W} + \lambda_V$; similarly, $\tau_{V,V^mW} = \tau_{V,W}
+
m\lambda_V$.

 The following formula is derived in \cite{PS} for any pair
(V,W) of marked
 generators for $G$:

 \begin{equation}
\label{eqn:bending}
 \cosh\frac{ \tau_{V,W}}{2} = \pm \cosh \frac{\lambda_{W}}{2}
\tanh
\frac{\lambda_V}{2}
\end{equation}
By our conventions, $\Re \lambda_{V},\, \Re \lambda_{W}
 >0$, so that we should choose the $+$ sign on $\F$ and hence
everywhere in
 $\QF$.

Applying this formula to the generators $(V,V^{-m}W)$ we find
\begin{equation}
\label{eqn:bending2}
 \cosh\frac{\tau -
  m\lambda_V}{2}=
\cosh\frac{\lambda_{V^{-m}W}}{2}\tanh\frac{\lambda_V}{2}.
\end{equation}
In particular, at $\tau = m\c$ we have
\begin{equation}
\label{eqn:critical}
\textstyle{ 1=  \cosh{\frac{\lambda_{\gamma^{-m}\gamma'}}{2}}
\tanh{\frac{\lambda_{\gamma}}{2}}}.
\end{equation}
or equivalently
\begin{equation}
\label{eqn:critical1}
\textstyle{\sinh{\frac{\c}{2}} \sinh
{\frac{\lambda_{\gamma^{-m}\gamma'}}{2}} = 1.}
\end{equation}

\begin{prop}
\label{prop:intrays}  Let $(\gamma,\gamma')$ be a marked pair of
generators for $\pi_1(\torus)$ and let $\c>0$.  Then for $m \in
\ZZ$,
$\Pr{\gamma,\gamma^{-m}\gamma',\c}$ and $
\Pr{\gamma^{-m}\gamma',\gamma,
\c}$ are the two line segments $\Re \tau = m\c, \, |\Im \tau| <
2\arccos
{\tanh{\frac{\c}{2}}}$ in $\Q_{\gamma, \c}$.  The two line
segments $
\Re\tau = m\c, |\Im\tau| \ge 2\arccos {\tanh{\frac{\c}{2}}}$ in
$\Q_{\gamma, \c}$ have empty intersection with $\QF$.
\end{prop}

\begin{remark}
{\rm Which of the two segments corresponds to
$\Pr{\gamma,\gamma^{-m}\gamma',\c}$ and which to $
\Pr{\gamma^{-m}\gamma',\gamma, \c}$ depends on our convention
for measuring
$\tau$ and is not important here.}
\end{remark}

\begin{proof}
Because $\tau_{V,V^{-m}W} = \tau_{V,W} -m \c$, we may restrict
ourselves to
the case $m = 0$.  From lemma~\ref{lemma:ratconncpts},
$\Pr{\gamma,\gamma'}$ is a union of connected components of the
hyperbolic
locus of $ \gamma' $ in $\Q_{\gamma, \c} - \F$, and by
theorem~\ref{thm:bendingaway} there is a unique component $K$
whose closure
meets the critical line $\F_{\gamma, \gamma'}$ in $p(\gamma,
\gamma', \c)$.

From equation~(\ref{eqn:bending}) ,
$$\cosh{\frac{\tau}{2}}= \cosh{\frac{\lambda_{\gamma'}}{2}}
\tanh{\frac{\lambda_{\gamma}}{2}}.$$ Thus the $\RR$-locus of $
\lambda_{\gamma'} $ in $\Q_{\gamma, \c}$ is the set defined by
$\cosh{\frac{\tau}{2}} \in \RR$, or equivalently, \{$\Re \tau =
0 \}\cup
\{\Im \tau = 0$\}.  The real axis $\Im \tau = 0$ corresponds to
$\E_{\gamma,\c} = \Q_{\gamma,\c} \cap \F $ and we see easily
(see
lemma~\ref{lemma:ratconncpts} ) that the connected components of
the
hyperbolic locus of $\gamma' $ in $\Q_{\gamma, \c} - \F$ which
meet the
real axis are the two segments $ 0 <|\Im \tau| < 2\arccos
{\tanh{\frac{\c}{2}}}$.  One of these segments must be the
component $K$
and the other is the corresponding component for
$\Pr{\gamma',\gamma}$.
Each of these segments is mapped bijectively by
$\lambda_{\gamma'}$ to $[0,
2\arccos {\tanh{\frac{\c}{2}}})$.

 Now on the imaginary axis, we have $
\cosh{\frac{\lambda_{\gamma'}}{2}}
 \le (\tanh{\frac{\c}{2}})^{-1}$, and hence by
lemma~\ref{lemma:image},
 $\Pr{\gamma,\gamma', \c}$ and $\Pr{\gamma',\gamma,\c}$ have no
other
 components.

Finally we have to show that that no other points on the
imaginary axis lie
 in $\QF$.  Equation~(\ref{eqn:bending2}) holds for groups in
 $\Q_{\gamma,\c}$ even when they are outside $\QF$.  On this
axis,
 therefore, we always have
   $$ -1 \le \cosh{\frac{\lambda_{\gamma'}}{2}}
\tanh{\frac{\c}{2}} \le 1
.$$

In \cite{PS} proposition~6.2, it is shown by a direct argument
that if
${\lambda_{\gamma'}} \in \RR$ and the above inequality is
strict, then the
group generated by $V,W$ is quasifuchsian and contained in
$\Pr{\gamma,\gamma'}$.  Moreover, in this situation, this group
is
determined by ${\lambda_{\gamma}} $ and ${\lambda_{\gamma'}} $
up to
conjugacy.  If equality holds, the group represents the unique
point
$p(\gamma, \gamma', \c) \in \F$.  These are the therefore the
groups we
have already discussed.

Since $\cosh{\frac{\lambda_{\gamma'}}{2} }\in \RR$, the only
other
possibility is that $\lambda_{\gamma'}$ is purely imaginary.  In
this case
the corresponding group element would have to be elliptic which
is
impossible in $\QF$.
\end{proof}

We can now obtain a bound on the pleating rays
$\Pr{\gamma,\nu,\c} $ for
 arbitrary $ \nu \in ML$.

\begin{cor}
\label{cor:retaubdd} Let $\nu \in ML$, $i(\nu,\gamma)>0$.
Then $|\Re \tau|$ is bounded on each component of
$ \Pr{\gamma,\nu,\c}$, where $\tau$
denotes the quakebend parameter $\tau_{\delta_{\gamma}}$ in
$\Q_{\gamma,\c}$.
\end{cor}

\begin{proof}  If along some component of $ \Pr{\gamma,\nu,\c} $
in $\Q_{\gamma,\c}$, $|\Re{\tau}| \to \infty$, the component
would have to
intersect infinitely many of the lines $\tau = m\c + i \theta,
\theta \in
\RR$.  According to proposition~\ref{prop:intrays}, however,
each such line
is the union of the integral pleating rays $
\Pr{\gamma,\gamma^{-m}\gamma,\c}, \,
\Pr{\gamma,\gamma^{-m}\gamma,\c}$, the
point $p(\gamma,\gamma^{-m}\gamma,\c) \in \F$, and points not in
$\QF$.
This is impossible.
\end{proof}

\smallskip
We can now prove theorem~\ref{thm:plraysI} on the structure of
the pleating
rays. Recall from section~\ref{sec:qbplanes} that $\Q_{\mu,\c}$
is the
quakebend plane along $\mu$ which meets $\F$ along the
earthquake path
$\E_{\mu,\c}$.

\begin{genericem}{Theorem~\ref{thm:plraysI}}
Let $\mu,\nu$ be measured laminations on $\torus$ with
$i(\mu,\nu)>0$ and
let $\c>0$.  Then the set $\Pr{\mu,\nu,\c} \subset \QF$ on which
$[pl^+]=[\mu]$, $[pl^-]=[\nu]$ and $l_{\mu} = \c$, is a
non-empty connected
non-singular component of the $\RR$-locus of the restriction of
$\lambda_{\nu}$ to $\Q_{\mu,\c}$.  The restriction of
$\lambda_{\nu}$ to
$\Pr{\mu,\nu,\c} $ is a diffeomorphism onto its image
$(0,f_{\mu,\nu}(\c))
\subset \RR^{+}$.
\end{genericem}

\begin{proof}
We assume first that $\mu \in ML_Q$; without loss of generality
we may take
$\mu = \delta_{\gamma},\, \gamma \in \S$.  Let $c>0$ and let $K$
be a
component of $ \Pr{\gamma,\nu,\c} $.  By
corollary~\ref{cor:retaubdd},
$|\Re \tau|$ is bounded on $K$.  By lemma~\ref{lemma:retoinfty},
$\lambda_{\nu}|_K$ is bounded and hence by
lemma~\ref{lemma:image} the
image is the interval $(0,\d)$ where $\d = f_{\gamma,\nu}(\c)$.
Moreover,
there exist points $ \tau_n \in K$ , $\tau_n \to
p(\gamma,\nu,\c) \in
\F_{\gamma,\nu}$.

Now by theorem~\ref{thm:bendingaway}, there is only one branch
of
$\lambda_{\nu}^{-1}(\RR^+)$ near $p(\gamma,\nu,\c)$; thus if the
degree of
$\lambda_{\nu}|_K$ were greater than one, there would be points
$\tau_n'
\in K$ with $\lambda_{\nu}(\tau_n') \to \d$, but with $\tau_n'
\to
q_{\infty} \in \QF-\F$.  Then, by lemma~\ref{lemma:image},
$\lambda_{\nu}
(K) \supset (0,\infty)$, which is impossible.

Now we remove the restriction that $\mu \in ML_Q$.  Suppose that
$q \in
\Pr{\mu,\nu,\c}$.  We have to replace the plane $\Q_{\gamma,\c}$
by the
plane $\Q_{\mu}^q$, in which we denote the quakebend parameter
$\tau_{\mu}$
by $\tau$.  Because there are no integral pleating rays if $\mu$
is
irrational, we need another argument to bound $\Re \tau$.

Choose a sequence $\nu_n \in ML_Q$ such that $\nu_n \to \nu$.
By
theorem~\ref{thm:welldefined} the holomorphic function
$\lambda_{\nu}(q)$
is continuous in $\nu$ and by
proposition~\ref{prop:genlnonconst} it is
nonconstant.  Thus we can apply Hurwitz's theorem in
$\Q^q_{\mu}$ to find
$q_n \in \Q^q_{\mu}$ such that $q_n \to q$ and
$\lambda_{\nu_n}(q_n) \in
\RR^+$.  By theorem~\ref{thm:P4}, for large enough $n$, $q_n \in
\Pr{\mu,\nu_n,\c}$.  Now because $\nu_n \in ML_Q$, we can apply
the
argument above with the roles of $\mu$ and $\nu_n$ reversed to
deduce that
$\lambda_{\mu}(q_{n}) < f_{\nu_n, \mu}(\lambda_{\nu_n}(q_n)$.
Thus, since
$f_{\nu_n,\mu}$ is monotonic decreasing, we have that
$\lambda_{\nu_n}(q_n)
< f_{\nu_n,\mu}^{-1}(\lambda_{\mu}(q_{n}))$.  Since
$f_{\nu_n,\mu}^{-1}=f_{\mu,\nu_n}$ we conclude that
$\lambda_{\nu_n}(q_n) <
f_{\mu,\nu_n}(c)$.

 Because $\nu_n \to \nu$, by corollary~\ref{cor:X} and
theorem~\ref{thm:welldefined} we have $f_{\mu, \nu_n}(\c) \to
f_{\mu,
\nu}(\c)$ so that $\{\lambda_{\nu_n}(q_n)\}$ is bounded by a
constant
depending only on $\mu,\nu$ and $\c$.  The remainder of the
argument is as
before.
\end{proof}

As an immediate corollary we have

\begin{cor}
\label{cor:uniqplane} If $q \in \Pr{\mu}$, then $G(q) $ is
obtained from a group $G(p), p \in \F$ by a quakebend
$\Q_{\mu}^p(\tau^*)$
along $\mu$.  Moreover, there is a quakebend path $ \sigma
:[0,1] \to \CC$
in $\QF$ from $p$ to $q$, or, in the coordinate of
$\Q_{\mu}^p(\tau)$,
$\sigma(0) = 0$, $\sigma(1) = \tau^*$ and $\Q_{\mu}^p(\sigma(t))
\in \QF, 0
\le t \le 1$.
\end{cor}

This settles the question about uniqueness of quakebend planes
raised at
the end of section~\ref{sec:genqbs}.

\begin{remark}
{\rm In~\cite{KStop}, we studied the {\em Maskit slice} for
punctured tori
in terms of pleating rays with a similar definition to the
above.  In
particular, theorem 7.2 of~\cite{KStop}, asserts a
non-singularity result
similar to that in theorem~\ref{thm:plraysI}.  It has been
pointed out to
us by Y. Komori that our proof in~\cite{KStop} in the case of
rays $\nu
\notin ML_Q$ is incorrect.  In fact, we need an openness result
like
theorem~\ref{thm:open} above.  The methods above also prove the
important
result, omitted in~\cite{KStop}, that the range of the length
function on
an irrational ray in the Maskit slice is $(0,\infty)$.  We refer
to~\cite{KomS} for a corrected version of the argument
in~\cite{KStop}.  }
\end{remark}

\subsection{Pleating planes}
\label{sec:plplanes}

We are finally able to prove theorem~\ref{thm:plplanesI} on the
structure
of the pleating varieties $\Pr{\mu,\nu}$.  As in the
introduction, let
$L_{\mu,\nu}:\QF \rightarrow \CC^2$ be the map $q \mapsto
(\lambda_{\mu}(q),\lambda_{\nu}(q))$.

\begin{genericem}{Theorem~\ref{thm:plplanesI}}
Let $(\mu,\nu)$ be measured laminations on $\torus$ with
$i(\mu,\nu)>0$.
Then the set $\Pr{\mu,\nu} \subset \QF$ on which $[pl^+]=[\mu],
[pl^-]=[\nu]$ is a non-empty connected non-singular component of
the
$\RR^2$-locus in $\QF-\F$ of the function $\L_{\mu,\nu}$.  The
restriction
of $\L_{\mu,\nu}$ to $\Pr{\mu,\nu}$ is a diffeomorphism to the
open region
under the graph of the function $f_{\mu,\nu}$ in $\RR^{+} \times
\RR^{+}$.
\end{genericem}

\begin{proof}
 By theorem~\ref{thm:realI}, the map
$L_{{\mu,\nu}|_{\Pr{\mu,\nu}}}$ takes
values in $\RR^+ \times \RR^+$.  That $L_{\mu,\nu}$ restricted
to
$\Pr{\mu,\nu}$ is injective follows immediately from the
injectivity of
$\lambda_{\nu}$ on each pleating ray $\Pr{\mu,\nu,\c}$. Hence, $
\Pr{\mu,\nu}$ is a non-singular $\RR^2$-locus in $\QF-\F$.  The
statement
about the image of $L_{\mu,\nu}$ follows from
theorem~\ref{thm:plraysI}.
\end{proof}

We remark that a similar proof shows that $\Pr{\mu,\nu}$ and
$\Pr{\nu,\mu}
$ are the unique connected components of the $\RR$-locus of
$L_{\mu,\nu}$
in $\QF- \F$ whose closure in $\QF$ meets $\F$ in
$\F_{\mu,\nu}$.

We also remark that if in theorem~\ref{thm:plplanesI} we replace
$\mu,\nu$
by $\mu'=s\mu,\nu'=t\nu, \, s,t,\in \RR^+$, then $\Pr{\mu,\nu}$
is
unchanged and the length function $L_{\mu',\nu'}$ is simply a
rescaling of
$L_{\mu,\nu}$:
 $$L_{\mu',\nu'}(q) = (s\lambda_{\mu}(q), t\lambda_{\nu}(q)).$$

 Our main result, theorem~\ref{thm:plinvarsI}, that a group in
$\QF$ is
characterized by its pleating invariants, uniquely up to
conjugation in
$PSL(2,\CC)$, is an immediate consequence of
theorem~\ref{thm:plplanesI}.

\subsection{Relation to Otal's theorem}
 
\label{sec:otal}

In \cite{otal} and later \cite{BonO}, Bonahon and Otal study spaces of various topological types of 3-manifolds with a
hyperbolic structure $\H3/G$ such that $\bch(G)$ is a pleated surface with (in our terminology) a fixed rational pleating lamination.  Translated to our situation, this means the study of a rational pleating plane $\Pr{\gamma,\gamma'} $ for fixed $\gamma,\gamma' \in {\cal S}$.
Write $pl^{+}= \theta \delta_{\gamma}$, $pl^{-}= \theta'
\delta_{\gamma'}$, $ \theta, \theta' \in \RR$.  A special case of their results shows that the map $\Theta(q) =
(\theta(q),\theta'(q))$ is a
homeomorphism from $\Pr{\gamma,\gamma'} $ to an open
neighborhood of
$(0,0)$ in $(0,\pi) \times (0,\pi) $.

Our methods prove that the map $\Theta$ is open and proper; we
have thus far however, been unable to derive injectivity by our methods.
(For the special case $i(\gamma,\gamma')=1$, see \cite{PS}, theorem~3.6.)

Note however that if $q_n \in \Pr{\gamma,\gamma'}$, $q_n \to p
\in \F$,
then $\Theta(q_n) \to (0,0)$ so the whole
critical
line $\F_{\gamma,\gamma'}$ appears on the boundary of this Bonahon-Otal
embedding as a single point.

\section{BM-slices}
\label{sec:BMslices}

In this section we study what happens when we fix the pleating
invariants
on one side of $\bch$. The slices thus defined turn out to be
the complex
extensions of the earthquake paths into $\QF$.

The space of marked conformal structures on $\torus$ can be
identified with
the space $\F$.  For $q \in \QF$, let $w^{\pm}(q)$ denote the
marked
conformal structures of $\Omega^{\pm}/G(q)$.  Bers used the
embedding $q
\mapsto (w^+,w^-)$ of $\QF$ into $ \F \times \overline{\F}$ to
find
holomorphic coordinates for $\F$ by fixing the second factor
$w^-$ and
proving that $w^+$ varies over $\F$; this is called the {\em
Bers
embedding} of $\F$. (Recall that the orientation and hence the
marking on
$\Omega^{-}/G(q)$ is reversed; this is why in the second factor
we write
$\overline{\F}$.)  Maskit, on the other hand, fixed a curve
$\gamma$ on
$\torus$ and studied the family of groups on $\partial\QF$ for
which
$\lambda_{\gamma}=0$ and the corresponding element $V \in G$ is
an
accidental parabolic.  These groups are known as {\em cusps}.
The
conformal structure $w^-$ is then fixed and represents a family
of thrice
punctured spheres; Maskit proved that the first coordinate $w^+$
varies so
as to define an embedding of $\F$ into $\CC$.  We studied the
pleating
invariants for this Maskit embedding of $\F$ in detail in
\cite{KStop}.
McMullen \cite{McMeq}, defines coordinates for Bers embeddings
of $\QF$
that extend to Maskit and generalized Maskit embeddings on
$\partial\QF$.
On the Maskit embeddings his coordinates agree with the pleating
invariants
of \cite{KStop}.

In terms of Minsky's ending invariants \cite{MinskyPT}, both
constructions
correspond to holding the ending invariant of one side fixed and
allowing
the other to vary.  It is thus natural to ask what happens when,
instead of
fixing an ending invariant, we fix the pleating invariants of
one side.

Let $\mu \in ML$, $\c \in \RR^+$ and set $$BM^+_{\mu,\c} = \{ q
\in
{\Pr{\mu}}^+ : \lambda_{\mu}(q) = \c \}.$$

 On $BM^+_{\mu,\c}$, neither the conformal structure on
$\Omega^+/G$ nor
the flat structure on $\bch^+/G$ are fixed.  They are, however,
constrained
by the condition $\lambda_{\mu}(q) = \c $.  We define
$$J:BM^+_{\mu,\c} \rightarrow (PML - \{[\mu]\})\times \RR^+,$$
by
$$J(q) = \Bigl( [pl^-(q)], \frac{l_{pl^-}(q)}{i(\mu,pl^-(q))}
\Bigr).$$

Since $[pl^-(q)]\neq [pl^+(q)]$, $i(\mu,pl^-(q)) > 0$.  The map
$J$ is
continuous by theorem~\ref{thm:plcont}.  Since for fixed $\mu
\in ML$, the
functions $l_{\nu}$ and $i(\mu,\nu)$ scale in the same way as we
vary $\nu$
in its projective class in $PML$, the entry in the second
coordinate of $J$
depends only on $[pl^-]$; it can therefore be written in terms
of our
pleating invariants as $\lambda_{\nu}(q)/i(\mu,\nu)$ for any
choice of $\nu
\in [pl^-]$.

Set
$$ {\cal X}(\mu,\c) = \Bigl\{ ([\nu],s) \in (PML-\{[\mu]\})
\times \RR^+ :
0 < s < \frac{f_{\mu,\nu}(c)}{i(\mu,\nu)} \Bigr\}.$$

 Identifying $PML-\{[\mu]\}$ with $\RR $ as in
section~\ref{sec:laminations}, we can think of ${\cal
X}(\mu,\c)$ as the
region in $\RR \times \RR^+$ under the graph of the function
$[\nu] \mapsto
\frac{f_{\mu,\nu}(c)}{i(\mu,\nu)}$.  As discussed above, this
function is
well defined and by corollary~\ref{cor:X}, it is continuous.

As before, we let $\Q_{\mu,\c}$ denote the quakebend plane along
$\mu$ that
 meets $\F$ along $\E_{\mu,\c}$.  Clearly
$\Q_{\mu,\c}=\Q^p_{\mu}$ for all
 $p \in \E_{\mu,\c}$.

\begin{genericem}{Theorem~\ref{thm:qbplanesI}}
Let $\mu \in ML$ and let $\c > 0$.  Then the closures in $\QF$
of precisely
two of the connected components of $\Q_{\mu,\c} \cap (\QF-\F)$
meet $\F$.
These components are the slices $BM^{\pm}_{\mu,\c}$.  The
intersection of
the closure of each slice with $\F$ is the earthquake path
$\E_{\mu,\c}$;
furthermore each slice is simply connected and retracts onto
$\E_{\mu,\c}$
and the map $J: BM^{\pm}_{\mu,\c} \to {\cal X}(\mu,\c)$ is a
homeomorphism.
\end{genericem}

\begin{proof}
  Noting that for $\nu \in ML$, the pleating ray
$\Pr{\mu,\nu,\c}$ depends
only on the projective class $[\nu]$ of $\nu$, it is clear from
the
definitions that
$$BM^+_{\mu,\c} = \bigcup_{ [\nu] \in
PML-\{[\mu]\}}{\Pr{\mu,\nu,\c}}.$$
Since for $[\nu] \in PML-\{[\mu]\}$, the closure of the pleating
ray
$\Pr{\mu,\nu,\c}$ in $\QF$ contains the point $p(\mu,\nu,\c)$,
the closure
of $BM^+_{\mu,\c}$ in $\QF$ contains $\E_{\mu,\c}$.  It follows
easily from
theorems~\ref{thm:thmq} and~\ref{thm:open} that $BM^+_{\mu,\c}$
is open and
closed in $ \Q_{\mu,\c} \cap (\QF -\F)$.  By
theorem~\ref{thm:bendingaway}
there are no other components of $\Q_{\mu,\c}$ whose closure
meets $\F$.

 For $[\nu] \in PML-\{[\mu]\}$, by lemma~\ref{lemma:image},
$\lambda_{\nu}|_{\Pr{\mu,\nu,\c}}$ is a homeomorphism to the
interval
$(0,f_{\mu,\nu}(c))$. This proves $J$ is a homeomorphism onto
$B^+_{\mu,\c}$.  Clearly therefore, $BM^+_{\mu,\c}$ is simply
connected and
retracts to $\E_{\mu,\c}$ along rays.
\end{proof}

In analogy with  theorem~\ref{thm:ratlplanesdenseI} we have
\begin{thm}
\label{thm:density2}
 The rational pleating rays $\Pr{\mu,\nu,\c}$ are dense in
  $BM^+_{\mu,\c}$.
\end{thm}

\begin{remark} {\rm  As discussed above, holding the Minsky
ending invariant of one side fixed and letting the ending
invariant of the
  other side vary over the full Teichm\"uller space $\F$, we
obtain the
  Bers and Maskit slices.  By contrast, the set of flat
structures $F^-(q)$
  for points $q \in BM^+_{\mu,\c}$ cannot be the full image of
$\F$.  In
  fact, on each ray $\Pr{\mu,\nu,\c}$, the length
$\lambda_{\nu}$ is
  bounded above by $f_{\mu,\nu}(\c)$.  Since by a theorem of
Sullivan,
  \cite{EpM}, lengths on $\bch^-$ and $\Omega^-$ are in bounded
ratio,
  those points on the earthquake path $\E_{\mu,\c}$ in $\F$ at
which
  $\lambda_{\nu}$ is very large will not occur as $F^-(q)$ for
points $q
  \in BM^+_{\mu,\c}$.  See also~\cite{McMeq} for related
phenomena.}
\end{remark}

\section{Rational pleating planes and computation}
\label{sec:ratlplplanes}

We can now easily prove
theorem~\ref{thm:bdryI} of the introduction.

\begin{genericem}{Theorem~\ref{thm:bdryI}}
  Let $\delta_{\gamma},\delta_{\gamma'} $ be rational
laminations
represented by non-conjugate elements $V,V' \in G$.  Then
$\Pr{\gamma,\gamma'}$ and $\Pr{\gamma',\gamma}$ are the unique
components
of the $\RR^2$-locus of the function $\tr V \times \tr V'$ in
$\QF-\F$
whose closures meet $\F$ in $\F_{\gamma,\gamma'}$.  On
$\Pr{\gamma,\gamma'}
\cup \Pr{\gamma',\gamma}$ the function $\tr V \times \tr V'$ is
non-singular and the boundary of $\Pr{\gamma,\gamma'} \cup
\Pr{\gamma',\gamma}$ can be computed by solving $\tr V=\pm 2$
and $\tr V'=
\pm 2$ on this component.
\end{genericem}

\begin{proof}
If $V,V' \in G$ represent $\gamma,\gamma'$ in ${\cal S}$, then
the
$\RR^+$-loci in $\QF$ of $\tr V, \tr V'$ and
$\lambda_{\gamma},\lambda_{\gamma'}$ agree.  As a consequence of
theorem~\ref{thm:plplanesI}, $\Pr{\gamma,\gamma'}$ can be
uniquely
identified as the component of the $\RR^+\times\RR^+$-locus of
$\tr V
\times \tr V'$ which meets $\F$ in the critical line
$\F_{\gamma,\gamma'}$.
\end{proof}

As a consequence of this theorem, given any embedding $ \QF \to
\CC^2$, we
 can compute the position of $\Pr{\gamma,\gamma'}$ and its
boundary
 exactly, provided we can express $\tr V$ and $\tr W$ as
holomorphic
 functions of the parameters and identify the critical line.

For the complex Fenchel Nielsen embedding this works as
follows.  We first note:

\begin{prop}
\label{prop:polynomial} Let  $(\lambda_V,\tau_{V,W})$
be complex Fenchel Nielsen coordinates for $\QF$ relative to a
marked pair
of generators $(V,W)$.  Suppose $\gamma' \in \S$ with
corresponding element
$ V' \in G$.  Then for fixed $\lambda_V$, the trace $\tr V'=\pm
2\cosh
\lambda_{\gamma'}$ is a polynomial in $\cosh{\tau_{V,W}/2}$ and
$\sinh{\tau_{V,W}/2}$.
\end{prop}

\begin{proof}
 From equation~(\ref{eqn:bending2}) we have
   $$ \cosh{\frac{\lambda_{W}}{2}} = {\cosh{\frac{\tau_{V,W}
}{2}}}/
   {\tanh{\frac{\lambda_{V}}{2}}}$$
and
 $$ \cosh{\frac{\lambda_{VW^{\pm 1}}}{2}} =
 {\cosh{\frac{\tau_{V,W} \pm \lambda_{V} }{2}}}/
   {\tanh{\frac{\lambda_{V}}{2}}}.$$

Expanding $\cosh{\frac{\tau_{V,W} \pm \lambda_{V} }{2}}$, the
result
follows in the special cases $V' =W$ and $V' =VW^{\pm1}$.  The
results for
general $V'$ follow from the recursive scheme in \cite{Wright},
see
also~\cite{KStop}, which allows us to express $\tr V'$ as a
polynomial
(with integer coefficients) in $\tr V, \tr W$ and either $ \tr
VW$ or $\tr
VW^{-1}$.
\end{proof}

To find the critical line  $\F_{\gamma,\gamma'}$ we proceed as
follows.
Fix $\c > 0$ and consider the function $\tr V'=\tr
V'(\lambda_V,\tau_{V,W})$.  Along the earthquake path
$\E_{\gamma,\c}$, $
t=\tau_{V,W}$ is real and varies over all of $\RR$; $\lambda_V$
is fixed
and equal to $\c$.  By Kerckhoff's theorem, the function
$\lambda_{\gamma'}$ has a unique critical point $p=
p(\gamma,\gamma',\c)
\in \F_{\gamma,\gamma'}$ along $\E_{\gamma,\c}$; clearly the
same is true
of the trace function $\tr V'$.  Using
proposition~\ref{prop:polynomial},
the position of this point can be computed as a function of $t$.
 Moreover
there are exactly two branches $\sigma^{\pm}$ of the $\RR$-locus
of $\tr
V'$ in $\QF -\F$ whose closures meet $\F$ at $p$.

 By theorem~\ref{thm:plplanesI}, the pleating plane
$\Pr{\gamma,\gamma'}$
is the union of the pleating rays $\Pr{\gamma,\gamma',\c}, \, \c
\in
\RR^+$.  By theorem~\ref{thm:plraysI}, the pleating ray
$\Pr{\gamma,\gamma',\c}$ is one of the two branches
$\sigma^{\pm}$, each of
which maps homeomorphically to $(0, 2\cosh
{{f_{\gamma,\gamma'}(c)} /{2}} )
$ under $\tr V'$.  Analytically continue $\tr V'$ along
$\sigma^{\pm}$.
Again by theorem~\ref{thm:plraysI}, these branches are
non-singular
$\RR$-loci and remain in $\QF$ until they reach points $\tau^*$
such that
$\tr V'(\tau^*)=\pm 2$.  The groups corresponding to such
$\tau^*$ are cusp
groups on $\partial\QF$ for which $\gamma'$ is pinched and $V'$
is an
accidental parabolic.

Drawing these rays for various $\c$'s, we get a picture of the
pleating
planes $\Pr{\gamma,\gamma'}$ and $\Pr{\gamma',\gamma}$.
Allowing $\gamma'$
to vary with $\c$ fixed gives us the slices $BM^{\pm}_{\gamma,
\c}$.  By
theorems~\ref{thm:ratlplanesdenseI} and~\ref{thm:density2}, we
can build up
an arbitrarily accurate picture of $\QF$.  Pictures of various
slices drawn
this way have been obtained in \cite{Wright} and \cite{PP}.

In~\cite{KomS}, similar ideas are used to draw a picture of the
{\em Earle
slice} of $\QF$.  This slice is an embedding of the
Teichm\"uller space of
$\torus$ into $\QF$ consisting of groups for which the
structures on
$\Omega^+$ and $\Omega^-$ are related by a conformal involution
which
induces the rhombus symmetry on $\pi_1(\torus)$.

\subsection{Examples}

We give two examples in which it is especially easy to
compute the pleating plane.

\medskip

\noindent {\em Example 1.}  Take ${\gamma,\gamma'}$ to be
generators of
$\pi_1(\torus)$, represented by the marked pair $V,W \in G$.  By
equation~(\ref{eqn:bending}), $\cosh( \lambda_W/2) = \cosh
(\tau_{V,W}/2)/\tanh (\lambda_V/2)$, so that on the earthquake
path
$\E_{\gamma,\c}$, $\cosh( \lambda_W/2) = \cosh(t/2)/\tanh (\c)$,
$t \in
\RR$.  This function clearly has a unique critical point at the
rectangular
torus $t=0$.  Therefore the critical line  $\F_{\gamma,\gamma'}$
is defined
by the equation $\sinh (\lambda_V/2) \sinh (\lambda_W/2 )= 1$
and the range
of $\lambda_V \times \lambda_W$ is the region
$$ \{ (\c,t) \in \RR^+ \times \RR^+: 0 < \c < 2 \sinh^{-1}
\Bigl(\frac
{1}{\sinh(\c/2)}\Bigr)\}.$$

Notice that under the rectangular symmetry $(V,W ) \to
(V,W^{-1})$ the
group is fixed but the marking is changed; clearly $\Omega^{+}(
G(V,W)) =
\Omega^{-} (G(V,W^{-1}))$.  Thus $\Pr{\gamma,\gamma'}$ maps
bijectively to
$\Pr{\gamma',\gamma}$ while $\F{\gamma,\gamma'}=
\F{\gamma',\gamma}$ is
fixed.  This implies $\cosh\frac{\lambda_{VW}}{2} = \cosh
\frac{\lambda_{VW^{-1}}}{2}$ on $\F{\gamma,\gamma'}$.  Solving
this
equation in $\F$ gives another way of finding the equation of
the critical
line.

\medskip

 \noindent  {\em Example 2. }
Let  $(V,W)$ be a marked pair of
generators for $G$ and let $\gamma,\gamma'$ be the
curves represented by $VW$ and $VW^{-1}$.
Since $G$ is a punctured torus group, the condition that the
commutator  $[V,W]$ be parabolic is expressed by the well known
Markov
equation

\begin{equation}
\label{eqn:Markov}
 \textstyle{ \tr^2 V + \tr^2 W +
\tr^2 VW = \tr V \tr W \tr VW}.
  \end{equation}

Writing $x=\tr V, y=\tr W$, we can solve for $z=\tr VW$ and
$z'=\tr
VW^{-1}$.  On the pleating plane $\Pr{\gamma,\gamma'}$, both $z$
and $z'$
are real so that $xy $ and $ x^2+y^2 $ are real.  It follows
that $x =
\bar{y}$.  Further, on $\Pr{\gamma,\gamma'}$, $x,y \in \RR$ if
and only if
$ G \in \F$.  Thus in the real $(z,z')$ plane, the critical line

$\F_{\gamma,\gamma'}$ has equation $ zz'= 2(z+z')$; in other
words the
hyperbola $(z -2)(z'-2)=4$.  Rewriting in terms of the lengths $
2
\cosh^{-1} \frac{z}{2},2 \cosh^{-1} \frac{z'}{2}$ we find the
region
$T_{\gamma,\gamma'} $ is of the shape claimed.

 We note that in this case, the critical line
$\F_{\gamma,\gamma'}$ is the
fixed line of the rhombic symmetry $(A,B) \to (B,A)$ in $\F$,
giving an
alternative proof that on this line, $\lambda_A = \lambda _B$.
It is also
interesting to note in this example that the Earle slice studied
in~\cite{KomS} is the holomorphic extension of the critical line
$\F_{\gamma,\gamma'}$ into $\QF$.

\section{Appendix}

\subsection{The convergence lemma}
\label{app:convergence}

For the proof of the convergence lemma~\ref{lemma:convergence},
we need to
recall some general facts about laminations.  Let $\Sigma$ be a
hyperbolic
surface and let $\alpha$ be a geodesic lamination on $\Sigma$.
We call a
set $R \subset \Sigma$ a {\em flow box} for $\alpha$ if:
\begin{enumerate}
   \item $R$ is a closed hyperbolic rectangle embedded in
$\Sigma$, with
one pair of opposite sides called ``horizontal'' and the other
pair ``
vertical''.  \item The horizontal sides $T,T'$ of $R$ are either
disjoint
from $\alpha$ or transversal to $\alpha$.  If a leaf $\gamma$ of
$\alpha$
intersects $R$ then it intersects both $T$ and $T'$ .
 \item The vertical
sides of $R$ are disjoint from $\alpha$.
\end{enumerate}

Label the sides of $R$ in counterclockwise order $1,2,3,4$ so
that $1,3$
are the horizontal sides and $2,4$ are the vertical ones.
Suppose that
$\beta \in ML$ is any measured lamination on $\Sigma$. The
underlying
lamination $|\beta|$ intersects $R$ in a family of pairwise
disjoint arcs.
If such an arc joins a vertical to a horizontal side, we call it
a corner
arc; if it joins the two horizontal sides we call it a vertical
arc and
otherwise it is a horizontal arc.  For $i,j \in \{ 1, \ldots,
4\}$, let
$\beta (i,j)= \beta(j,i)$ denote the total transverse measure of
the arcs
joining side $i $ to side $j$.  Clearly, $\alpha (1,3) = \alpha
(3,1) =
\alpha(T) = \alpha(T')$, the transverse measure of the
transversal $T$,
while $\alpha (i,j) =0$ otherwise.

The following simple lemma  applies to any hyperbolic surface
$\Sigma$.

\begin{genericem}{Lemma 11.1}
\label{lemma:thinbox} Let $\nu_0 \in ML$ and let $R$ be a flow
box for $|\nu_0|$. Suppose $\nu_0(T) \neq 0$.  Then for $\nu \in
ML$
sufficiently near $\nu_0$, the lamination $|\nu|$ has a vertical
arc.
\end{genericem}

\begin{proof}
Note that because $|\nu|$ consists of pairwise disjoint simple
geodesics,
it does not have both horizontal and vertical arcs.  Let $V,V'$
denote the
vertical sides.  Since $\nu_0(V)=\nu_0(V')=0$, both $\nu(V)$ and
$\nu(V')$
can be assumed arbitrarily small by taking $\nu$ sufficiently
close to
$\nu_0$ in $ML$.  We can write $\nu(V) = \nu(4,1)+
\nu(4,2)+\nu(4,3)$ and
$\nu(V') = \nu(2,1)+ \nu(2,4)+\nu(2,3)$.  All the terms on the
right in
these relations are non-negative so each is arbitrarily small.

If we assume $|\nu|$ has no vertical arc we have $\nu(T) =
\nu(4,1)+
\nu(2,1)$, $\nu(T') = \nu(3,2)+ \nu(3,4)$ and by the above we
deduce that
both are arbitrarily small.  But this is a contradiction because
$\nu(T)$
and $\nu(T')$ are both near $\nu_0(T)$ which is a definite
positive value.
\end{proof}

Now we need some facts specific to laminations on a punctured
torus (see
\cite{ThuN}, 9.5.2).  Let $\alpha \in {\cal S}$ and cut $\torus$
along
$\alpha$ to obtain a punctured annulus $A$ with boundary curves
$\alpha_1$
and $\alpha_2$.  The leaves of any measured lamination $\nu$,
$|\nu| \neq
\alpha$ intersect $A$ in a union of arcs that either join
$\alpha_1$ to
$\alpha_2$ or join one of the boundary components to itself.  It
is easy to
show, (see~\cite{ThuN}), that the set of arcs joining a
component
$\alpha_i$ to itself has zero transverse measure.  In
particular, by
minimality any transversal to any leaf of $|\nu|$ carries
non-zero measure,
so that all arcs of $|\nu|$ in $A$ join $\alpha_1$ to
$\alpha_2$.

We also recall that on $\torus$, if $\nu \not\in ML_Q$, the
complement of
$|\nu|$ is a punctured bigon $B$, and also that there is a
horocyclic
neighborhood of definite size about the cusp disjoint from the
support of
any measured lamination.

Now we can prove the convergence lemma~\ref{lemma:convergence}.

\begin{genericem}{Lemma~\ref{lemma:convergence}}
\label{lemma:appconvergence}  Suppose that $\nu_0 \in ML -ML_Q$,
and that $\nu $ and $\nu_0$ are close in $ML$. Then $|\nu|$ and
$|\nu_0|$
are close in the Hausdorff topology on $GL$.
\end{genericem}

\begin{proof}
First we show that given a long arc in $|\nu_0|$ there exists a
long nearby
arc in $|\nu|$.  Let $L, \epsilon>0$ be given.  Since $\nu_0 \in
ML - ML_Q$,
all leaves have infinite length.  Thus, given $x \in |\nu_0|$,
by choosing
sufficiently short transversals we can find a flow box for which
the leaf
of $|\nu_0|$ through $x$ is a vertical arc, the segments of
length $L$ on
either side of $x$ are contained in $R$, and the horizontal
sides of $R$
have length less than $\epsilon$.  We call a flow box of this
kind, a good
$\epsilon, L$-flow box for $x$.  Now standard hyperbolic
geometry estimates
show, that if two geodesics are a bounded distance apart over a
long
distance $t$, then in fact they are close to order $ e^{-t}$
along a large
fraction of their length.  Thus any vertical arc in a good
$\epsilon,L$-flow box is certainly close to leaves of $|\nu_0|$
over
distance at least $2L$.  Clearly, $|\nu_0|$ can be covered by a
finite
number of flow boxes of this kind.

Now suppose we are given a long arc $\lambda$ of a leaf of
$|\nu_0|$.  Let
$x$ be the midpoint of $\lambda$ and let $R$ be a good
$\epsilon,L$-flow
box for $x$.  By lemma~\ref{lemma:thinbox}, we deduce that if
$\nu \in ML$
is near $\nu_0$, then $\nu$ has a vertical arc in $R$ so that by
the above,
$|\nu|$ has long arc of a leaf near $\lambda$ as required.

Next we claim conversely, that given a long arc in $|\nu|$ there
exists a
long nearby arc in $|\nu_0|$.  For a lamination $\lambda$, let
$T_1(\lambda)$ denote the the set of unit tangent vectors to
leaves
pointing along leaves of $\lambda$.  Since there is a horocyclic
neighborhood of definite size about the cusp disjoint from the
support of
any measured lamination on $\torus$, the set $\cup_{\lambda \in
GL}
T_1(\lambda)$ is a compact subset of the unit tangent bundle
$T_1(\torus)$.
Clearly, laminations $\lambda$ and $\lambda'$ are close in the
Hausdorff
topology on closed subsets of $GL$ if and only if $T_1(\lambda)$
and
$T_1(\lambda')$ are close in the Hausdorff topology on closed
subsets of
$T_1(\torus)$.

If our claim is false, then there is a sequence of points $\bar
v_n \in
 T_1(|\nu_n|)$, $\nu_n \in ML_Q$ with $\nu_n \rightarrow \nu_0$
in $ML$,
 for which there are no nearby points of $T_1(|\nu_0|)$.  A
geodesic
 $\beta$ through a limit point of the vectors $\bar v_n$ will be
a limit of
 leaves of $|\nu_n|$, but will not be a leaf of $|\nu_0|$.

 If $\beta \cap |\nu_0| \neq \emptyset$, we obtain a
contradiction.  For if
$x \in \beta \cap |\nu_0|$, the tangent directions to $\beta$
and $|\nu_0|$
at $x$ are distinct.  Therefore we can find a good $|\nu_0|$
flow box $R$
for $x$, such that the arc of $\beta$ through $x$ is only close
to the leaf
of $\nu_0$ through $x$ for a short distance and thus cannot be
either a
vertical or a corner arc in $R$.  But then all laminations
$|\nu|$ with
leaves close to $\beta$ also contain arcs which must intersect
$R$ in
horizontal arcs, contradicting lemma~\ref{lemma:thinbox}.

To complete the proof we must show $\beta \cap |\nu_0| \neq
\emptyset$.  If
not, then $\beta$ is contained in the complement of $|\nu_0|$ in
$\torus$.
Since $\nu_0 \not\in ML_Q$, the complement of $|\nu_0|$ is a
punctured
bigon $B$.  If $\beta$ enters $B$ through one vertex and leaves
through the
other it is homotopic to, and therefore coincides with, a leaf
of
$|\nu_0|$; thus $\beta$ must come in from one vertex of the
bigon, go
around the puncture and return back to the {\em same} vertex.
Let $\alpha$
be a simple closed curve that intersects $\beta$ and as above,
cut $\torus$
along $\alpha$ to obtain a punctured annulus $A$ with two
boundary curves
$\alpha_1, \alpha_2$.  Since $\beta$ goes around the puncture,
it crosses
one of the $\alpha_i$ and returns through the same side of
$\alpha_i$ (see
the figure in \cite{ThuN}, 9.5.2).  It follows that any closed
simple
geodesic sufficiently close in the Hausdorff topology to $\beta$
would also
have an arc entering and leaving $A$ across the same $\alpha_i$.
But any
arc of a simple closed geodesic carries a non-zero transverse
measure, and
by the fact stated above, must join $\alpha_1$ to $\alpha_2$.
Hence $\beta
\cap |\nu_0| \neq \emptyset$.
\end{proof}

\subsection{Proof of proposition~\ref{prop:P2}.}

Before beginning the proof, we need to review the definitions of
the
bending measure and intrinsic metric for paths on $\bch$ as
given in
\cite{KSconvex}.  We suppose that $ q \in \QF$, and that as
usual $\bch=
\bch(q)$ is the convex hull boundary of $\HH^3/G(q)$.  We shall
only
indicate the dependence on $q$ when needed in the proof.  In
fact, we shall
only need to apply what follows to the component $\bch^+$.

A {\em support plane} for $\bch$ at a point $x \in \bch$ is a
hyperbolic
plane $P$ containing $x$ such that $\C$ is contained entirely in
one of the
two half spaces cut out by $P$.  The {\em bending angle} between
two
intersecting support planes $P_1,P_2$ at points $x_1,x_2 \in
\bch$ is the
absolute value of the angle $\theta(P_1,P_2)$ between their
outward normals
from $\bch$.

Let $\Pi(x)$ denote the set of oriented support planes at $x \in
\bch$ and
let
$$Z= \{(x, P(x)) : x \in \bch, P(x) \in \Pi(x) \}, $$ with
topology induced
from ${\cal G} = \H3 \times {\cal G}_2(\H3)$, where ${\cal
G}_2(\H3)$ is
the Grassmanian of 2-planes in $\H3$. Let $Z^+$ be the obvious
restriction
of $Z$ to $\bch^+$ and call the it {\em approximating set} for
$\bch^+$.

To define the bending measure and intrinsic metric, it suffices
to define
the measure and length of any path $\bar \omega$ on $\bch$.  Any
such path
lifts to a path $\omega \colon [0,1] \rightarrow Z$ as follows.
Suppose $
x \in \bar \omega$.  Either $\Pi(x)$ consists of a unique point,
in which
case there is nothing to do, or we add to the path an arc in
which the
first coordinate $x$ is fixed but the second moves continuously
on the line
in $\cal G$ from the left to the right extreme support planes at
$x$.

A {\em polygonal approximation} to $\omega$ is a sequence
$$\P=\{\omega(t_i)= (x_{i},P_i) \in Z\} ; 0=t_0 < t_1 <
\ldots < t_n=1,$$
such that $P_i \cap P_{i+1} \neq \emptyset, i=0,\ldots,n-1$.

Let $\theta_i=\theta(P_{i-1},P_i)$ be the bending angle between
$P_{i-1}$
and $P_{i}$, $i=1,\ldots,n$ and let $d_i$ be the hyperbolic
length of the
shortest path from $x_{i-1}$ to $x_i$ in the planes $P_{i-1}
\cup P_i$.

The intrinsic metric on $\bch$ is given by

\begin{equation}
\label{eqn:length}
 l(\omega) = \inf_{\P} \sum_{i=1}^n d_i
\end{equation}
and the bending measure $\beta$ on $\bch$  by

\begin{equation}
\label{eqn:bendingmeasure}
\beta(\omega)= \inf_{\P} \sum_{i=1}^n \theta_i
\end{equation}
where $\P$ runs over all
polygonal approximations to $\omega$.

\medskip

In order to prove proposition~\ref{prop:P2}, we shall also make
similar
polygonal approximations to the pleated surface
$\PP_{\gamma}(q)$.  We
shall prove the proposition by showing that polygonal
approximations in
$Z^+= Z^+(q)$ to the convex hull boundary $\bch^+$ can be
replaced by
polygonal approximations to the pleated surface
$\PP_{\gamma}(q)$, and that
the above approximating sums are simultaneously good
approximations to the
intrinsic metric of the flat structure $\phi_{\gamma}^+(q)$ and
the
transverse measure $b_{\gamma}(q)$.  Thus we also need to
discuss polygonal
approximations for $\PP_{\gamma}(q)$.

 \medskip

The surface $\PP_{\gamma}(q)$ is made up of planar pieces,
precisely two of
which meet along each bending line $\alpha$ (which projects to $
\gamma$ on
$\torus$.  Call a plane $P$ a {\em pseudo-support plane} to
$\PP_{\gamma}(q)$ if either it is one of these planar pieces, or
if it
meets $\PP_{\gamma}(q)$ along $\alpha $ and lies in the half
space cut out
by the planar pieces of $\PP_{\gamma}(q)$ through $\alpha$.  The
pseudo-support planes of $\PP_{\gamma}(q)$ inherit natural
orientations
from the pleated surface map under which $\PP_{\gamma}(q)$ is an
immersed
image of the hyperbolic disk $\DD$ in $\HH^3$.

Let $\tilde \Pi(x)$ denote the set of oriented pseudo-support
planes at $x
\in \PP_{\gamma}(q)$ and let
$$W=W(q)= \{(x, P(x)) | x \in \PP_{\gamma}(q), P(x) \in \tilde
\Pi(x) \},
$$ with topology induced from $\cal G$ as before.  We define
polygonal
approximations in $W(q)$ in the obvious way, and call $W(q)$ the
approximating set for $\PP_{\gamma}(q)$.

We claim that the flat metric $\phi_{\gamma}(q)$ and the measure
$b_{\gamma}$ on $\PP_{\gamma}(q)$ are defined by sums similar to
those
in~(\ref{eqn:length}) and~(\ref{eqn:bendingmeasure}), where the
infimum is
taken now over polygonal approximations in $W(q)$.

Let $\omega$ be a path in $W$ and let $\{(x_i,Q_i)\}$ be such a
$W$-polygonal approximation.  As in the proof of Proposition~4.8
of~\cite{KSconvex}, we consider the segment of path $\omega_i$
in $W(q)$
between $x_{i-1}$ and $x_i$, and we work in a hyperbolic plane
$H$ through
$x_{i-1}$ and $x_i$, such that the shortest path $\sigma$ from
$x_{i-1}$ to
$x_i$ in the planes $Q_{i-1} \cup Q_i$ is contained in the
intersections of
these planes with $H$.  Let the segments of $\sigma$ in
$Q_{i-1}$ and $Q_i$
have lengths $a_1$ and $a_2$ respectively, so that $a_1+a_2$ is
an upper
bound for the contribution to the sum giving the length of
$\omega_i$.
Notice that even though we do not know that $\PP_{\gamma}(q)$
bounds a
convex half space, it follows easily from Gauss-Bonnet that
$\omega_i$ does
not intersect $\sigma$.  Thus it is easy to check that inserting
an extra
pair $(x,Q) \in W$ between $x_{i-1}$ and $x_i$, the
approximating sum for
the length of $\omega_i$ decreases.  Since by assumption
$[\gamma] \in
ML_Q$, there are in fact sufficiently fine polygonal
approximations for
which the sum in~(\ref{eqn:length}) actually {\em equals} the
intrinsic
metric on $\PP_{\gamma}(q)$. A similar argument, on the lines of
that in
Proposition~4.8 of~\cite{KSconvex}, shows that the
sums~(\ref{eqn:bendingmeasure}) decrease on inserting extra
support planes
and that there are sufficiently fine sums which actually equal
the measure
$b_{\gamma}$.

We are now ready to prove proposition~\ref{prop:P2}.

\label{app:deferred2}
\begin{genericem}{Proposition~\ref{prop:P2}}
 \label{prop:P20} Given $\nu_0 \in ML - ML_Q$, and $q_0 \in
\Pr{\nu_0}^+\cup\F$, let $F^+(q_0)\in \F$ and $pl^+(q_0)$ be the
flat
structure and bending measure of $\bch^+(q_0)$ respectively. (If
$q_0 \in
\F$, then $pl^+(q_0)=0$ and $F^+(q_0)$ is the Fuchsian group
representing
$q_0$.)  Then, given neighborhoods $V$ of $F^+(q_0)$ in $\F$ and
$W$ of
$pl^+(q_0)$ in $ML$, there exist neighborhoods $U$ of $q_0$ in
$\QF$ and
$X$ of $[\nu_0]$ in $PML$ such that if $q \in U$,
$[\delta_{\gamma}] \in X
\cap PML_Q$ and $\lambda_{\gamma}(q) \in \RR^+$, then the flat
structure
$\phi_{\gamma}(q)$ of $\PP_{\gamma}(q)$ is in $V$ and the
transverse
measure $b_{\gamma}(q)$ is in $W$.
\end{genericem}

\begin{proof}
Let $\nu_0,q_0$ be as in the statement of the proposition.
Suppose that
for some $q$ near $q_0$ and $ [ \delta_{\gamma}]$ near
$[\nu_0]$, we have
$\lambda_{\gamma}(q) \in \RR^+$.  Let $\PP_{\gamma}(q)$ be the
associated
pleated surface with approximating set $W(q) \subset {\cal G}$
as above.
Let $Z^+(q_0)$ and $Z^+(q)$ be the approximating sets for
$\bch^+(q_0),
\bch^+(q)$ respectively.

We claim that for every $(x,P(x)) \in Z^+(q_0)$ and $q \in \QF$
near
$q_0$, there is a nearby pair $(y,P(y)) \in W(q)$, and
conversely.  This
will follow immediately if we can show that, for every geodesic
in
$|pl^+(q)|$, there is a geodesic in the bending locus of
$\PP_{\gamma}(q)$
with nearby endpoints in $\HH^3$, and vice versa.  Now, the
crucial
condition $\lambda_{\gamma}(q) \in \RR^+$ implies that the
bending locus of
$\PP_{\gamma}(q)$ is {\em exactly} $\gamma=\gamma(q) $.  Thus,
applying
lemma~\ref{lemma:convergence} to the laminations $ \nu_0$ and
$\k
\delta_{\gamma}$ for a suitable choice of $ \k>0$ on the surface
$\bch^+(q_0)$, we see that $|\nu_0(q_0)|$ and $ \gamma(q_0)$ are
close in
the Hausdorff topology on closed subsets of $\bch^+(q_0)$.
Lifting to
$\HH^3$, this means that the endpoints $x_0,x_0'$ of any lift of
a leaf of
$|\nu_0(q_0)|$ are close to the endpoints $x,x'$ of a lift of
$\gamma(q_0)$
and vice versa since the geodesic representative of $\gamma$ on
$\bch^+$
has the same endpoints as the geodesic $\gamma$ in $\HH^3$.  It
follows
that the $\HH^3$ geodesics with $x_0,x_0'$ and $x,x'$ also have
long close
arcs.  Finally, moving to a nearby point $q$ in $\QF$, the
endpoints of
geodesics which project to the leaves of $|\nu_0(q)|$ are close
to the
endpoints of geodesics which project to $|\nu_0(q_0)|$, and
similarly for
endpoints of geodesics which project to $\gamma(q)$ and
$\gamma(q_0)$.  The
claim follows.

We now consider the key estimates which were the basis of the
continuity
results proved in~\cite{KSconvex}.  Call a polygonal
approximation an {\em
$(\alpha,s)$-approximation} if $$\max_{1 \leq i \leq n}
\theta(P_{i-1},P_i)
< \alpha$$ and $$\max d_{\omega} (x_{i-1},x_i) < s$$ where
$d_{\omega}$ is
distance along $\omega$ measured in the intrinsic metric on
$\bch$.  We
have

\begin{genericem}{Proposition \cite{KSconvex}, Prop.~4.8}
\label{prop:error estimate}
There is a universal constant $K$, and a function $s(\alpha)$
with values
in $(0,1)$, such that if $\P$ is an
$(\alpha,s(\alpha))$-approximation to a
path $\omega$ in $Z$, where $ \alpha < \pi/2$, then
 $$ | \sum_{\P} d_i - l(\omega)| < K \alpha l(\omega) $$ and
$$ | \sum_{\P} \theta_i - \beta(\omega)| < K \alpha l(\omega)
.$$
\end{genericem}

To complete the present proof, it suffices to check that similar
estimates
hold if polygonal approximations in $Z^+(q_0)$ are replaced by
approximations in $W(q)$.  The estimates work in exactly the
same way; the
only point to note is that we need the same local convexity
property
implied by Gauss Bonnet as above.
\end{proof}

\subsection{Proof of proposition~\ref{prop:genbendinglimit}}
\label{app:deferred3}

\begin{genericem}{Proposition~\ref{prop:genbendinglimit}}
\label{prop:genbendinglimit0}
 Suppose $\mu \in ML$, $q \in \Pr{\mu} \cup \F$ and consider the
quakebend
plane $\Q^q_{\mu}$ along $\mu$ based at $q$ with parameter
$\tau_{\mu}$.
Given $K > 0$, there exists $B > 0$ such that if $|\Re
\tau_{\mu}| <K$ and
$|\Im \tau_{\mu}| > B$, then $\Q^q_{\mu}(\tau_{\mu}) \notin
\Pr{\mu}$.
\end{genericem}

\begin{proof}
Our proof will show that if $\tau_{\mu}$ is outside
the range described the proposition, then the pleated surface
obtained by bending by $\tau$  along $\mu$
cannot be embedded and thus that
$\Q^q_{\mu}(\tau_{\mu}) \notin \Pr{\mu}$.
The group $\Q^q_{\mu}(\tau_{\mu})$ may or may not be in $\QF$.

We use the definitions of  support planes  and bending angles
from the proof of proposition~\ref{app:deferred2}.
 From the definition, the bending angle between two intersecting
support planes $P_1,P_2$ to $\bch$ at points $x_1,x_2$
 is an upper bound for the bending measure of a transversal
to $|\mu|$ joining $x_1,x_2$ which lies
 between the ``roof'' formed by
$P_1$ and $P_2$ and the $\HH^3$ geodesic from $x_1$ to $x_2$.

We make the following claims.
\begin{enumerate}
\item There exists $\epsilon >0$ such that if $x_1,x_2,x_3 \in
\bch$ lie in
a ball of radius $\epsilon$ in $\BB^3$, and if $P_1,P_2,P_3$ are
support
planes to $\bch$ at $x_1,x_2,x_3$ respectively, then either $P_1
\cap P_3
\neq \emptyset$, or both $P_1 \cap P_2\neq \emptyset$ and $P_2
\cap P_3\neq
\emptyset$.

\item Given $\epsilon >0$, $\mu \in ML$, $\mu \ne 0$, and a
compact subset
$V \subset \cal F$, there is a constant $a>0$ such that if $\phi
\in V$,
then there is a transversal $\kappa$ to $|\mu|$ with hyperbolic
length
$l(\kappa) <\epsilon$ in the structure $\phi$ and transverse
measure
$\mu(\kappa) > a$.
\end{enumerate}

\noindent{\em Proof of claim 1.}  A support plane $P$ to $\bch$
meets
 $\Chat$ in a circle which contains points of the limit set
$\Lambda$ and
 which bounds a disk $D(P)$ containing no points of $\Lambda$.
Therefore
 if $P_1 \cap P_3 = \emptyset$, the discs $D(P_1)$ and $D(P_3)$
are
 disjoint.  To prove the claim amounts to showing that in this
case, both
 $D(P_1) \cap D(P_2)$ and $D(P_3) \cap D(P_2)$ are non-empty.
Without loss
 of generality, we may suppose that $x_1,x_2,x_3$ are within
hyperbolic
 distance $\epsilon$ of the origin $O$ in $\BB^3$ so that the
planes $P_i$
 are close to equatorial planes through $O$.  The result is then
obvious.

\smallskip
\noindent{\em Proof of claim 2.}  Choose $\gamma \in \S$ with
$i(\gamma,\mu) > 0$.  There are constants $c_1,c_2,d_1,d_2$ such
that $c_1<
\mu (\gamma) < c_2$ and $d_1< l(\gamma) < d_2$ for $\phi \in
V$. Subdividing $\gamma$ into $N$ segments with $d_2 /N <
\epsilon$, the
result is clear with $a= c_1 /N$.

\smallskip
Now, working in the quakebend plane $\Q^q_{\mu}$, with parameter
$\tau=
\tau_{\mu}$, consider the set of groups for which $|\Re \tau
|<K$.  The
corresponding flat structures $F^+(\tau)$ are independent of
$\Im \tau$ and
thus lie in a compact set $V \subset \F$.  Choose a transversal
$\kappa$ as
in claim (2). Let $x_1,x_3$ be its initial and final points and
$x_2$ its
midpoint, and let $P_i$ be a support plane at $x_i$.  Using
claim 1, either
$P_1,P_3$, or both pairs $P_1,P_2$ and $P_2,P_3$, intersect.
Thus at least
one of the segments $(x_1,x_3)$, $(x_1,x_2)$ or $(x_2,x_3)$ of
$\kappa$,
for definiteness say the segment $\kappa_1$ joining $(x_1,x_2)$,
has
$\mu(\kappa_1) >a/2$.

Consider the point in $\Q^q_{\mu}$ with parameter $\tau$. The
bending
measure $pl^{+}(\tau)(\kappa_1)$ of $\kappa_1$ is $k+\Im \tau
\mu(\kappa_1)$, where $k=pl^{+}(q)(\kappa_1)$ is the bending
measure of
$\kappa_1$ at the base point $q$.  The bending angle between
$P_1,P_2$ is
bounded above by $\pi$.  As in the first paragraph, this gives
an upper
bound for $pl^{+}(\tau)(\kappa_1)$, and we obtain the required
bound on
$|\Im \tau |$.
\end{proof}

\end{document}